\documentclass[11pt]{article}
\usepackage[letterpaper]{geometry}
\usepackage{amsfonts,amssymb,amsmath,dsfont}
\usepackage{graphicx,color,mathrsfs}
\usepackage{times,latexsym}
\usepackage{fullpage}
\usepackage{enumerate}
\usepackage[toc,page]{appendix}

\numberwithin{equation}{section}

\newtheorem{theorem}{Theorem}

\newtheorem{corollary}[theorem]{Corollary}

\newcommand{\R}{\mathbb{R}}

\newcommand{\E}{\mathbb{E}}

\newcommand{\cL}{\mathcal{L}}

\newcommand\1{\mathds{1}}
\newcommand{\indi}[1]{\1_{\{#1\}}}

\begin{document}

\begin{center}
\Large \bf Bounded-Velocity Stochastic Control for Dynamic Resource Allocation%
  \footnote{An abbreviated preliminary form of this paper has appeared in a conference proceedings (without copyright transfer) \cite{GaLuSh+13}.}
\end{center}
%\funding{Funding information goes here.}}}
\author{}
\begin{center}
{Xuefeng Gao}\,%
  \footnote{Department of Systems Engineering and Engineering Management, The Chinese University of Hong Kong, Shatin, Hong Kong; (xfgao@se.cuhk.edu.hk).},%
  {Yingdong Lu}\,%
  \footnote{Mathematical Sciences Department, IBM Research AI~-–~Science, Yorktown Heights, NY 10598, USA; (yingdong@us.ibm.com).},
 { Mayank Sharma}\,%
  \footnote{AI and Blockchain Solutions Department, IBM Research, Yorktown Heights, NY 10598, USA; (mxsharma@us.ibm.com).},
  {Mark S.~Squillante}\,%
   \footnote{Mathematical Sciences Department, IBM Research AI~-–~Science, Yorktown Heights, NY 10598, USA; (mss@us.ibm.com).},
  {Joost W.~Bosman}\,%
  \footnote{Centrum Wiskunde \& Informatica, 1098 XG Amsterdam, The Netherlands; (J.W.Bosman@cwi.nl)}
\end{center}

%\begin{center}
%{Xuefeng
%  Gao}\,\footnote{Corresponding Author. Department of Systems
%    Engineering and Engineering Management, The Chinese University of Hong Kong, Shatin, N.T. Hong Kong;
%    xfgao@se.cuhk.edu.hk},
%  Lingjiong Zhu\,\footnote{Department of Mathematics, Florida State University, 1017 Academic Way, Tallahassee, FL-32306, United States of America; zhu@math.fsu.edu.
%  }
%\end{center}

% REQUIRED
\begin{abstract}
We consider a general class of dynamic resource allocation problems within
a stochastic optimal control framework.
This class of problems arises in a wide variety of applications, each of
which intrinsically involves resources of different types and demand with
uncertainty and/or variability.
The goal involves dynamically allocating capacity for every resource type
in order to serve the uncertain/variable demand, modeled as Brownian motion,
and maximize the discounted expected net-benefit over an infinite time horizon
based on the rewards and costs associated with the different resource types,
subject to flexibility constraints on the rate of change of each type of resource capacity.
We derive the optimal control policy within a bounded-velocity stochastic
control setting, which includes efficient and easily implementable
algorithms for governing the dynamic adjustments to resource allocation
capacities over time.
Computational experiments investigate various issues of both theoretical and
practical interest, quantifying the benefits of our approach over recent
alternative optimization approaches.
\end{abstract}

% REQUIRED
%\begin{keywords}
%  example, \LaTeX
%\end{keywords}
%
%% REQUIRED
%\begin{AMS}
%  68Q25, 68R10, 68U05
%\end{AMS}

\section{Introduction}\label{sec:intro}
A general class of canonical forms of dynamic resource allocation problems arises naturally
across a broad spectrum of computer systems, communication networks and business applications.
As
%the complexities of these systems, networks and applications
their complexities
continue to grow, together
with ubiquitous advances in technology, new approaches and methods are required to effectively
and efficiently solve canonical forms of general dynamic resource allocation problems in such
complex system, network and application environments.
These environments often consist of different types of resources that are allocated in combination
to serve demand whose behavior over time includes diverse types of uncertainty and variability.
Each type of resource has a different reward and cost structure that ranges from the best of a
set of primary resource allocation options~--~having the highest reward, highest cost and highest
net-benefit~--~to the next best primary resource allocation option~--~having the next highest
reward, next highest cost and next highest net-benefit~--~and so on down to a secondary resource
allocation option~--~having the lowest reward, lowest cost and lowest net-benefit.
Each type of resource also has different degrees of flexibility and different cost structures
with respect to making changes to the allocation capacity.
The resource allocation optimization problem we consider consists of adaptively determining the
vector of primary resource capacities and the secondary resource capacity that serve the
uncertain/variable demand and that maximize the expected net-benefit over a time horizon of
interest based on the foregoing structural properties of the different types of resources.

Motivated by this general class of dynamic resource allocation problems, we take a stochastic control
approach that manages future risks associated with resource allocation
decisions and uncertain demand, including the reward, cost and flexibility structures of the
primary and secondary resource allocation options.
Specifically, we consider the underlying fundamental stochastic optimal control problem where
the dynamic control policy that allocates the set of primary resource capacities to serve
uncertain/variable demand, modeled as Brownian motion, is a vector of absolutely continuous
stochastic processes with constraints on its element-wise rates of change with respect to time
(``bounded-velocity''), which in turn determines the secondary resource allocation capacity.
The ultimate objective is to maximize the expected discounted net-benefit over an infinite horizon
based on the structural properties of the different resources types, with the desired outcome
rendering an explicit characterization of the optimal dynamic control policy that includes
efficient and easily implementable algorithms for governing the dynamic adjustments to resource
allocation capacities over time.

\subsection{Motivating applications}\label{sec:examples}
The wide variety of system, network and application domains in which arise the general class of
canonical forms of dynamic resource allocation problems of interest in this study include cloud
computing and data center environments, computer and communication networks, energy-aware computing
and smart power grid environments, and business analytics and optimization, among many others.
For example, large-scale cloud computing and data center environments often involve resource
allocation over different server options (ranging from the fastest performance and most expensive
to the slowest performance and least expensive) and different network bandwidth options (ranging
from the highest guaranteed performance at the highest cost to opportunistic options at no cost,
such as the Internet); e.g., refer to~\cite{AmCiGu+10,ChHeLi+08,ClFrHa+05,JaMeMo+12,PaSpTa+05}.

In related energy-aware computing environments, the dynamic resource allocation problem concerns
the effective and efficient management of the consumption of energy by resources in the face of
time-varying uncertain system demand and (especially in large-scale environments) energy prices;
e.g., see~\cite{GuJaWi11,LiWiAn+11}.
In particular, the system control policy dynamically adjusts the allocation of very high-performance,
very high-power servers (best primary resource option), the next highest performance, next highest
power servers (next primary resource option), and so on to satisfy the uncertain/variable system
demand where any remaining demand is satisfied by low-performance, low-power servers (secondary
resource), with the objective of maximizing expected (discounted) net-benefit over time.
Here, the rewards (costs) are based on the performance (energy) properties of the type of servers
allocated to satisfy demand over time, together with the additional per-unit costs incurred for
respectively increasing or decreasing the allocation of the primary resource capacities over time.
The resulting optimal dynamic control policy determines the adaptive control of the primary resource
capacities at every time, subject to certain constraints and to the demand process.

Another motivating application is based on strategic business provisioning and workforce sourcing in
human capital supply chains; refer to, e.g.,~\cite{Wagner-2011}.
The demand for product and services offerings is satisfied through resource allocation from among a diversity
of available supply options.
These sourcing options include internal resources with the highest reward, highest cost and highest
net-benefit, business-partner resources with the next highest reward, next highest cost and next
highest net-benefit, and so forth down to contractor or crowdsourcing resources with the lowest reward,
lowest cost and lowest net-benefit.
For each supply option, examples of the rewards can be revenue and quality of service and examples
of the costs can be salary and other compensation.
The goal of the strategic human capital sourcing problem is to determine the portfolio of various
supply options to meet the time-varying and uncertain demand in the sense of maximizing expected
(discounted) net-benefit over time, including the costs incurred in hiring, reskilling, promoting
and incentivizing to reduce attrition for the primary human capital resources over time.
This also involves constraints on the rates of change in the capacities of primary options to accommodate
the non-instantaneous adjustment of human capital resources and the limited availability of such resources
in the labor market (see, e.g.,~\cite{DaDeSe+87}).

%Related issues arise in smart power grids where resource allocation is required across a diversity
%of available energy sources, including a long-term market (stable and less expensive, but rather
%inflexible), local generation (significant operating constraints and limited capacity) and a real-time
%spot market (readily available and responsive, but at a premium price), together with the additional
%uncertainty from renewables such as wind and solar; e.g., refer to~\cite{GhKaKa+11,VaWuBi11}.
%Across these and many other domain-specific resource allocation problems, there is a common need
%for the dynamic adjustment of allocations among multiple types of resources, each with different
%structural properties, to satisfy time-varying and uncertain demand.

\subsection{Related work}
The research literature contains a great diversity of studies of resource allocation problems,
with differing objective functions, control policies, and reward, cost and flexibility structures.
A wide variety of approaches and methods have been developed and applied to (approximately) solve
this diversity of resource allocation problems including, for example, online algorithms and
dynamic programming.
It is therefore important to compare and contrast our problem formulation and solution approach
with some prominent and closely related alternatives.
One classical instance of a dynamic resource allocation problem is the multi-armed bandit
problem~\cite{mab2007} where the rewards are associated with tasks and the goal is to determine
under uncertainty which tasks the resource should work on, rather than the other way around.
Another widely studied problem is the \emph{ski-rental} or \emph{lease-or-buy}
problem~\cite{linial1999} where there is demand for a resource, but it is initially not known as
to how long the resource would be required.
In each decision epoch, the choice is between two options: either lease the resource for a fee,
or purchase the resource for a price much higher than the leasing fee.
Our resource allocation problem differs from this situation in that there are multiple types of
resources each with an associated reward and cost per unit time of allocation, since the
resources cannot be purchased outright.

From a methodological perspective, the general resource allocation problem we consider in this
paper is closely related to the vast literature on stochastic control; refer to, e.g.,~\cite{Pham09,YonZho99}.
Of particular relevance is the so-called \emph{bounded velocity
follower problem}~\cite{BeShWi80,KarOco02} where the control is an absolutely continuous process with bounded derivative.
For example, Bene{\v{s}} et al.~\cite{BeShWi80} consider this {bounded velocity
follower problem} with a quadratic running cost objective function, where the authors propose a
smooth-fit principle to characterize the optimal policy.
In comparison with our study, the paper only considers a single resource, does not consider any
costs associated with the actions taken by the control policy, and deals with a smoother objective
function.
Karatzas and Ocone~\cite{KarOco02} study an alternative bounded velocity follower problem where
discretionary stopping is allowed and the objective is to choose a control law and a stopping
time that minimizes the expected sum of a running cost and a termination cost.
Davis et al.~\cite{DaDeSe+87} consider a control problem where a decision maker determines the
rate of investment in capacity expansion to satisfy a Poisson demand process, showing that the
optimal control is ``bang-bang'' type (either carry out construction at the maximal speed or take no
action) and using computational methods to compute the value function.
In contrast, we derive the optimal solution to a stochastic optimal
control problem that does not allow discretionary stopping and that seeks to maximize the expected
discounted net-benefit over an infinite horizon under a Brownian motion demand process, where we
analytically characterize the corresponding value function and explicitly characterize the optimal
dynamic control policy.
%Taking an alternative variational inequality approach that does not use the smooth-fit principle,
%Menaldi and Robin~\cite{MenRob83} study multidimensional stochastic control problems,
%including the one-dimensional bounded velocity follower as a special case, and show that the value
%function is continuous.
Our work is also related to, but different from, the problem of \emph{drift control} for Brownian motion
(see, e.g., \cite{Ata2005,Matoglu2011,Matoglu2015}) where the controller can, at some cost, shift the drift among a finite set of alternatives, and
the problem of \emph{bounded variation control} for diffusion processes (see, e.g., \cite{Weerasinghe2005})
where the available control is an added bounded variation process (but without bounded velocity constraints).
Additional related studies include reversible/irreversible investment and capacity planning using a stochastic control approach,
for which we refer the interested reader to \cite{abel1996,bentolila1990firing,federico2012smooth,guo2005optimal,guoandtomecek09,MerhiandZervos07,shepp1996hiring}
as well as \cite{van2003commissioned} for a survey and a list of references.

From an applications perspective, there is a growing interest and vast literature in the computer system,
communication network and operations research communities to address allocation problems involving various
types of resources associated with computation, memory, bandwidth and/or energy;
refer to, e.g., \cite{CioFar12,Dieker2012,JaMeMo+12,Kasbekar2014,LiWiAn+11,Urgaonkar2010} and the references therein.
We limit our discussions here to the works of Lin et al.~\cite{LiWiAn+11} and Ciocan and Farias~\cite{CioFar12}
since we will compare our dynamic control policy through computational experiments with the optimization algorithms
from both of these studies that are cast within discrete time-interval models.
Lin et al.~\cite{LiWiAn+11} study the problem of dynamically adjusting the number of
active servers in a data center as a function of demand to minimize operating costs.
They consider average demand over small intervals of time, subject to system constraints, and
propose an optimal offline algorithm together with an online algorithm that is shown to be within a constant factor
worse than the proposed optimal offline policy.
Ciocan and Farias~\cite{CioFar12} study model predictive control for a large class of dynamic resource allocation
problems with stochastic demand rate processes.
They develop
%a simple
an online algorithm for demand allocation within appropriately selected time intervals that relies on frequent reoptimization using
suitably updated demand forecasts.
Our work differs from these two and the aforementioned studies in that we take a stochastic control approach
for the dynamic resource allocation optimization problem where the control has bounded velocity constraints and associated costs.
%Jagannathan et al.~\cite{JaMeMo+12} consider the decision-making process of ``secondary users'' in
%dynamic spectrum access communication systems who can choose between acquiring dedicated spectrum
%(which provides immediate yet costly transmission) or sharing spectrum-holes of ``primary users''
%(which provides free yet delayed and unreliable transmission), where the authors characterize the
%resulting Nash equilibrium behavior for the single-band case.

\subsection{Our contributions}
The goal of our study is to determine the optimal dynamic risk hedging strategy for managing the
portfolio of primary resource allocation options and secondary resource allocation option in order to
maximize the expected discounted net-benefit over an infinite time horizon based on the structural
properties of the different primary and secondary resources, which we show to be equivalent to a
minimization problem involving a piecewise-linear running cost and a proportional cost for making
adjustments to the control policy process.
Our solution approach is based on explicitly constructing a twice continuously differentiable (with the exception of at most three points) solution to the corresponding Hamilton-Jacobi-Bellman
equation.
Our theoretical results also include an explicit characterization of the dynamic control policy,
which is of threshold type,
and then we verify that this control policy is optimal through a martingale argument.
In contrast to an optimal static allocation strategy, in which a single primary allocation vector capacity
is determined to maximize expected net-benefit over the entire time horizon, our theoretical results
establish that the optimal dynamic control policy adjusts its allocation decisions in primary and
secondary resources to hedge against the risks of under-allocating primary resource capacities
(resulting in lost reward opportunities) and over-allocating primary resource capacities (resulting
in incurred cost penalties).
%We present extensive computational experiments that quantify the effectiveness of our optimal dynamic
%control policy, including comparisons of how our control policy convincingly outperforms a
%standard offline optimal policy which has often been used in previous studies.

Our study provides important methodological contributions and new theoretical results by deriving
the solution of a fundamental stochastic optimal control problem.
This stochastic optimal control solution approach highlights the importance of timely and adaptive
decision making in the allocation of a mixture of different resource options with distinct features
in optimal proportions to satisfy time-varying and uncertain demand.
Our study also provides important algorithmic contributions through a new class of online policies for
dynamic resource allocation problems arising across a wide variety of application domains.
Computational experiments quantify the effectiveness of our optimal online dynamic control
algorithm over recent work in the area, including comparisons demonstrating how our optimal online
algorithm significantly outperforms the type of optimal offline algorithm within a discrete-time framework
recently proposed in~\cite{LiWiAn+11}, which appears to be related to the optimal online model predictive
control algorithm proposed in~\cite{CioFar12} within a different discrete-time stochastic optimization framework.
This includes relative improvements up to $90\%$ and $130\%$ in comparison with the optimal offline
algorithm considered in~\cite{LiWiAn+11}, and even larger relative improvements of more than $150\%$
and $230\%$ in comparison with the optimal online algorithm in~\cite{CioFar12}.

An abbreviated preliminary form of this paper has appeared in a conference proceedings (without copyright transfer) \cite{GaLuSh+13} in which we presented
a subset of our results on the optimal control policy for resource allocation problems with only a single primary resource option and a secondary resource option,
together with limited computational experiments demonstrating some of the benefits of our optimal dynamic control policy over the optimal offline
algorithm proposed in~\cite{LiWiAn+11}.
The current paper significantly extends the preliminary conference paper in four important aspects:
First, we present herein our complete derivation of the optimal control policy for resource allocation problems with single primary and secondary resource options;
Second, we present herein a more thorough set of computational experiments that demonstrate the significant benefits of our optimal dynamic control policy over both
the optimal offline algorithm in~\cite{LiWiAn+11} and the optimal online algorithm proposed in~\cite{CioFar12};
Third, we extend herein our theoretical results beyond the case of a single primary resource allocation option to a general case with multiple primary resource allocation options
and provide computational experiments that compare the single and multiple primary resource allocation models;
Finally, we provide herein rigorous proofs of all our main results for both resource allocation models under all possible conditions on the model parameters.

The remainder of this paper is organized as follows.
Section~\ref{sec:math} defines our mathematical models of the stochastic resource allocation optimization problem,
first for the case of a single primary resource and then for a version of the multiple primary resources case.
Our mathematical formulations and main results for the corresponding stochastic control problems are presented
in Sections~\ref{sec:main} and \ref{sec:main2}, respectively.
A representative sample of numerous computational experiments are discussed in Section~\ref{sec:computational}.
%Section~\ref{sec:main2} discusses the case of multiple primary resources, followed by some concluding remarks.
%Section~\ref{sec:main2} discusses the use of model predictive control to determine the parameters of our optimal dynamic control policy over time, followed by some concluding remarks.
Concluding remarks are provided in Section~\ref{sec:conclusion},
including a discussion of the use of model predictive control and learning to determine the parameters of our optimal dynamic control policy over time.
All of our technical proofs are provided in the appendices.
%The appendices provide additional proofs and supporting theoretical material.

\section{Mathematical Models}\label{sec:math}
We investigate a general class of fundamental resource allocation problems in
which a set of \emph{primary} resource allocation options and a \emph{secondary}
resource allocation option are available to satisfy demand whose behavior over
time includes uncertainty and/or variability.
There is an important reward and cost ordering among these resource options
where the first primary resource allocation option has the highest net-benefit,
followed by the second primary resource allocation option having the next highest
net-benefit, and so on, with the (single) secondary resource allocation option
having the lowest net-benefit.
In addition, the set of primary resource capacities are somewhat less flexible
in the sense that their rates of change at any instant of time are bounded,
whereas the secondary resource capacity is more flexible in this regard (as made
more precise below).
Beyond these differences, the set of primary resources and the secondary resource
are capable of serving the demand and all of this demand needs to be served, i.e.,
there is no loss of demand.

To elucidate the exposition of our analysis, we first consider the single primary
and secondary instance of our general resource allocation model.
This instance captures the key aspects of the fundamental trade-offs among the
net-benefits of the various resource allocation options together with their associated risks.
We then consider our general resource allocation model comprising multiple primary resources
and a secondary resource together with an important relationship maintained among the primary resources.
%The system models for each case are next defined in turn.

\subsection{System model: Single primary resource}\label{sec:model:Single}
We consider the stochastic optimal control problems underlying our
resource allocation model in which uncertain and/or variable demand needs to
be served by the primary resource allocation capacity and the secondary
resource allocation capacity.
A control policy defines at every time $t \in \R_+$ the level of capacity allocation for the primary resource, denoted by $P(t)$,
and the level of secondary resource capacity allocation, denoted by $S(t)$, that
are used in combination to satisfy the uncertain/variable demand, denoted by $D(t)$.

The demand process $\{D(t): t \ge 0 \}$ is modeled by
\begin{eqnarray}
dD(t) &=& b dt +\sigma d{W}(t),
\label{eqn:demand}
\end{eqnarray}
where $b \in \R$ is the demand growth/decline rate (which can be extended to a deterministic function of time),
%but we do not consider this further in the present paper),
$\sigma>0$ is the demand volatility, and ${W}(t)$ is a one-dimensional
standard Brownian motion, whose sample paths are continuous but
nondifferentiable almost everywhere~\cite{KarShr91}.
This model is natural when the demand forecast involves Gaussian noises; see, e.g., \cite{ChoMeyn2010,Yao2010}
for similar Brownian demand models.
The demand process is served by the combination of primary and secondary resource
allocation capacities $P(t) + S(t)$.
Given the higher net-benefit structure of the primary resource option, the
optimal dynamic control policy seeks to determine an absolutely continuous stochastic process $P(\cdot)$ describing the
primary resource allocation capacity
to serve the demand $D(\cdot)$ such that any remaining demand is served by the
secondary resource allocation capacity $S(\cdot)$.
% All the parameters can become policy dependent in the generalized model.

Let $R_p(t)$ and $C_p(t)$ respectively denote the reward and cost
associated with the primary resource allocation capacity $P(t)$ at time $t$.
The rewards $R_p(t)$ are linear functions of the primary resource
capacity and the demand,
whereas the costs $C_p(t)$ are linear functions of the primary
resource capacity.
We therefore have
\begin{equation}
R_p(t) \; = \; \mathcal{R}_p\times [P(t) \wedge D(t)]
%& \qquad \mbox{and} \qquad &
\qquad\qquad \mbox{and} \qquad\qquad
C_p(t) \; = \; \mathcal{C}_p\times P(t) , \quad \label{eq:benefitP:Single}
\end{equation}
where $x\wedge y := \min\{x, y\}$, $\mathcal{R}_p \geq 0$ captures all
per-unit rewards for serving demand with the primary resource capacity,
$\mathcal{C}_p \geq 0$ captures all per-unit costs for the primary
resource capacity, and $\mathcal{R}_p > \mathcal{C}_p$.
Observe that the rewards for the primary resource are linear in $P(t)$
as long as $P(t) \leq D(t)$, otherwise any
primary resource capacity exceeding the demand will solely incur
costs without rendering rewards.
Hence, from a risk hedging perspective, the risks associated with the
primary resource allocation position at time $t$, $P(t)$, concern lost
reward opportunities whenever $P(t) < D(t)$ on the one hand
and concern incurred cost penalties whenever $P(t) > D(t)$
on the other hand.

In addition, any adjustments to the primary resource allocation capacities
have associated costs, where we write $\mathcal{I}_p$ and $\mathcal{D}_p$
to denote the per-unit costs of increasing and decreasing the decision process
$P(t)$, respectively.
Namely, $\mathcal{I}_p$ represents the per-unit cost whenever the allocation
of the primary resource capacity is being increased while $\mathcal{D}_p$
represents the per-unit cost whenever the allocation of the primary resource
capacity is being decreased.

Since all the remaining demand is served by the
secondary resource allocation capacity, we therefore have
\begin{equation}
S(t) \; = \; [D(t) - P(t)]^+ .
\label{eq:secondary:Single}
\end{equation}
The corresponding reward function $R_s(t)$ and cost function $C_s(t)$ are
then given by
\begin{equation}
R_s(t) \; = \; \mathcal{R}_s\times [D(t) - P(t)]^+
%& \qquad\qquad \mbox{and} \qquad\qquad &
\qquad\qquad \mbox{and} \qquad\qquad
C_s(t) \; = \; \mathcal{C}_s\times [D(t) - P(t)]^+ , \quad
\label{eq:benefitS:Single}
\end{equation}
where $x^+ := \max\{x,0\}$, $\mathcal{R}_s \geq 0$ captures all per-unit
rewards for serving demand with the secondary resource capacity,
$\mathcal{C}_s \geq 0$ captures all per-unit costs for the secondary resource
capacity, and $\mathcal{R}_s > \mathcal{C}_s$.
%Modified Cost Function: $C_s(t) = c_2([D(t) -P(t)]^+)$, here, we extend our model to include availability, $c_2(\cdot)$ now becomes a convex function.
Hence, from a risk hedging perspective, the secondary resource
allocation position at time $t$, $S(t)$, is riskless in the sense
that rewards and costs are both linear in the resource capacity
actually used.

\subsection{System model: Multiple primary resources}\label{sec:model:Multiple}
We next consider the stochastic optimal control problems underlying our
resource allocation models in which uncertain and/or variable demand needs to
be served by the set of primary resource allocation capacities and the secondary
resource allocation capacity.
Let $\mathcal{P}$ denote the number of primary resource options.
A control
policy then defines at every time $t \in \R_+$ the level of capacity allocation for all
primary resources, respectively denoted by $P_1(t), P_2(t), \ldots, P_{\mathcal{P}}(t)$,
and the level of secondary resource capacity allocation, denoted by $S(t)$, that
are used in combination to satisfy the uncertain/variable demand, denoted by $D(t)$.

The demand process $\{D(t): t \ge 0 \}$ continues to be modeled as in \eqref{eqn:demand}.
This demand process is served by the combination of primary and secondary resource
allocation capacities $P_1(t)+P_2(t)+\ldots+P_{\mathcal{P}}(t)+S(t)$, where $P_i(\cdot)$ has absolutely continuous paths for each $i=1, \ldots, \mathcal{P}$.
Given the higher net-benefit structures of the primary resource options, the
optimal dynamic control policy seeks to determine at every time $t \in \R_+$ the
primary resource allocation capacity vector $\mathbf{P}(t) = [P_1(t), \ldots, P_{\mathcal{P}}(t)]$
to serve the demand $D(t)$ such that any remaining demand is served by the
secondary resource allocation capacity $S(t)$.
% All the parameters can become policy dependent in the generalized model.

The corresponding stochastic control problem becomes high-dimensional in the presence of multiple primary resources, and thus it is inherently prohibitive to solve analytically or computationally in general.
% due to the curse of dimensionality.
%One may set up the dynamic programming equation for the stochastic optimization problem, and obtain a second order PDE (partial differential equation) in $\R^{\mathcal{P}}$ with unknown boundary.
%Computationally solving such free boundary PDE in high-dimension efficiently and related convergence analysis is challenging.
%On the other hand, alternative computational approximation for stochastic control problems based on Markov Chain methods~\cite{KusDup01} also suffer the curse of dimensionality.
We therefore introduce in our system model definition the notion of coordination through a contractual agreement among the multiple primary resources,
namely the model definition includes a contract that fixes the ratio of capacities among the primary resources.
From a mathematical perspective, this contract-based system model definition makes it possible to derive explicit solutions for the corresponding stochastic control problem.
In particular, by introducing such a contractual agreement among the multiple primary resource options, the dimensionality of the stochastic control problem is relaxed
and our mathematical derivations exploit and build on the results we obtain for the single primary resource model.
%the optimization problem then becomes how to control the aggregate capacity of primary resources and secondary options to meet the volatile demand, and determine the optimal fraction of demand should be served by each primary
%resource option.
From an applications perspective, our contract-based model definition is consistent with coordinated contractual agreements that have been adopted in various resource allocation problems in practice,
such as strategic sourcing in human capital supply chains (see, e.g., \cite{bulmash2010} and the strategic sourcing example in the introduction)
because they capture important business relationships and are easy to implement.

%Assuming the ratio of employees within each
%level of the workforce is fixed
%
%\textcolor{red}{Assuming the ratio of primary resource capacities is fixed (e.g., assume the same emp within each
%level of the workforce is fixed; Assuming the same employee distribution from 2009; }

More formally, we introduce a nonnegative vector $w=(w_1, \ldots, w_{\mathcal{P}}) \in \R_+^{\mathcal{P}}$, with $\sum_{i=1}^{\mathcal{P}} w_i =1$,
to represent the contract-based fixed distribution of capacities among the multiple primary resources.
Then, for each $t \ge 0$, we set
\begin{eqnarray} \label{eq:multip}
P(t) := \sum_{i=1}^{\mathcal{P}} P_i(t), \qquad\qquad \text{and} \qquad\qquad P_i(t) = w_i P(t), \quad \text{for $i=1, \ldots, \mathcal{P}$},
\end{eqnarray}
where $P(t)$ (with a slight abuse of notation) represents the aggregate capacity of primary resources and $w_i$ is set to maintain the initial (agreed upon) percentage $P_i(0)/P(0)$.
In other words, for a given contract vector $w$, all $\mathcal{P}$ primary resource allocations must maintain the relationship \eqref{eq:multip} at all time $t$.

For each primary resource type $i=1, \ldots, \mathcal{P}$,
let $R_{p,i}(t)$ and $C_{p,i}(t)$ respectively denote the reward and cost associated with the $i$th primary resource allocation capacity at time $t$.
When the collection of primary resource allocations exceeds the demand, the rewards $R_{p,i}(t)$ are linear functions of the fraction
of the demand served by the $i$th primary resource allocation capacity, namely $w_i \cdot D(t)$, from \eqref{eq:multip}.
However, when the demand exceeds the collection of primary resource allocations,
then the rewards $R_{p,i}(t)$ are linear functions of the $i$th primary resource allocation capacity, which is less than $w_i \cdot D(t)$ due to \eqref{eq:multip}.
The costs $C_{p,i} (t)$ are linear functions of the $i$th primary resource allocation capacity.
We therefore have
\begin{eqnarray}
R_{p,i}(t) \; = \; \mathcal{R}_{p,i}\times [P_i(t) \wedge (w_i \cdot D(t))]
& \qquad \mbox{and} \qquad &
C_{p,i}(t) \; = \; \mathcal{C}_{p,i}\times P_i(t) , \quad \label{eq:benefitP:Multiple}
\end{eqnarray}
where $x\wedge y := \min\{x, y\}$, $\mathcal{R}_{p,i} \geq 0$ captures all
per-unit rewards for serving demand with the $i$th primary resource capacity,
$\mathcal{C}_{p,i} \geq 0$ captures all per-unit costs for the $i$th primary
resource capacity, and $\mathcal{R}_{p,i} > \mathcal{C}_{p,i}$.
The per-unit reward and the per-unit cost for the aggregate primary resource capacity $P(t)$, under a given contract vector $w$, are then respectively given by
\begin{eqnarray} \label{eq:Rpw}
\mathcal{R}_{p}(w) = \sum_{i=1}^{\mathcal{P}} w_i \mathcal{R}_{p,i}, \qquad
\mathcal{C}_{p}(w) = \sum_{i=1}^{\mathcal{P}} w_i \mathcal{C}_{p,i}.
\end{eqnarray}

Observe that the rewards for the $i$th primary resource allocation are linear in $P_i(t)$
as long as $P_i(t) \leq w_i D(t)$, or equivalently $P(t) \leq D(t)$;
otherwise the fraction of the $i$th primary resource capacity exceeding $w_i D(t)$ will
solely incur costs without rendering rewards.
Hence, from a risk hedging perspective, the risks associated with the collection of
primary resource allocations at time $t$ concern lost reward opportunities whenever $P(t) < D(t)$
on the one hand and concern incurred cost penalties whenever $P(t) > D(t)$ on the other hand.

Similarly to \eqref{eq:Rpw},
as any adjustments to the primary resource allocation capacities have associated costs,
let $\mathcal{I}_{p,i}$ and $\mathcal{D}_{p,i}$ respectively denote the per-unit cost associated
with increasing and decreasing the allocation of the $i$th primary resource capacity.
The per-unit costs of increasing and decreasing the allocation of the aggregate primary
resource capacity are then respectively given by
\begin{eqnarray} \label{eq:Ipw}
\mathcal{I}_{p}(w) = \sum_{i=1}^{\mathcal{P}} w_i \mathcal{I}_{p,i}, \qquad \mathcal{D}_{p}(w) = \sum_{i=1}^{\mathcal{P}} w_i \mathcal{D}_{p,i} .
\end{eqnarray}

Since the optimal dynamic control policy serves all remaining demand with
secondary resource allocation capacity, we therefore have
\begin{equation}
S(t) \; = \; [D(t) - \sum_{i=1}^{\mathcal{P}} P_i(t)]^+ .
\label{eq:secondary:Multiple}
\end{equation}
The corresponding reward function $R_s(t)$ and cost function $C_s(t)$ are
then given by
\begin{eqnarray}
R_s(t) \; = \; \mathcal{R}_s\times [D(t) - \sum_{i=1}^{\mathcal{P}} P_i(t)]^+
& \qquad \mbox{and} \qquad &
C_s(t) \; = \; \mathcal{C}_s\times [D(t) - \sum_{i=1}^{\mathcal{P}} P_i(t)]^+ , \quad
\label{eq:benefitS:Multiple}
\end{eqnarray}
where $x^+ := \max\{x,0\}$, $\mathcal{R}_s \geq 0$ captures all per-unit
rewards for serving demand with secondary resource capacity,
$\mathcal{C}_s \geq 0$ captures all per-unit costs for secondary resource
capacity, and $\mathcal{R}_s > \mathcal{C}_s$.
%Modified Cost Function: $C_s(t) = c_2([D(t) -P(t)]^+)$, here, we extend our model to include availability, $c_2(\cdot)$ now becomes a convex function.
Hence, from a risk hedging perspective, the secondary resource
allocation position at time $t$, $S(t)$, is riskless in the sense
that rewards and costs are both linear in the resource capacity
actually used.

\section{Main Results: Single Primary Resource}\label{sec:main}
In this section we consider our main results on the optimal dynamic control policy for the stochastic optimal control problem
when there is a single primary resource.
%, corresponding to the system model of Section~\ref{sec:model:Single}.
%To simplify the presentation, we shall suppress the single primary resource index
%$i=1$ in our notation throughout this section and, as an example, simply write
%$P(t)$ for the primary resource allocation at time $t$.
After providing a formulation of the stochastic control problem and some technical preliminaries,
we present our main results under different conditions for the values of $\mathcal{I}_p \geq 0$ and $\mathcal{D}_p \geq 0$,
including the special case in which both are zero.
%Consideration of the proofs of our main results in these cases is postponed until
%the next section.

\subsection{Problem formulation}\label{sec:formulation:Single}
The stochastic optimal control problem associated with the system model of Section~\ref{sec:model:Single}
allows the dynamic control policy to adjust its allocation positions in
primary and secondary resource capacities based on the demand realization
observed up to the current time,
which we call the risk-hedging position of the dynamic control policy.
More formally, in the single primary resource case, the decision process $P(t)$,
is adapted to the filtration $\mathcal{F}_t$ generated by $\{D(s): s \le t\}$.
Then the objective of the optimal dynamic control policy is to maximize the
expected discounted net-benefit over an infinite horizon, where net-benefit at
time $t$ consists of the difference between the rewards and costs from the
primary resource allocation capacity and the secondary resource allocation
capacity minus the additional costs for adjustments to $P(t)$.

In formulating the corresponding stochastic optimization problem, we impose
a pair of additional conditions on the decision process
$\{P(t): t \ge 0\}$ based on practical aspects of the diverse application
domains motivating our study.
The control policy cannot
%instantaneously change the primary resource
%allocation capacity in an attempt to directly follow the demand
%$D(t)$; i.e., some time is required (even if only
%a very small amount of time) to adjust $P(t)$.
%Moreover, the control policy cannot
make unbounded adjustments in the
primary resource allocation capacity at any instant in time; i.e., the amount
of change in $P(t)$ at time $t$ is restricted (even if only to a very small
extent) by various factors.
We therefore assume that the rate of change in the
primary resource allocation capacity by the control policy is bounded.
More precisely, for an absolutely continuous process $P(\cdot)$, there is a pair of finite constants $\theta_\ell < 0$ and
$\theta_u > 0$ such that \text{for each $t \ge 0$}
\begin{eqnarray} \label{eq:adjbound:Single}
\theta_\ell \; \le \; \dot P(t) \; \le \; \theta_u,
\end{eqnarray}
where $\dot P(t)$ denotes the derivative of the decision variable $P(\cdot)$ (absolutely continuous process)
with respect to time.
On the other hand, the ability of the control policy to make adjustments to
the secondary resource capacity in response to changes in the primary resource
capacity tends to be more flexible such that \eqref{eq:secondary:Single} holds at all time $t$.

Now we can present the mathematical formulation of our stochastic optimization problem
for the case of a single primary resource.
Defining
\begin{eqnarray*}
N_p(t) \; := \; R_p(t) - C_p(t)  & \qquad\qquad \mbox{and} \qquad\qquad & N_s(t) \; := \; R_s(t) - C_s(t) ,
\end{eqnarray*}
we seek to determine the optimal dynamic control policy that solves the problem (SC-OPT:S)
\begin{eqnarray}
\max_{P(\cdot)} \quad && \E \int_0^\infty e^{-\alpha t} \left[ N_p(t) + N_s(t) \right] dt -
\E \int_0^\infty e^{-\alpha t} [ \mathcal{I}_p \cdot \indi{\dot P(t)>0} ] dP(t) \nonumber \\
&& \qquad\qquad\qquad\qquad\qquad\qquad\quad - \; \E \int_0^\infty e^{-\alpha t} [ \mathcal{D}_p \cdot \indi{\dot P(t)<0} ] d(-P(t)) \qquad \label{opt1:obj:Single} \\
\mbox{s.t.} && -\infty \; < \; \theta_\ell \; \le \; \dot P(t) \; \le \; \theta_u \; < \; \infty , \qquad\qquad \mbox{for $t \ge 0$}, \qquad \label{opt1:st1:Single} \\
&& \qquad dD(t) \; = \; b dt +\sigma d{W}(t), \qquad\qquad \mbox{for $t \ge 0$}, \label{opt1:st2:Single}
\end{eqnarray}
where $\alpha$ is the discount factor and $\indi{A}$ denotes the indicator
function returning $1$ if $A$ is true and $0$ otherwise.
The control variable is the rate of change in the primary resource capacity
by the control policy at every time $t$ subject to the lower and upper bound
constraints on $\dot P(t)$ in \eqref{opt1:st1:Single}.
Note that the second (third) expectation in \eqref{opt1:obj:Single} causes a decrease
with rate $\mathcal{I}_p$ ($\mathcal{D}_p$) in the value of the
objective function whenever the control policy increases (decreases) $P(t)$.

\subsection{Preliminaries}\label{sec:main:prelim:Single}
For notational convenience, we define the constants
\begin{eqnarray}
r_1 \; := \; \frac{b+ \sqrt{b^2 +2 \alpha \sigma^2 }}{\sigma^2} \; > \; 0, \label{eq:r1}
& \quad &
r_2 \; := \; \frac{b- \sqrt{b^2 +2 \alpha \sigma^2 }}{\sigma^2} \; < \; 0, \label{eq:r2} \\
s_1 \; := \; \frac{b- \theta_u + \sqrt{(b-\theta_u)^2 +2 \alpha \sigma^2 }} {\sigma^2} \; > \; 0, \label{eq:s1}
& \quad &
s_2 \; := \; \frac{b- \theta_u - \sqrt{(b-\theta_u)^2 +2 \alpha \sigma^2 }} {\sigma^2} \; < \; 0 , \label{eq:s2} \\
t_1 \; := \; \frac{b- \theta_\ell + \sqrt{(b-\theta_\ell)^2 +2 \alpha \sigma^2 }} {\sigma^2} \; > \; 0, \label{eq:t1}
& \quad &
t_2 \; := \; \frac{b- \theta_\ell - \sqrt{(b-\theta_\ell)^2 +2 \alpha \sigma^2 }} {\sigma^2} \; < \; 0 . \label{eq:t2}
\end{eqnarray}
These quantities are the roots of the quadratic equation
\[ \frac{\sigma^2}{2} y^2 + (\theta -b) y - \alpha =0,\]
when $\theta$ takes on the values of $\theta_\ell$, $0$ or $\theta_u$.

%Since $\theta_\ell<0$ and $\theta_u>0$, we conclude that $B_k$ and $J_k$
%are all positive for $k=1,2,3$, and that $A$ and $K$ are both negative.

Next, we turn to consider the first expectation in the objective function
\eqref{opt1:obj:Single} of the stochastic optimization problem (SC-OPT:S),
which can be simplified as follows.
Define
\[
X(t) \; := \; P(t) - D(t) , \qquad\qquad \mathcal{N}_p \; := \; \mathcal{R}_p - \mathcal{C}_p , \qquad\qquad \mathcal{N}_s \; := \; \mathcal{R}_s - \mathcal{C}_s ,
\]
and $x^- := -\min\{x,0\}$.
Upon substituting \eqref{eq:benefitP:Single} and \eqref{eq:benefitS:Single}
into the first expectation in \eqref{opt1:obj:Single},
and making use of the fact that
\begin{eqnarray*}
[P(t) \wedge D(t)] &=& D(t) - [ D(t) - P(t) ]^{+} ,
%{[D(t) - P(t)]}^+ &=& {[P(t) - D(t)]}^- ,
\end{eqnarray*}
we obtain
\begin{equation}
\E \left[ \int_0^\infty e^{-\alpha t} [-\mathcal{C}_p X(t) + (\mathcal{N}_s-\mathcal{R}_p) X(t)^-] dt \right]
+ \mathcal{N}_p \E \left[ \int_0^\infty e^{-\alpha t} D(t) dt \right] . \label{eq:simplify}
\end{equation}
Note that the second expectation in \eqref{eq:simplify} represents the expected discounted cumulative demand over the infinite horizon.
Since this second summand in \eqref{eq:simplify} does not depend on the control variable $P(t)$,
this term plays no role in determining the optimal dynamic control policy.
%Noting that $x = x^+ - x^-$ and that $\max y = \min -y$,
%in combination with the above results,
Together with the above results,
we derive the following stochastic optimization problem which is
equivalent to the original problem formulation (SC-OPT:S):
\begin{eqnarray}
\min_{P(\cdot)} \quad && \E_x \bigg[ \int_0^\infty e^{-\alpha t} \Big\{ \left( \mathcal{C}_+ X(t)^+ + \mathcal{C}_- X(t)^{-} \right) dt
+ \left( \mathcal{I}_p \indi{\dot P(t)>0} - \mathcal{D}_p \indi{\dot P(t)<0} \right) dP(t) \Big\} \bigg] \qquad \label{optHF:obj} \\
\mbox{s.t.} && -\infty \; < \; \theta_\ell \; \le \; \dot P(t) \; \le \; \theta_u \; < \; \infty , \qquad\qquad \mbox{for $t \ge 0$}, \label{optHF:st1} \\
&& dX(t) \; = \; dP(t) - b dt - \sigma d{W}(t),  \qquad\qquad \mbox{for $t \ge 0$}, \label{optHF:st2} \\
&& X(0) \; = \; x , \qquad\qquad \mathcal{C}_+ \; = \; \mathcal{C}_p , \qquad\qquad \mathcal{C}_- \; = \; \mathcal{N}_p - \mathcal{N}_s,  \label{optHF:st3}
\end{eqnarray}
where $\E_x[\cdot]$ denotes expectation with respect to the initial state
distribution (i.e., state at time $t=0$) being $x$ with probability one.

We use $V(x)$ to represent the optimal value of the objective function
\eqref{optHF:obj}; namely, $V(x)$ is the value function of the corresponding
stochastic dynamic program.
Given its equivalence with the original optimization problem (SC-OPT:S), the
remainder of this section will focus on the stochastic dynamic program
formulation in \eqref{optHF:obj}~--~\eqref{optHF:st3}.

Finally, for additional convenience in stating our main results, we further
define the following constants
\begin{eqnarray*}
B_1 := ({\mathcal{C}_+}-\alpha \mathcal{D}_p) (t_2 -r_2), & \qquad &
    B_2 := (\mathcal{C}_- - \alpha \mathcal{I}_p) (s_1 -r_2), \qquad\qquad
    B_3 := (\mathcal{C}_+ + \mathcal{C}_-) (-r_2), \\
J_1 := ({\mathcal{C}_+}-\alpha \mathcal{D}_p) (r_1 -t_2), & \qquad &
    J_2 := (\mathcal{C}_- - \alpha \mathcal{I}_p) (r_1 -s_1), \qquad\qquad
    J_3 := (\mathcal{C}_+ + \mathcal{C}_-) r_1, \\
A := (\mathcal{C}_+ + \alpha \mathcal{I}_p) (r_2-r_1), & \qquad &
    K := (\mathcal{C}_- + \alpha \mathcal{D}_p)(r_2-r_1) ,
\end{eqnarray*}
where $r_1, r_2, s_1, s_2, t_1, t_2$ are given in \eqref{eq:r1}~--~\eqref{eq:t2} and ${\mathcal{C}_+}, {\mathcal{C}_-}$ are given in \eqref{optHF:st3}.

\subsection{Case 1: $\boldsymbol{\mathcal{D}_p<\mathcal{C}_+/\alpha}$ and $\boldsymbol{\mathcal{I}_p<\mathcal{C}_-/\alpha}$}\label{sec:main:case1}
Let us first briefly explain the conditions of this subsection,
which are likely to be the most relevant case in practice.
Observe from the objective function~\eqref{optHF:obj} that
$\mathcal{C}_+/\alpha$ reflects the discounted overage cost associated
with the primary resource capacity and $\mathcal{C}_-/\alpha$ reflects
the corresponding discounted shortage cost, recalling that $\alpha$ is
the discount rate.
In comparison, $\mathcal{D}_p$ represents the cost incurred for decreasing
$P(t)$ when in an overage position while $\mathcal{I}_p$ represents the
cost incurred for increasing $P(t)$ when in a shortage position.

To elucidate the exposition, we denote

{Condition (1a):} \quad $\mathcal{I}_p + \mathcal{D}_p >0$, \quad $0 \; < \; {B_3- B_2} \; < \; {B_1} \quad \mbox{and} \quad \left(\frac{B_3- B_2}{B_1}\right)^{\frac{r_2}{r_1}} \; \ge \; \frac{J_3-J_2}{J_1}$.

{Condition (1b):} \quad $\mathcal{I}_p + \mathcal{D}_p >0$, \quad \mbox{and} \quad $B_3 \le B_2$.

{Condition (1c):} \quad $\mathcal{I}_p = \mathcal{D}_p =0$, \quad \mbox{and} \quad $B_3 -B_2 -B_1 \le 0$.

{Condition (2a):} \quad $\mathcal{I}_p + \mathcal{D}_p >0$, \quad ${B_3-B_2- B_1} \; > \; 0 \qquad \mbox{and} \qquad \left(\frac{B_3-B_1}{B_2}\right)^{\frac{r_2}{r_1}} \;  \ge \; \frac{J_3 -J_1}{J_2}$.

{Condition (2b):} \quad $\mathcal{I}_p = \mathcal{D}_p =0$, \quad \mbox{and} \quad $B_3 -B_2 -B_1 \ge 0$.

\noindent
We are now ready to state our main result for Case 1.

\begin{theorem}\label{THM:CASE1}
Suppose the adjustment costs satisfy $\mathcal{D}_p<\mathcal{C}_+/\alpha$
and $\mathcal{I}_p<\mathcal{C}_-/\alpha$.
Then there are two threshold values $L$ and $U$ with $L \le U$
such that the optimal dynamic control policy is given by: For each $t \geq 0$,
\begin{align*}
\dot P(t)= \left\{\begin{array}{ll}
\theta_u, & \qquad \text{if} \quad P(t)-D(t)<L,\\
 0,       & \qquad \text{if} \quad P(t)-D(t) \in [L, U], \\
\theta_\ell, & \qquad \text{if} \quad P(t)-D(t)>U.
\end{array}\right.
\end{align*}
Moreover, the values of $L$ and $U$ can be characterized by the following
three cases.
\begin{enumerate}[I.]
\item If either Condition (1a), (1b) or (1c) hold, we have $U \ge L \ge 0$ where $L$ and $U$ are uniquely determined by
\begin{eqnarray}
B_1 e ^{r_1 (L-U)} + J_1 e ^{r_2 (L-U)}+ A &=& 0, \label{eq:xhminusxf} \\
\frac{B_1 r_2}{r_1-r_2} e ^{r_1 (L-U)} + \frac{J_1 r_1}{r_1-r_2} e^{r_2 (L-U)} &=&
(r_1+r_2-s_1)(\alpha \mathcal{I}_p+{\mathcal{C}_+}) + (\mathcal{C}_+ + \mathcal{C}_-) s_1 e^{s_2 L}. \qquad \label{eq:xH}
\end{eqnarray}
\item If either Condition (2a) or (2b) hold, we have $L \; \le \; U \; \le \; 0,$
where $L$ and $U$ are uniquely determined by
\begin{eqnarray}
B_2 e ^{r_1 (U-L)} + J_2 e^{r_2 (U-L)}+ K &=& 0, \label{eq:xfminusxh} \\
\frac{B_2 r_2}{r_1-r_2} e^{r_1(U-L)} + \frac{J_2 r_1}{r_1-r_2} e^{r_2 (U-L)} &=&
(r_1+r_2-t_2)(\alpha \mathcal{D}_p+\mathcal{C}_-) + (\mathcal{C}_+ + \mathcal{C}_- )t_2 e^{t_1 U}. \qquad \label{eq:xF}
\end{eqnarray}
\item If none of the above conditions hold, we then have $U \; \ge \; 0 \; \ge \; L,$
where $L$ and $U$ are uniquely determined by
\begin{eqnarray}
B_1 e^{-r_1 U} + B_2 e^{-r_1 L} &=& B_3, \label{eq:xhxf2} \\
J_1 e^{-r_2 U} + J_2 e^{-r_2 L} &=& J_3. \label{eq:xhxf2c}
\end{eqnarray}
\end{enumerate}
\end{theorem}

Theorem~\ref{THM:CASE1} can be interpreted as follows.
The optimal dynamic control policy seeks to maintain $X(t)=P(t)-D(t)$
within the risk-hedging interval $[L,U]$ at all time $t$, taking no action
(i.e., making no change to $P(t)$) as long as $X(t) \in [L,U]$.
Whenever $X(t)$ falls below $L$, the optimal dynamic control policy pushes
toward the risk-hedging interval as fast as possible, namely at rate
$\theta_u$, thus increasing the primary resource capacity allocation.
Similarly, whenever $X(t)$ exceeds $U$, the optimal dynamic control policy
pushes toward the risk-hedging interval as fast as possible, namely at rate
$\theta_\ell$, thus decreasing the primary resource capacity allocation.
In each of the cases \textit{I}, \textit{II} and \textit{III},
the optimal threshold values $L$ and $U$ are uniquely
determined by two nonlinear equations.

%\textcolor{red}{$\alpha \rightarrow 0$? Average cost}

%\subsection{Remaining cases.}\label{sec:main:remaining}
%%
%We further establish our main results, analogous to those of Theorem~\ref{THM:CASE1},
%for all remaining possible conditions on the adjustment costs $\mathcal{D}_p$
%and $\mathcal{I}_p$.
%Due to space restrictions, however, these results are provided in the
%supplemental appendix.

\subsubsection{Special Case $\boldsymbol{\mathcal{I}_p=\mathcal{D}_p=0}$}\label{sec:main:special}
In the special case where the dynamic control policy incurs
no costs for making adjustments,
which may be of particular interest in some application domains,
Theorem~\ref{THM:CASE1} has the following reduced form.
\begin{corollary}\label{THM:SPECIAL}
Suppose there are no adjustment costs, namely $\mathcal{I}_p=\mathcal{D}_p=0$.
%in the stochastic optimization problem \eqref{optHF:obj}~--~\eqref{optHF:st3}.
Then there exists a constant $\delta$ such that the optimal dynamic
control policy is given by
\begin{eqnarray*}
\dot P(t)= \left\{\begin{matrix}
\theta_u, && \qquad \text{if} \quad P(t)-D(t)<\delta,\\
0 && \qquad \text{if} \quad P(t)-D(t)= \delta, \\
\theta_\ell, && \qquad \text{if} \quad P(t)-D(t)>\delta.
\end{matrix}\right.
\end{eqnarray*}
Moreover, $\delta$ can be given explicitly by
\begin{eqnarray} \label{eq:delta}
\delta = \left\{\begin{matrix}
\frac{1}{s_2}\ln (\frac{\mathcal{C}_+}{\mathcal{C}_++\mathcal{C}_-} \frac{s_1-t_2}{s_1})>0, &&
\qquad \text{if} \quad
B_1 + B_2 -B_3 >0, \\
 0, &&
\qquad \text{if} \quad
B_1 + B_2 -B_3  =0,  \\
\frac{1}{t_1}\ln (\frac{\mathcal{C}_-}{\mathcal{C}_++\mathcal{C}_-} \frac{s_1-t_2}{-t_2})<0,
&&
\qquad \text{if} \quad
B_1 + B_2 -B_3 <0 .
\end{matrix}\right. \qquad
\end{eqnarray}
%%where the constants $\mathcal{C}_+, \mathcal{C}_-, t_1,t_2, s_1, s_2$
%%are given in \eqref{eq:cpcm}, \eqref{eq:t1}, \eqref{eq:t2}, \eqref{eq:s1}
%%and \eqref{eq:s2}.
%where the constants $s_1,s_2$ and $t_1, t_2$ are as in
%%\eqref{eq:s1}, \eqref{eq:s2}, \eqref{eq:t1}, \eqref{eq:t2}.
%\eqref{eq:s2} and \eqref{eq:t2}, respectively.
\end{corollary}

The interpretation of this corollary is the same as that for
Theorem~\ref{THM:CASE1} where the risk-hedging interval collapses to a
single point $\delta$.
Hence, the optimal dynamic control policy seeks to maintain $X(t)=P(t)-D(t)$
at the position $\delta$ at all time $t$, pushing toward this point as fast
as possible with rate $\theta_u$ when below and with rate $\theta_\ell$ when
above.
%Note that Corollary~\ref{THM:SPECIAL} does not specify any action whenever
%$P(t)=D(t)$ because the set $\{t: P(t)=D(t)\}$ has measure 0.
%Note that the set $\{t: P(t)=D(t)\}$ has measure 0.

\subsection{Case 2: $\boldsymbol{\mathcal{D}_p \ge \mathcal{C}_+ / \alpha}$ and $\boldsymbol{\mathcal{I}_p <  \mathcal{C}_- / \alpha}$}\label{sec:main:case2}
We next consider the case in which the per-unit cost for decreasing the
decision variable $P(t)$ is at least as large as the discounted overage
cost associated with this decision variable.
Our main result for this case can be expressed as follows.
\begin{theorem}\label{THM:CASE2}
Suppose the adjustment costs satisfy $\mathcal{D}_p \ge \mathcal{C}_+ / \alpha$
and $\mathcal{I}_p < \mathcal{C}_- / \alpha$.
Then there exists a threshold $L$ such that the optimal policy is given by
\begin{eqnarray*}
\dot P(t)= \left\{\begin{matrix}
\theta_u, && \qquad \text{if} \quad P(t)-D(t)<L,\\
0, && \qquad \text{if} \quad P(t)-D(t) \ge L.
\end{matrix}\right.
\end{eqnarray*}
Moreover, the value of $L$ can be characterized by
\begin{eqnarray}
L = \left\{\begin{array}{ll} \frac{1}{s_2} \ln \left[ {\frac{(\alpha
\mathcal{I}_p +\mathcal{C}_+)(s_1 -r_2)}{ (\mathcal{C}_+ +\mathcal{C}_-) s_1}}\right]
  \ge 0, & \qquad \text{if} \quad B_3 \le B_2, \label{eq:thm2Lpos}\\
%&\\
\frac{1}{-r_1} \ln \left[{\frac{(\mathcal{C}_++\mathcal{C}_-)( -r_2)}{ (\mathcal{C}_- - \alpha \mathcal{I}_p)
(s_1-r_2)}}\right] < 0, & \qquad \text{if} \quad B_3 > B_2 \label{eq:thm2Lneg}.
\end{array}\right.
\end{eqnarray}
\end{theorem}

This result is closely related with Theorem~\ref{THM:CASE1}. One readily checks that in Theorem~\ref{THM:CASE2} when $B_3 \le B_2$, $L \ge 0$ satisfies
\begin{eqnarray*}
 \frac{r_1}{r_1-r_2} (-A) &=&
(r_1+r_2-s_1)(\alpha \mathcal{I}_p+{\mathcal{C}_+}) + (\mathcal{C}_+ + \mathcal{C}_-) s_1 e^{s_2 L},
\end{eqnarray*}
which has a similar structure as \eqref{eq:xH}. When $B_3 >B_2 $, $L<0$ solves
\begin{eqnarray*}
B_2  e^{-r_1 L} = B_3,
\end{eqnarray*}
which is the same as Equation \eqref{eq:xhxf2} if we regard $U=\infty$. Therefore, given the relatively larger cost for decreasing $P(t)$, this theorem
essentially provides a one-sided version of Theorem~\ref{THM:CASE1}
in which the optimal dynamic control policy seeks to maintain
$X(t)=P(t)-D(t)$ at or above the threshold $L$ at all time $t$.
Whenever $X(t)$ falls below $L$, the optimal dynamic control policy
pushes toward the risk-hedging threshold as fast as possible,
namely at rate $\theta_u$.
Otherwise, the optimal dynamic control policy takes no action,
because taking action to decrease an overage position costs more
than the benefits from such an action.

\subsection{Case 3: $\boldsymbol{\mathcal{D}_p < \mathcal{C}_+ / \alpha}$ and $\boldsymbol{\mathcal{I}_p \ge \mathcal{C}_- / \alpha}$}\label{sec:main:case3}
We now consider the case in which the per-unit cost for increasing the
decision variable $P(t)$ is at least as large as the discounted shortage
cost associated with this decision variable.
Our main result for this case can be expressed as follows.
\begin{theorem}\label{THM:CASE3}
Suppose the adjustment costs satisfy $\mathcal{D}_p < \mathcal{C}_+ / \alpha$
and $\mathcal{I}_p \ge \mathcal{C}_- / \alpha$.
Then there exists a threshold $U$ such that the optimal policy is given by
\begin{eqnarray*}
\dot P(t)= \left\{\begin{matrix}
0, && \qquad \text{if} \quad P(t)-D(t) \le U,\\
\theta_\ell, && \qquad \text{if} \quad P(t)-D(t)>U.
\end{matrix}\right.
\end{eqnarray*}
Moreover, the value of $U$ can be characterized by
\begin{eqnarray} \label{eq:thm3U}
U= \left\{\begin{array}{ll} \frac{1}{-r_2} \ln \left[ \frac{(\mathcal{C}_+ +
\mathcal{C}_-)r_1}{ (\mathcal{C}_- -\alpha \mathcal{D}_p)(r_1 -t_2)} \right]
  \ge 0, & \qquad \text{if} \quad J_1 \le J_3,\\
%&\\
\frac{1}{t_1} \ln \left[\frac{(\alpha \mathcal{D}_p + \mathcal{C}_-)(r_1-t_2)}{(\mathcal{C}_+ + \mathcal{C}_-)
(-t_2)}\right] < 0, & \qquad \text{if} \quad J_1 > J_3.
\end{array}\right.
\end{eqnarray}
\end{theorem}

We also observe the connection of this result with Theorem~\ref{THM:CASE1}. It can be readily verified from Theorem~\ref{THM:CASE3} that when $J_1 \le J_3$, $U \ge 0$ satisfies
\begin{eqnarray*}
J_1  e^{-r_2 U} = J_3,
\end{eqnarray*}
and thus is equivalent to \eqref{eq:xhxf2c} if we regard $L = -\infty$. When $J_1 >J_3$, we have $U<0$ solving
\begin{eqnarray*}
\frac{r_2}{r_1-r_2} (-K) &=&
(r_1+r_2-t_2)(\alpha \mathcal{D}_p+\mathcal{C}_-) + (\mathcal{C}_+ + \mathcal{C}_- )t_2 e^{t_1 U},
\end{eqnarray*}
which is closely related to \eqref{eq:xfminusxh} and \eqref{eq:xF}. Hence, given the relatively larger cost for increasing $P(t)$, this theorem
essentially provides a one-sided version of Theorem~\ref{THM:CASE1}
in which the optimal dynamic control policy seeks to maintain
$X(t)=P(t)-D(t)$ at or below the threshold $U$ at all time $t$.
Whenever $X(t)$ exceeds $U$, the optimal dynamic control policy
pushes toward the risk-hedging threshold as fast as possible,
namely at rate $\theta_\ell$.
Otherwise, the optimal dynamic control policy takes no action,
because taking action to increase a shortage position costs more
than the benefits from such an action.

\subsection{Case 4: $\boldsymbol{\mathcal{D}_p \ge \mathcal{C}_+ / \alpha}$ and $\boldsymbol{\mathcal{I}_p \ge  \mathcal{C}_- / \alpha}$}\label{sec:main:case4}
Lastly, we consider the case in which the per-unit costs for adjusting
the decision variable $P(t)$ are at least as large as the corresponding
discounted overage and shortage costs associated with this decision variable.
We now state our main result for this case.
\begin{theorem}\label{THM:CASE4}
Suppose the adjustment costs satisfy $\mathcal{D}_p \ge \mathcal{C}_+ / \alpha$
and $\mathcal{I}_p \ge  \mathcal{C}_- / \alpha$.
Then the optimal policy consists of taking no action.
Specifically,
\[P (t) \; \equiv \; P(0), \qquad \text{for all $t$.} \]
\end{theorem}

Given the relatively larger costs for adjusting $P(t)$, this theorem
essentially consists of the inaction sides of both Theorems~\ref{THM:CASE2}
and~\ref{THM:CASE3}.
Intuitively, the theorem characterizes the conditions under which the
costs of any control policy action exceeds the resulting benefit, namely
taking an action to decrease an overage position or increase a shortage
position costs more than the benefits from such an action.

\section{Main Results: Multiple Primary Resources with Contract}\label{sec:main2}
In this section we consider our main results on the optimal dynamic control policy for the stochastic optimal control problem
when there are multiple primary resources under a contract-based relationship.
%, corresponding to the system model of Section~\ref{sec:model:Multiple}.
%To simplify the presentation, we shall suppress the single primary resource index
%$i=1$ in our notation throughout this section and, as an example, simply write
%$P(t)$ for the primary resource allocation at time $t$.
After providing a formulation of the stochastic optimal control problem and some technical preliminaries,
we present our main results analogous to those in Section~\ref{sec:main}.
%under different
%conditions for the values of $\mathcal{I}_p \geq 0$ and $\mathcal{D}_p \geq 0$,
%including the special case in which both are zero.
%%Consideration of the proofs of our main results in these cases is postponed until
%%the next section.

\subsection{Problem formulation}\label{sec:formulation:Multiple}
The stochastic optimal control problem associated with the system model of Section~\ref{sec:model:Multiple}
allows the dynamic control policy to adjust its allocation positions in primary and secondary resource capacities,
while maintaining the contract-based relationship $w$, based on the demand realization observed up to the current time.
More formally, the decision processes $P_1(t), \ldots, P_{\mathcal{P}}(t)$, are adapted to the filtration $\mathcal{F}_t$ generated by $\{D(s): s \le t\}$.
%Furthermore, any adjustments to the primary resource allocation capacities
%have associated costs, where we write $\mathcal{I}_{p,i}$ and $\mathcal{D}_{p,i}$
%to denote the per-unit costs of increasing and decreasing the decision process
%$P_i(t)$, respectively;
%namely, $\mathcal{I}_{p,i}$ represents the per-unit cost for increasing
%allocation of the $i$th primary resource capacity while $\mathcal{D}_{p,i}$
%represents the per-unit cost for decreasing allocation of the $i$th primary
%resource capacity.
Then the objective of the optimal dynamic control policy is to maximize the expected discounted net-benefit over an infinite horizon,
subject to the contract-based constraints \eqref{eq:multip}, where net-benefit at time $t$ consists of the difference between rewards and costs
from both the set of primary resource allocation capacities and the secondary resource allocation capacity and the additional costs for adjustments
to $P_1(t), \ldots, P_{\mathcal{P}}(t)$.

%In formulating the corresponding stochastic optimization problem, we impose
%a pair of additional conditions on each of the decision processes
%$\{P_i(t): t \ge 0\}$ based on practical aspects of the diverse application
%domains motivating our study.
%The control policy cannot instantaneously change the $i$th primary resource
%allocation capacity in an attempt to directly follow the remaining demand
%$D(t) - \sum_{j=1}^{i-1} P_j(t)$; i.e., some time is required (even if only
%a very small amount of time) to adjust $P_i(t)$.
Similar to the single primary resource formulation, the control policy cannot make unbounded adjustments in the $i$th
primary resource allocation capacity at any instant in time; i.e., the amount of change in $P_i(t)$ at time $t$ is
restricted (even if only to a very small extent) by various factors.
We therefore assume that the rate of change in the
$i$th primary resource allocation capacity by the control policy is bounded.
More precisely, there are pairs of finite constants $\theta_{\ell,i} < 0$ and
$\theta_{u,i} > 0$ such that \text{for each $t \ge 0$}
\begin{eqnarray} \label{eq:adjbound:Multiple}
\theta_{\ell,i} \; \le \; \dot P_i(t) \; \le \; \theta_{u,i},
\end{eqnarray}
where $\dot P_i(t)$ denotes the derivative of the decision variable $P_i(t)$
with respect to time, for all $i=1, \ldots, \mathcal{P}$.
On the other hand, the ability of the control policy to make adjustments to
the secondary resource capacity,
in response to changes in the primary resource capacities,
tends to be more flexible such that \eqref{eq:secondary:Multiple} holds at all time $t$.

Now we can present the mathematical formulation of our stochastic optimization problem
for the case of multiple primary resources with a contract-based relationship among these resources.
Let us fix a given contract vector $w$.
Defining
\begin{eqnarray*}
N_{p,i}(t) \; := \; R_{p,i}(t) - C_{p,i}(t)  & \qquad\qquad \mbox{and} \qquad\qquad & N_s(t) \; := \; R_s(t) - C_s(t) ,
\end{eqnarray*}
we seek to determine the optimal dynamic control policy for the contract vector $w$ that solves the problem (SC-OPT:M)
\begin{eqnarray}
\max_{P_1(\cdot), \ldots, P_{\mathcal{P}}(\cdot)} \quad && \E \int_0^\infty e^{-\alpha t} \left[ \sum_{i=1}^{\mathcal{P}} N_{p,i}(t) + N_s(t) \right] dt -
\sum_{i=1}^{\mathcal{P}} \E \int_0^\infty e^{-\alpha t} [ \mathcal{I}_{p,i} \cdot \indi{\dot P_i(t)>0} ] dP_i(t) \nonumber \\
&& \qquad\qquad\qquad\qquad\qquad\qquad\qquad\quad - \; \sum_{i=1}^{\mathcal{P}} \E \int_0^\infty e^{-\alpha t} [ \mathcal{D}_{p,i} \cdot \indi{\dot P_i(t)<0} ] d(-P_i(t)) \qquad \label{opt1:obj:Multiple} \\
\mbox{s.t.}
&& -\infty \; < \; \theta_{\ell,i} \; \le \; \dot P_i(t) \; \le \; \theta_{u,i} \; < \; \infty , \qquad\qquad \forall i = 1, \ldots, \mathcal{P} \quad \mbox{and $t \ge 0$}, \label{opt1:st1:Multiple} \\
&& P_i(t) \; = \; w_i P(t), \qquad\qquad\qquad\qquad\qquad\qquad \forall i = 1, \ldots, \mathcal{P} \quad \mbox{and $t \ge 0$}, \label{opt1:st0a:Multiple} \\
&& \sum_{i=1}^{\mathcal{P}} w_i =1 , \qquad\qquad\qquad\qquad\qquad\qquad\qquad w=(w_1, \ldots, w_{\mathcal{P}}) \in \R_+^{\mathcal{P}} , \label{opt1:st0b:Multiple} \\
&& dD(t) \; = \; b dt +\sigma d{W}(t), \label{opt1:st2:Multiple}
\end{eqnarray}
where $\alpha$ is the discount factor and $\indi{A}$ again denotes the indicator function associated with event $A$.
The control variables are the rates of change in the primary resource capacities by the control policy at every time $t$ subject to
the lower and upper bounds on each $\dot P_i(t)$ in \eqref{opt1:st1:Multiple} and the contract-based relationship among the primary resources in
\eqref{opt1:st0a:Multiple} and \eqref{opt1:st0b:Multiple}.
Note that the second (third) expectation in \eqref{opt1:obj:Multiple} causes a decrease with rate $\mathcal{I}_{p,i}$ ($\mathcal{D}_{p,i}$) in the value of the
objective function whenever the control policy increases (decreases) $P_i(t)$.

\subsection{Preliminaries}\label{sec:main:prelim:Multiple}
Consider the first expectation in the objective function \eqref{opt1:obj:Multiple} of the stochastic optimization problem (SC-OPT:P).
From \eqref{eq:multip}, \eqref{eq:Rpw} and \eqref{eq:benefitP:Multiple}, we can rewrite this expectation in terms of the aggregate primary
resource capacity $P(t)$ for a given contract vector $w$.
Upon analogously applying the derivations of Section~\ref{sec:main:prelim:Single}, we then can simplify and reduce this expectation to the
single primary resource setting as follows:
\begin{equation*}
\E \bigg[ \int_0^\infty e^{-\alpha t} \Big\{ \left( \mathcal{C}_+ (w)X(t)^+ + \mathcal{C}_-(w) X(t)^{-} \right) dt \Big\} \bigg]
\end{equation*}
where $X(t) := P(t) - D(t)$, $\mathcal{C}_+(w) = \mathcal{C}_p(w)$, $\mathcal{C}_-(w) = \mathcal{N}_p(w) - \mathcal{N}_s$,
$\mathcal{N}_p(w) := \mathcal{R}_p(w) - \mathcal{C}_p(w)$, and $\mathcal{N}_s := \mathcal{R}_s - \mathcal{C}_s$.
Similarly, the second and third expectations in \eqref{opt1:obj:Multiple} can be rewritten with respect to \eqref{eq:Ipw} in terms of the aggregate
primary resource capacity $P(t)$ for the given contract vector $w$.
Define $\mathcal{N}_{p,i} := \mathcal{R}_{p,i} - \mathcal{C}_{p,i}$, for $i=1,\ldots,\mathcal{P}$.

Next, we deduce from \eqref{eq:multip}, (\ref{eq:adjbound:Multiple}) and the system model definition of Section~\ref{sec:model:Multiple} that the aggregate
control process $P(t)$ satisfies the constraint
\begin{eqnarray} \label{eq:constraintw}
{\tilde \theta}_\ell  \; \le \; \dot P(t) \; \le \; {\tilde \theta}_u, \quad \mbox{for each $t \ge 0$},
\end{eqnarray}
where
\begin{align}\label{eq:tilde-theta}
{\tilde \theta}_\ell := \max_{i=1, \ldots, \mathcal{P}} \frac{\theta_{\ell,i}}{w_i} \qquad \mbox{and} \qquad {\tilde \theta}_u:= \min_{i=1, \ldots, \mathcal{P}} \frac{\theta_{u,i}}{w_i}.
\end{align}
Then, for a given contract vector $w$, we have the following stochastic optimization problem in terms of the aggregate control process $P(t)$
that is equivalent to the original problem formulation (SC-OPT:M):
\begin{align}
\min_{P(\cdot)} & \quad \E_x \bigg[ \int_0^\infty e^{-\alpha t} \Big\{ \left( \mathcal{C}_+ (w)X(t)^+ + \mathcal{C}_-(w) X(t)^{-} \right) dt
+ \left( \mathcal{I}_p(w) \indi{\dot P(t)>0} - \mathcal{D}_p(w) \indi{\dot P(t)<0} \right) dP(t) \Big\} \bigg] \qquad \label{optHF:obj_w} \\
 \mbox{s.t.} & \quad -\infty \; <{\tilde \theta}_\ell  \; \le \; \dot P(t) \; \le \; {\tilde \theta}_u \; < \; \infty , \qquad\qquad \mbox{for $t \ge 0$}, \nonumber \\
 & \quad dX(t) \; = \; dP(t) - b dt - \sigma d{W}(t), \qquad\qquad \mbox{for $t \ge 0$}, \nonumber \\
& \quad X(0) \; = \; x , \qquad \mathcal{C}_+ (w) \; = \; \mathcal{C}_p (w) , \qquad \mathcal{C}_- (w) \; = \; \mathcal{N}_p(w) - \mathcal{N}_s  , \nonumber
\end{align}
where $\E_x[\cdot]$ again denotes expectation with respect to the initial state
distribution (i.e., state at time $t=0$) being $x$ with probability one.

Once again, we use $V_w(x)$ to represent the optimal value of the objective function in
\eqref{optHF:obj_w} for a given contract vector $w$; namely, $V_w(x)$ is the value function of the corresponding stochastic dynamic program.
Given its equivalence with the original optimization problem (SC-OPT:M), the remainder of this section will
focus on the stochastic dynamic program formulation in \eqref{optHF:obj_w}.

For notational convenience, we define the following constants which represent modifications of some of the constants
in Section~\ref{sec:main:prelim:Single} due to the differences in the parameters and the problem setting:
\begin{eqnarray*}
%r_1 \; := \; \frac{b+ \sqrt{b^2 +2 \alpha \sigma^2 }}{\sigma^2} \; > \; 0, \label{eq:r1}
%& \qquad &
%r_2 \; := \; \frac{b- \sqrt{b^2 +2 \alpha \sigma^2 }}{\sigma^2} \; < \; 0, \label{eq:r2} \\
{\tilde s}_1 \; := \; \frac{b- {\tilde \theta}_u+ \sqrt{(b-{\tilde \theta}_u)^2 +2 \alpha \sigma^2 }} {\sigma^2} \; > \; 0, \label{eq:s1_w}
& \qquad &
{\tilde s}_2 \; := \; \frac{b- {\tilde \theta}_u - \sqrt{(b-{\tilde \theta}_u)^2 +2 \alpha \sigma^2 }} {\sigma^2} \; < \; 0 , \label{eq:s2_w} \\
{\tilde t}_1 \; := \; \frac{b- {\tilde \theta}_\ell + \sqrt{(b-{\tilde \theta}_\ell)^2 +2 \alpha \sigma^2 }} {\sigma^2} \; > \; 0, \label{eq:t1_w}
& \qquad &
{\tilde t}_2 \; := \; \frac{b- {\tilde \theta}_\ell - \sqrt{(b-{\tilde \theta}_\ell)^2 +2 \alpha \sigma^2 }} {\sigma^2} \; < \; 0 . \label{eq:t2_w}
\end{eqnarray*}
and
\begin{eqnarray*}
	B_1(w) := ({\mathcal{C}_+}(w)-\alpha \mathcal{D}_p(w)) ({\tilde t}_2 -r_2), & \qquad &
	B_2(w) := (\mathcal{C}_-(w) - \alpha \mathcal{I}_p(w)) ({\tilde s}_1 -r_2), \\
	B_3(w) := (\mathcal{C}_+(w) + \mathcal{C}_-(w)) (-r_2),  & \qquad &
	J_1(w) := ({\mathcal{C}_+}(w)-\alpha \mathcal{D}_p(w)) (r_1 -{\tilde t}_2), \\
	J_2(w) := (\mathcal{C}_-(w) - \alpha \mathcal{I}_p(w)) (r_1 -{\tilde s}_1),  & \qquad &
	J_3(w) := (\mathcal{C}_+(w) + \mathcal{C}_-(w)) r_1, \\
	A(w) := (\mathcal{C}_+(w) + \alpha \mathcal{I}_p(w)) (r_2-r_1), & \qquad &
	K(w) := (\mathcal{C}_-(w) + \alpha \mathcal{D}_p(w))(r_2-r_1) .
\end{eqnarray*}

\subsection{Case 1: $\boldsymbol{\mathcal{D}_{p}(w) < \mathcal{C}_{+}(w) / \alpha}$ and $\boldsymbol{\mathcal{I}_p(w) < \mathcal{C}_-(w) / \alpha}$}\label{sec:main:case1:Multiple}
Let us first briefly explain the conditions of this subsection, for a given contract vector $w \in \Omega$ where
\begin{equation*}
%\label{eq:omega}
\Omega  := \{ w = (w_1, \ldots, w_{\mathcal{P}}) \in \R^{\mathcal{P}}: w \ge 0, \sum_{i=1}^{\mathcal{P}} w_i =1\}.
\end{equation*}
The conditions of this subsection are likely to be the most relevant case in practice.
Observe that,
if $\mathcal{D}_{p,i} < \mathcal{C}_{p,i} / \alpha$ and $\mathcal{I}_{p,i} < (\mathcal{N}_{p,i} -\mathcal{N}_{s})/ \alpha$ for each primary resource type $i=1, \ldots, \mathcal{P}$,
then the conditions $\mathcal{D}_{p}(w) < \mathcal{C}_{+}(w) / \alpha$ and $\mathcal{I}_p(w) < \mathcal{C}_-(w) / \alpha$ hold for any $w \in \Omega$.
Further observe from the objective function in \eqref{optHF:obj_w} that $\mathcal{C}_+(w)/\alpha$ reflects the discounted overage cost associated
with the aggregate primary resource capacity and $\mathcal{C}_-(w)/\alpha$ reflects the discounted shortage cost associated with the
aggregate primary resource capacity, recalling that $\alpha$ is the discount rate.
In comparison, $\mathcal{D}_p(w)$ represents the cost incurred for decreasing $P(t)$ when in an overage position while $\mathcal{I}_p(w)$
represents the cost incurred for increasing $P(t)$ when in a shortage position.

The conditions of this subsection represent counterparts of the conditions (1a), (1b), (1c), (2a) and (2b) of Section~\ref{sec:main:case1}.
To elucidate the exposition, we denote

{Condition (1a'):} \quad $\mathcal{I}_p(w) + \mathcal{D}_p(w) >0$, \quad $0 \; < \; {B_3(w)- B_2(w)} \; < \; {B_1(w)} \quad \mbox{and} \quad \left(\frac{B_3(w)- B_2(w)}{B_1(w)}\right)^{\frac{r_2}{r_1}} \; \ge \; \frac{J_3(w)-J_2(w)}{J_1(w)}$.

{Condition (1b'):} \quad $\mathcal{I}_p(w)(w) + \mathcal{D}_p(w) >0$, \quad \mbox{and} \quad $B_3(w) \le B_2(w)$.

{Condition (1c'):} \quad $\mathcal{I}_p(w) = \mathcal{D}_p(w) =0$, \quad \mbox{and} \quad $B_3(w) -B_2(w) -B_1(w) \le 0$.

{Condition (2a'):} \quad $\mathcal{I}_p(w) + \mathcal{D}_p(w) >0$, \quad ${B_3(w)-B_2(w)- B_1(w)} \; > \; 0 \qquad \mbox{and} \qquad \left(\frac{B_3(w)-B_1(w)}{B_2(w)}\right)^{\frac{r_2}{r_1}} \;  \ge \; \frac{J_3(w) -J_1(w)}{J_2}$.

{Condition (2b'):} \quad $\mathcal{I}_p(w) = \mathcal{D}_p(w) =0$, \quad \mbox{and} \quad $B_3(w) -B_2(w) -B_1(w) \ge 0$.

\noindent
We are now ready to state our main result for Case 1 under this setting.
Recall that ${\tilde \theta}_\ell$ and ${\tilde \theta}_u$ are defined in \eqref{eq:tilde-theta}.

\begin{theorem}\label{THM:CASE1:Multiple}
Fix $w \in \Omega$.
Suppose the adjustment costs satisfy $\mathcal{D}_p(w)<\mathcal{C}_+ (w)/\alpha$ and $\mathcal{I}_p (w)<\mathcal{C}_- (w)/\alpha$.
Then there are two threshold values $L(w)$ and $U(w)$ with $L(w) \le U(w)$
such that the optimal dynamic control policy is given by
\begin{align*}
\dot P(t)= \left\{\begin{array}{ll}
\tilde \theta_u, & \qquad \text{if} \quad P(t)-D(t)<L(w),\\
0,       & \qquad \text{if} \quad P(t)-D(t) \in [L(w), U(w)], \\
\tilde \theta_\ell, & \qquad \text{if} \quad P(t)-D(t)>U(w).
\end{array}\right.
\end{align*}
Moreover, the values of $L(w)$ and $U(w)$ can be characterized by the following
three cases.
\begin{enumerate}[I.]
\item If either Condition (1a'), (1b') or (1c') hold, we have $U(w) \ge L(w) \ge 0$ where $L(w)$ and $U(w)$ are uniquely determined by
\begin{eqnarray}
B_1(w) e ^{r_1 (L(w)-U(w))} + J_1(w) e ^{r_2 (L(w)-U(w))}+ A(w) &=& 0, \label{eq:xhminusxf_w} \\
\frac{B_1(w) r_2}{r_1-r_2} e ^{r_1 (L(w)-U(w))} + \frac{J_1(w) r_1}{r_1-r_2} e^{r_2 (L(w)-U(w))} &=& \nonumber
(r_1+r_2-{\tilde s}_1)(\alpha \mathcal{I}_p(w)
+{\mathcal{C}_+(w)})\\ && + (\mathcal{C}_+(w) + \mathcal{C}_-(w)) {\tilde s}_1 e^{{\tilde s}_2 L(w)}. \qquad \label{eq:xH_w}
\end{eqnarray}
\item If either Condition (2a') or (2b') hold, we have $L(w) \; \le \; U(w) \; \le \; 0,$
where $L(w)$ and $U(w)$ are uniquely determined by
\begin{eqnarray}
B_2(w) e ^{r_1 (U(w)-L(w))} + J_2 e^{r_2 (U(w)-L(w))}+ K(w) &=& 0, \label{eq:xfminusxh_w} \\
\frac{B_2(w) r_2}{r_1-r_2} e^{r_1(U(w)-L(w))} + \frac{J_2(w) r_1}{r_1-r_2} e^{r_2 (U(w)-L(w))} &=&
(r_1+r_2-{\tilde t}_2)(\alpha \mathcal{D}_p(w)+\mathcal{C}_-(w)) \nonumber \\ &&+ (\mathcal{C}_+(w) + \mathcal{C}_-(w) ){\tilde t}_2 e^{{\tilde t}_1 U(w)}. \qquad \label{eq:xF_w}
\end{eqnarray}
\item If none of the above conditions hold, we then have $U(w) \; \ge \; 0 \; \ge \; L(w),$
where $L(w)$ and $U(w)$ are uniquely determined by
\begin{eqnarray}
B_1(w) e^{-r_1 U} + B_2(w) e^{-r_1 L(w)} &=& B_3(w), \label{eq:xhxf2_w} \\
J_1(w) e^{-r_2 U} + J_2(w) e^{-r_2 L(w)} &=& J_3(w). \label{eq:xhxf2c_w}
\end{eqnarray}
\end{enumerate}
\end{theorem}

%\begin{corollary} \label{THM:multiCASE1}
%Fix $w \in \Omega$. Suppose $\mathcal{D}_p(w)<\mathcal{C}_+ (w)/\alpha$
%and $\mathcal{I}_p (w)<\mathcal{C}_- (w)/\alpha$. Then there are two threshold values $L(w)$ and $U(w)$
%such that the optimal policy is given by
%\begin{align*}
%\dot P(t)= \left\{\begin{array}{ll}
%\min_{i=1, \ldots, \mathcal{P}} \frac{\theta_{u,i}}{w_i} , & \qquad \text{if} \quad P(t)-D(t)<L(w),\\
% 0,       & \qquad \text{if} \quad P(t)-D(t) \in [L(w), U(w)], \\
%\max_{i=1, \ldots, \mathcal{P}} \frac{\theta_{\ell,i}}{w_i}, & \qquad \text{if} \quad P(t)-D(t)>U(w),
%\end{array}\right.
%\end{align*}
%and
%\[ \dot P_i(t) = w_i \dot P(t), \quad \text{for $i=1, \ldots, \mathcal{P}$}. \]
%The two threshold values $L(w)$ and $U(w)$ can be uniquely determined by two nonlinear equations similarly as in Theorem~\ref{THM:CASE1}.
%\end{corollary}

Theorem~\ref{THM:CASE1:Multiple} can be interpreted similar to Theorem~\ref{THM:CASE1}.
Given a fixed contract $w$ among primary resource options, the optimal dynamic control policy seeks to maintain $X(t)$
within the risk-hedging interval $[L(w),U(w)]$ at all time $t$.
When outside this risk-hedging interval $[L(w), U(w)]$, the optimal dynamic control policy pushes toward this interval as fast as possible in a
synchronized way such that the contract condition (\ref{eq:multip}) is always maintained among the primary resource allocations $P_1(t), \ldots, P_{\mathcal{P}}(t)$;
namely, $\dot P_i(t) = w_i \dot P(t)$ for each $i$.
The optimal threshold values $L(w)$ and $U(w)$ are uniquely determined by two nonlinear equations for each of the cases $I$, $II$ and $III$.

%\textcolor{blue}{TODO: add Computational computation of $L(w)$, $U(w)$ and maybe value function $V_w(x)$ for different $w$ in two dimensions.}
We now explore the dependence of the value function on the contract $w$.
For a fixed contract $w \in \Omega$, and $x=P(0)-D(0)$,
we write $J_w(x)$ for the optimal value of the objective function \eqref{opt1:obj:Multiple}, i.e.,
the maximal expected discounted net-benefit over an infinite time horizon when one optimally adjusts the aggregate primary resource option
$P(t)=\sum_{i=1}^{\mathcal{P}} P_i(t)$ under the fixed contract $w$ among these primary options together with the secondary resource option in order to meet the Brownian demand.
Then it is easily seen from Section~\ref{sec:main:prelim:Single} (refer to \eqref{eq:simplify}) that
\begin{eqnarray}\label{eq:J-V}
J_w(x) = (\mathcal{R}_p (w) - \mathcal{C}_p (w)) \cdot \E \left[ \int_0^\infty e^{-\alpha t} D(t) dt \right] - V_{w}(x),
\end{eqnarray}
where ${{V}}_w(x)$ is the value function of the stochastic dynamic program \eqref{optHF:obj_w} with parameters given in
\eqref{eq:multip}~--~\eqref{eq:benefitS:Multiple} and $\mathcal{R}_p (w), \mathcal{C}_p (w)$ in \eqref{eq:Rpw} are linear in $w$.
A closer examination of the expression we obtained for Theorem~\ref{THM:CASE1} and its proof imply that the value function continuously depends on $w$.
Such continuity is the consequence of the continuity of the solution of the corresponding ordinary differential equation with respect to the initial condition and parameters,
as well as the ``smooth-fit'' principle.

%now the optimal value of the original optimization problem \eqref{opt1:obj:Multiple}--\eqref{opt1:st2:Multiple} with the additional constraint \eqref{eq:Rpw} depends on $w$. We use $J_w(x)$ to denote this value function where. That is, Mathematically, if we write we infer from

% where Thus it is enough to show $V_w(x)$ is continuous with respect to $w \in \Omega$.

%Here $P(0)$ is given. Suppose one has the option to select contracts, or equivalently, the initial confi
%%stochastic dynamic program \eqref{optHF:obj}-\eqref{optHF:st3}.
%
%The key idea of the proof is to show that the optimal threshold values $L(w)$ and $U(w)$ are continuous functions of $w$. %, which leads to the continuity of the value function $\mathcal{V}_w(x)$ with respect to $w$.
%Fix $x \in \R$. Since $\Omega$ defined in (\ref{eq:omega}) is a compact set in $\R^{\mathcal{P}}$, we deduce that in order to prove the existence of the optimal contract $w^{*}(x)$, it suffices to show that $\mathcal{V}_w(x)$ is continuous with respect to $w \in \Omega$.
%
%
%We are now ready to state our next result. %(\textcolor{red}{unique?})

\begin{theorem} \label{thm:contract}
Given any fixed $x \in \R_+$, the optimal threshold values $L(w)$ and $U(w)$ in Theorem~\ref{THM:CASE1:Multiple} are continuous functions of $w \in \Omega$. As a consequence,
$V_w(x)$ and $J_w(x)$ are continuous with respect to $w \in \Omega$. %As a consequence,
%there exists a optimal contract $w^*(x) \in \Omega$
%such that
%\begin{eqnarray*}
% \max_{w \in \Omega } J_w(x) = J_{w^*(x)} (x).
%\end{eqnarray*}
%Moreover, there exists $x^*$ and $w^*(x^*)$ such that
%\begin{eqnarray*}
%\max_{x} \max_{w \in \R^{\mathcal{P}}_{+}, \sum_{i=1}^{\mathcal{P}} {w}_i =1} \mathcal{V}_w(x) = \mathcal{V}_{w^*(x^*)} (x^*).
%\end{eqnarray*}
\end{theorem}

This result also suggests a simple possible scheme for our dynamic resource allocation problem in the presence of multiple primary resource options.
Given the initial imbalance $x$ between the aggregate primary resource capacity and the demand,
i.e., $x=P(0)-D(0)$, one can first solve offline an optimal contract $w^*(x)$ given the characteristics of the demand, the reward and the cost associated with each sourcing option.
Such an optimal contract exists due to the continuity result in Theorem~\ref{thm:contract}.
Next, in order to meet the uncertain and volatile demand over time, one fixes this optimal contract $w^*(x)$ among primary resource options throughout the time horizon and then
dynamically adjusts the primary sourcing capacities according to the threshold policy given in Theorem~\ref{THM:CASE1:Multiple}.
Note that the capacity of each individual primary resource is aligned according to the contract vector, i.e., the ratio $P_i(t)/P_j(t) \equiv w_i^*(x)/w_j^*(x)$ is fixed for all $t\ge 0$.
Any remaining demand will be served by the secondary resource option.

\subsection{Remaining Cases}
It is easy to see that similar results can be obtained for all remaining possible conditions on the adjustment costs $\mathcal{D}_p(w)$ and $\mathcal{I}_p(w)$.
In particular, a corollary to Theorem~\ref{THM:CASE1:Multiple} can be expressed and established analogous to Corollary~\ref{THM:SPECIAL} applied to the aggregate control process $P(t)$.
Similarly, the theorems corresponding to Theorems~\ref{THM:CASE2}~--~\ref{THM:CASE4} can be expressed and established analogous to these theorems applied to the aggregate control process $P(t)$.

\section{Computational Experiments}\label{sec:computational}
The foregoing sections establish the explicit optimal dynamic control policy among
all admissible nonanticipatory control processes $P(t)$ within a stochastic optimal
control setting that maximizes the stochastic dynamic programs (SC-OPT:S) and (SC-OPT:M).
These optimal dynamic control policies render a new class of online algorithms for
general dynamic resource allocation problems that arise across a wide variety of application domains.
The resulting online algorithms are easily implementable in computer systems and communication
networks (among others) at runtime and consists of maintaining
$X(t) = P(t)-D(t)$ respectively
within the risk-hedging intervals $[L,U]$ and $[L(w),U(w)]$ at all time $t$,
where $L$, $U$, $L(w)$ and $U(w)$ are easily obtained from application parameters.
In this section, we present a representative sample of computational experiments
conducted across a broad spectrum of application
environments to investigate various issues of both theoretical and practical interest by
comparing our online optimal dynamic control algorithm against alternative optimization approaches
from recent work in the research literature.

\begin{figure}[ht!]
\centering
\includegraphics[width=.43\textwidth]{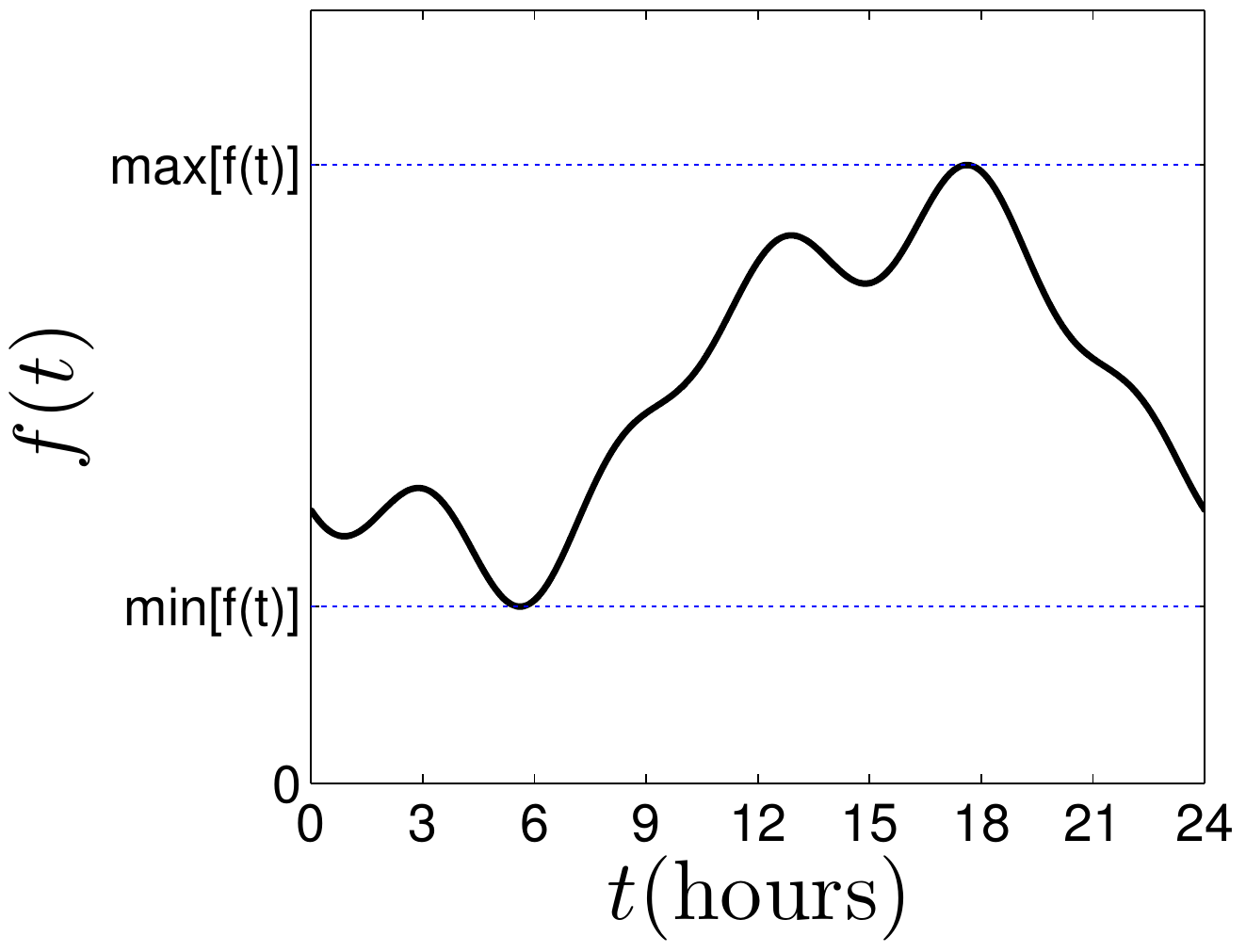}
\hspace*{0.6in}
\includegraphics[width=.43\textwidth]{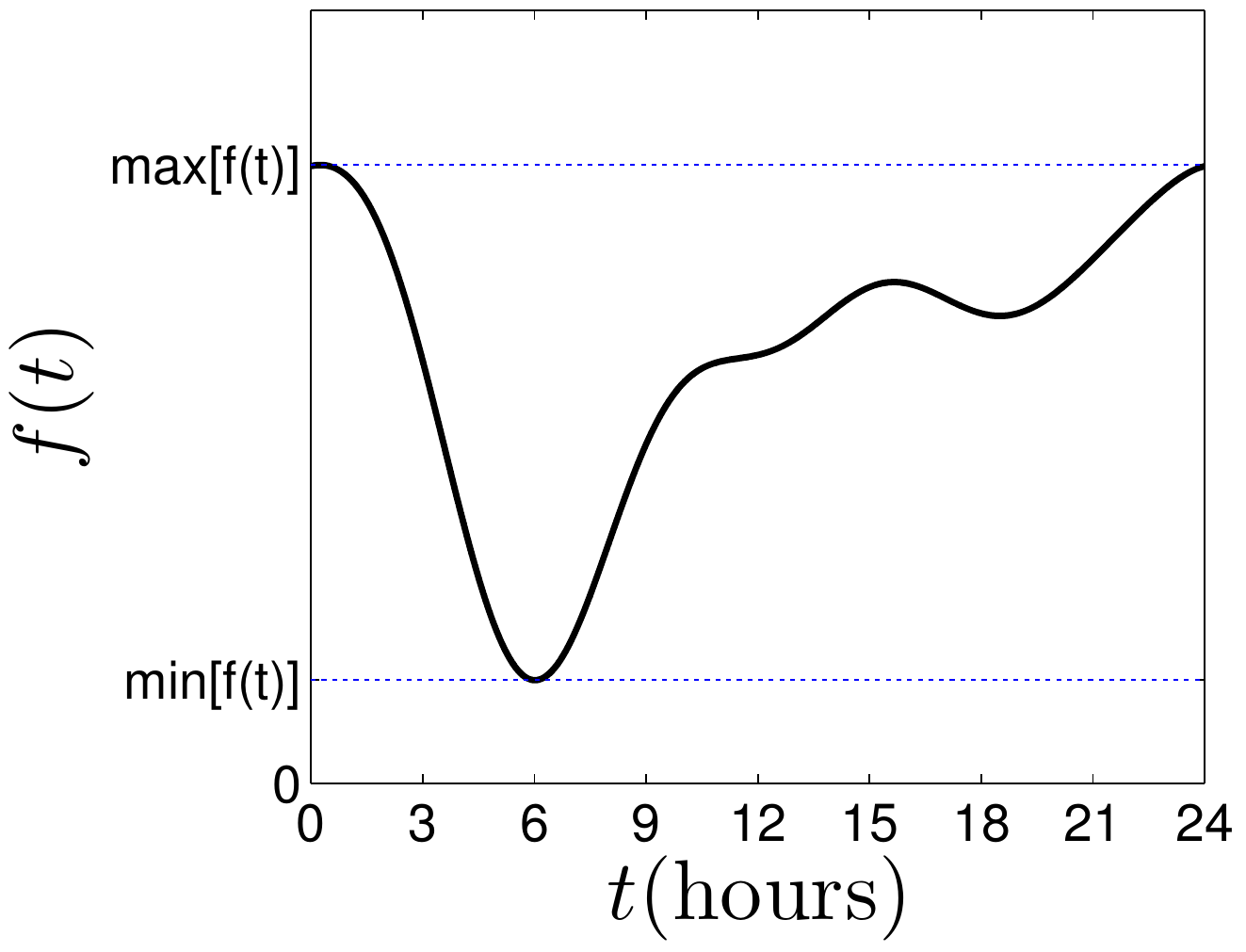}
\caption{Representative average daily demand patterns $f^1(t)$ (left) and $f^2(t)$ (right).}
\label{fig:fx}
\end{figure}

Through a detailed analysis of real-world trace data~\cite{GaLuSh+13}, we fitted
average daily demand processes for different environments by
smooth functions $f^1(t)$ and $f^2(t)$, depicted in Figure~\ref{fig:fx}.
In addition, our detailed analysis of real-world data revealed a wide range of volatility in the
demand process over time, as well as from one environment to another.
We therefore focus in the remainder of this section on the average daily demand patterns
$f^1(t)$ and $f^2(t)$
%in Figure~\ref{fig:fx}
for the drift parameter $b$ of the demand process while varying its volatility parameter $\sigma$
(as made more precise below),
thus representing a broad spectrum of application environments.

\subsection{Single Primary Resource}
\label{sec:exp:single}
For comparison with our optimal online dynamic control algorithm, we consider two alternative optimization approaches
that have recently appeared in the research literature.
First, Lin et al.~\cite{LiWiAn+11} propose an optimal offline algorithm that consists of making
optimal provisioning decisions in a clairvoyant anticipatory manner based on the {\it{known average demand}}
within each slot of a discrete-time model where the slot length is chosen to match the timescale
at which a data center can adjust its resource capacity and so that demand activity within a slot is
sufficiently nonnegligible in a statistical sense.
Applying this particular optimal offline algorithm within our mathematical framework, we partition the daily
time horizon into $T$ slots of length $\gamma$ such that
%\begin{eqnarray*}
%h_i \; = \; (t_{i-1},t_i] , & \qquad\qquad & \gamma \; = \; t_i-t_{i-1} ,
%\end{eqnarray*}
$h_i = (t_{i-1},t_i]$, $\gamma = t_i-t_{i-1}$,
$i=1,\ldots,T$, $t_0 := 0$, and we compute the average demand
$g_i := \gamma^{-1} \int_{h_{i}}f(t)dt$
within each slot $i$ yielding the average demand vector
$(g_1, g_2, \ldots, g_T)$.
Define
$\Delta(P_i) := P_i - P_{i-1}$,
where $P_i$ denotes the primary resource allocation capacity for slot $i$.
The optimal solution under this offline algorithm is then obtained by solving the following
linear program (LP) for each sample path:
\begin{eqnarray}
\min_{\Delta(P_1),\ldots,\Delta(P_T)} && \sum_{i=1}^T \; \mathcal{C}_{+} (P_i-g_i)^+ \; + \; \mathcal{C}_{-} (P_i-g_i)^- \; + \;
    \mathcal{I}_p (P_i-P_{i-1})^+ \; + \; \mathcal{D}_p (P_i-P_{i-1})^- \label{offline:obj} \\
\mbox{s.t.} && -\infty \; < \; \theta_\ell \; \le \; \Delta(P_i)/\gamma \; \le \; \theta_u \; < \; \infty ,
    \qquad \forall i = 1, \ldots, T , \label{offline:st1}
\end{eqnarray}
where the constraints on the control variables $\Delta(P_i)$ in \eqref{offline:st1} correspond to \eqref{optHF:st1}.
We refer to this solution as the offline LP algorithm.

Second, we consider a related optimal online algorithm proposed by Ciocan and Farias~\cite{CioFar12}.
Although they focus on dynamic resource allocations with stochastic demand rate, their allocation scheme based on re-optimization heuristics (Section~3.1 in \cite{CioFar12}) is quite general and can be applied to our setting with stochastic demand process and deterministic demand rate.
The main idea of their algorithm is that at each discrete point in time, one uses demand information realized up to that point and assumes the demand rate over the remaining time horizon is unchanged, and then employs the allocation rule that is optimal for such a scenario over the interval of time until the next re-solve.
We refer to this as the online CF algorithm.

%A typical approach to solving planning problems of dynamic resource management is as follows.
%Suppose we have a mechanism to obtain estimates for the average consumption or usage of the
%resource in question over a given time period. Assume a discrete-time formulation in which the
%duration of a time slot $\gamma$ is carefully chosen to be long enough so as to: (i) exceed
%the minimum time it takes for the system to make changes in its resource levels, and (ii)
%allow sufficient averaging of system stochastics. Then the planner is asked to optimize, in an
%\emph{online} fashion, the resource levels in each time slot based on average resource requirements
%that are revealed time slot by time slot. Adopting the computer science methodology of \emph{competitive
%analysis}, the planner proposes an online approximation algorithm to solve the problem and
%provides performance guarantees in terms of a competitive ratio of the online solution relative
%to the \emph{offline} optimal solution, where all the averages are assumed to be available to the
%planner at once. We now present the results of comparing the performance of our algorithm against
%the offline algorithm.

The sample paths of demand for our computational experiments are generated from a
linear diffusion process for the entire time horizon such that the drift of the demand process is
obtained as the derivative of $f(t)$ (i.e., $b(t)=df(t)$) and the corresponding volatility term
is set to match $\sigma(t)$.
Since the volatility pattern $\sigma(t)$ tended to be fairly consistent with respect to time
within each daily real-world trace for a specific environment and since the volatility pattern
tended to vary considerably from one daily real-world trace to another, our linear diffusion
demand process is assumed to be governed by the following model
$$dD(t)=b(t)dt+\sigma dW(t),$$
where we vary the volatility term $\sigma$ to investigate different application environments.
Each workload then consists of a set of sample paths generated
from the Brownian demand process $D(t)$ defined in this manner.

Given such a demand process, we calibrate our optimal online dynamic control algorithm by
first partitioning the average daily demand function $f(t)$ into piecewise linear segments,
then correspondingly setting the drift function $b(t)$ of the demand process $D(t)$, and
finally computing the threshold values $L$ and $U$ for each per-segment drift and $\sigma$ according to Theorem~\ref{THM:CASE1}.
This (fixed) version of our optimal online dynamic control algorithm is applied to every daily
sample path of the Brownian demand process $D(t)$ and the time-average value of net-benefit is
computed over this set of daily sample paths.
For comparison under the same set of Brownian demand process sample paths, we compute the average
demand vector $(g_1,\ldots,g_T)$ and the corresponding solution under the offline LP algorithm for
each daily sample path by solving the linear program \eqref{offline:obj},\eqref{offline:st1} with
respect to $(g_1,\ldots,g_T)$, and then we calculate the time-average value of net-benefit over
the set of daily sample paths.
The corresponding computational experiments for the CF algorithm are performed within this discrete-time
framework and the corresponding time-average net-benefit is computed in a similar manner.
All of our computational experiments were implemented in Matlab using, among other functionality,
the econometrics toolbox.

For our first set of results based on the average daily demand pattern $f^1(t)$ illustrated in the
leftmost plot of Figure~\ref{fig:fx}, the base parameter settings are given by $\alpha=0.02$, $\sigma=0.4$,
$\theta_{l}=-10$, $\theta_{u}=10$, $\mathcal{C}_{+}=20$, $\mathcal{C}_{-}=2$,
$\mathcal{D}_{p}=0.5$, $\mathcal{I}_{p}=0.5$, $f^1_{\mbox{\tiny min}}=2$, $f^1_{\mbox{\tiny max}}=7$,
$f^1_{\mbox{\tiny avg}}=4.5$ and $x = X(0) = P(0) - D(0) = 0$, where
$f_{\mbox{\tiny min}} := \min_t\{f(t)\}$, $f_{\mbox{\tiny max}} := \max_t\{f(t)\}$
and
$f_{\mbox{\tiny avg}} := T^{-1} \int_0^T f(t)dt$.
In addition to these base settings, we vary the parameter values one at a time for
$\sigma \in [0.01, 1.0]$, $\mathcal{C}_{+} \in [10,40]$, $\mathcal{C}_{-} \in [1,10]$,
$f^1_{\mbox{\tiny min}} \in [1,5]$ and $f^1_{\mbox{\tiny max}} \in [4,25]$, in order to
investigate the impact and sensitivity of these parameters on the performance of the various optimization algorithms.
For each computational experiment under a given set of parameters, we generate $N=10,000$ daily
sample paths using a timescale of a couple of seconds and a $\gamma$ setting of five minutes,
noting that a wide variety of experiments with different timescale and $\gamma$ settings provided
the same performance trends as those presented herein.
We then apply our optimal dynamic control policy and the two alternative optimization approaches
to this set of $N$ daily sample paths as described above.
%, where our performance evaluation comparison
%is based on the expectation of net-benefit realized under each of the three algorithms, also as
%described above.
%In particular, the expected net-benefit is computed as the time-average value of the rewards minus
%the costs from the primary and secondary resource allocation capacities and minus the costs for
%adjustments to the primary resource allocation capacity, taken over all $N$ daily sample paths
%under each of our optimal online dynamic control algorithm and the offline LP and online CF
%algorithms.

\begin{figure}[ht!]
\centering
\includegraphics[width=0.49\textwidth]{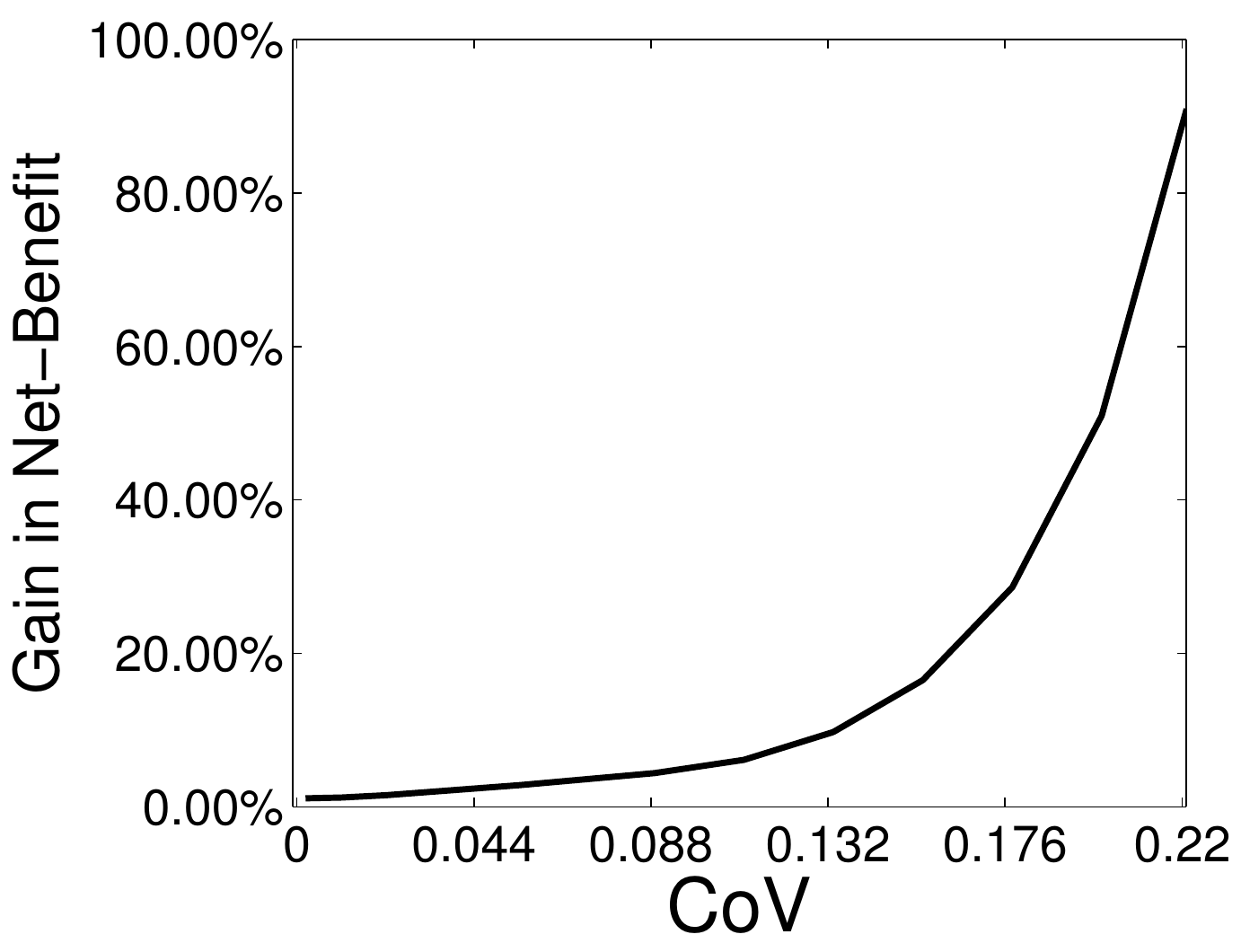}
\hspace*{0.0005in}
\includegraphics[width=0.49\textwidth]{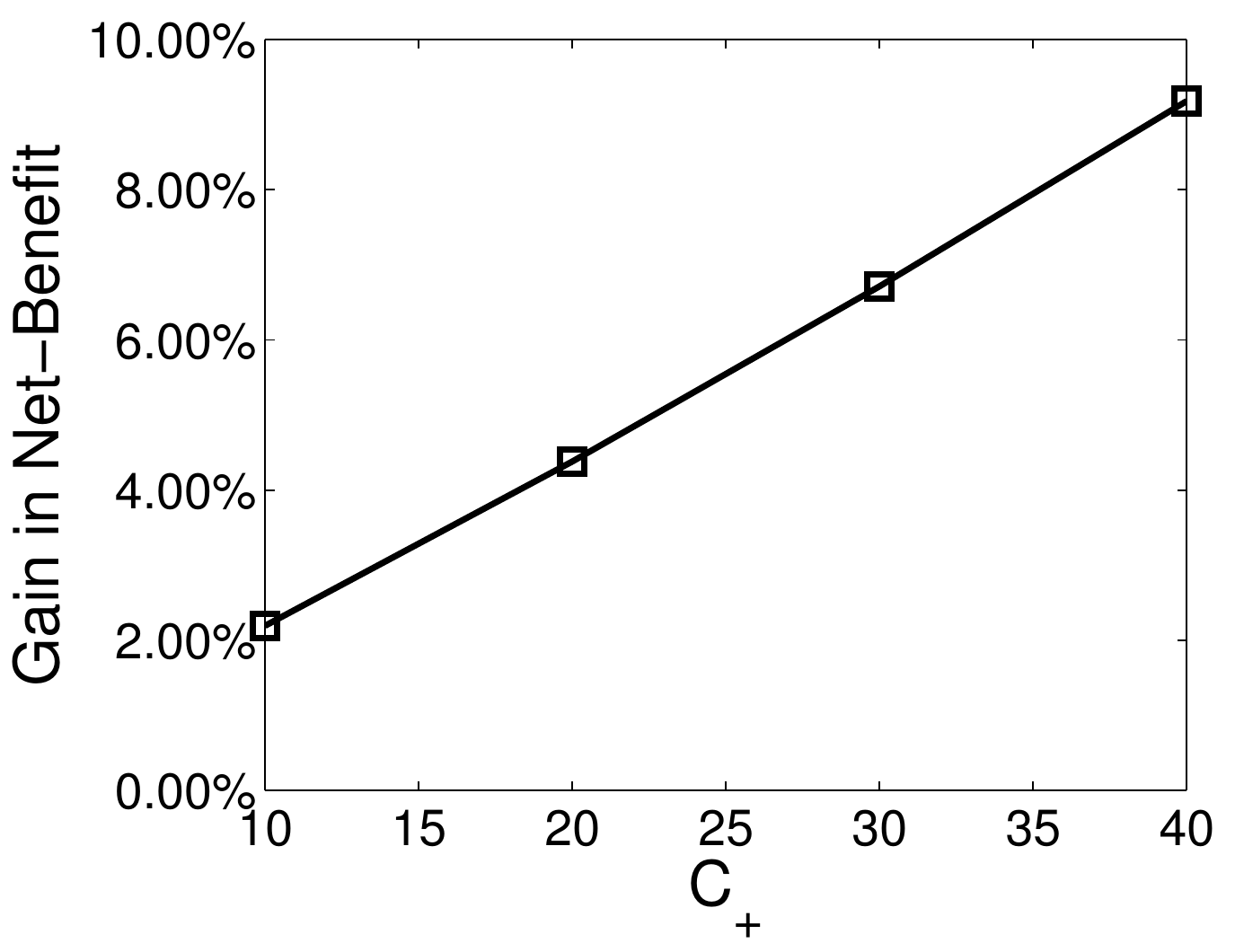}
\hspace*{0.0005in}
\includegraphics[width=0.49\textwidth]{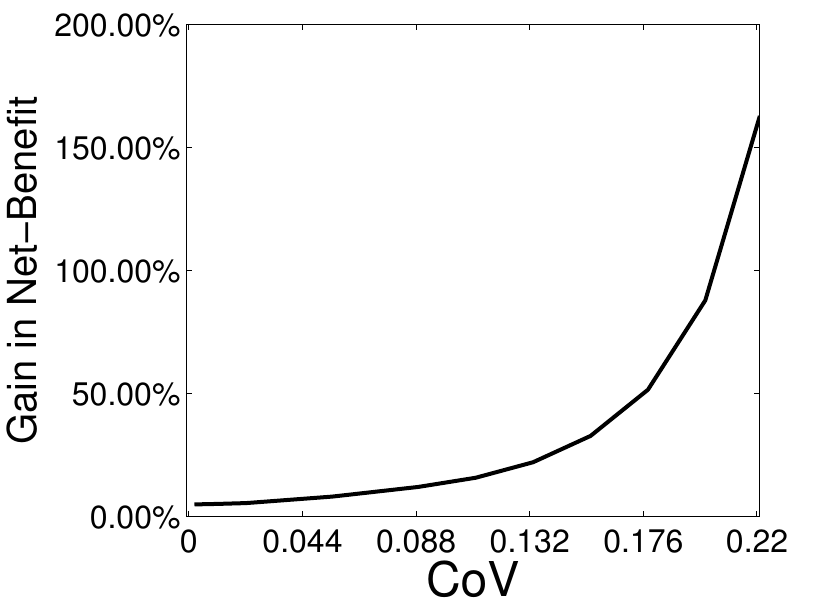}
\hspace*{0.0005in}
\includegraphics[width=0.49\textwidth]{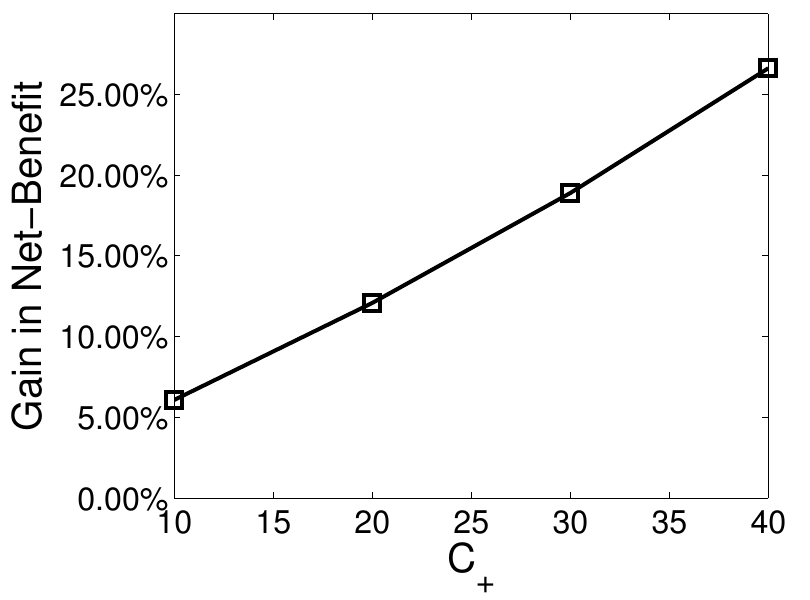}
\caption{Improvement in expected net-benefit under our optimal dynamic control policy relative to
the alternative offline LP (top two plots) and online CF algorithms (bottom two plots) for the first
set of workloads based on $f^1(t)$ and for varying values of $\sigma$ and $\mathcal{C}_{+}$.}
\label{fig:results1}
\end{figure}

Figure~\ref{fig:results1} presents a representative sample of our computational results for
the first demand process based on $f^1(t)$.
The two leftmost graphs provide performance comparisons of our optimal online dynamic control
algorithm against the alternative offline LP and online CF algorithms, respectively,
where both comparisons are based on the relative improvements in expected net-benefit under our
optimal control policy as a function of $\sigma$;
the relative improvement is defined as the difference in expected net-benefit under our optimal
dynamic control policy and under the alternative optimization approach, divided by the expected
net-benefit of the alternative approach.
For the purpose of comparison across sets of workloads with very different $f_{\mbox{\tiny avg}}$
values, we plot both of these graphs as a function of the coefficient of variation
$\mbox{\textsf{CoV}} = \sigma / f_{\mbox{\tiny avg}}$.
The two rightmost graphs provide similar comparisons of relative improvement in expected net-benefit
between our optimal dynamic control policy and the two alternative optimization approaches as a
function of $\mathcal{C}_{+}$, both with $\sigma$ fixed to be $0.4$.

\begin{figure}[ht!]
\centering
\includegraphics[width=0.7\textwidth]{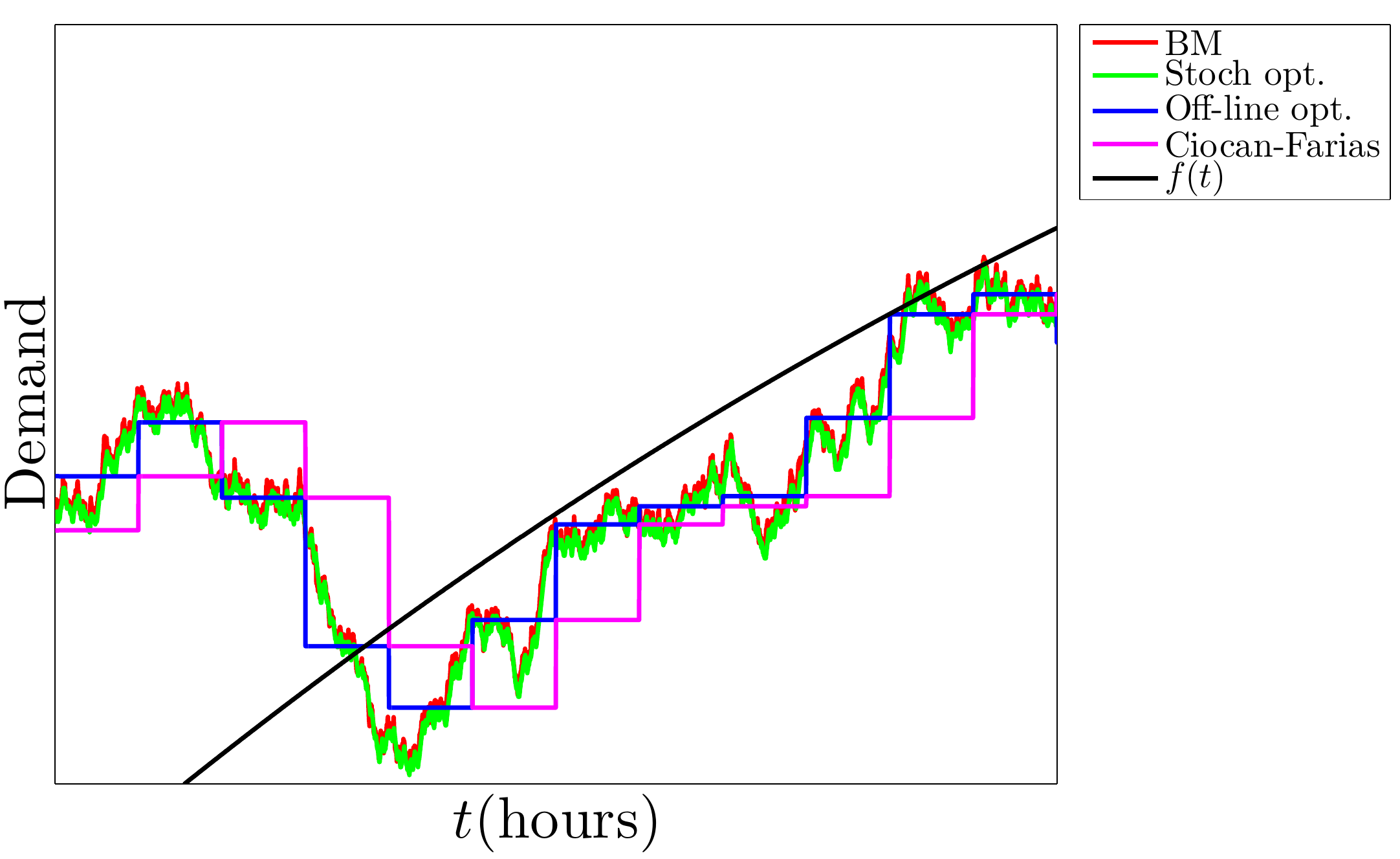}
\caption{Performance of all three algorithms over a representative interval of a single sample path.}
\label{fig:results1:sp}
\end{figure}

We first observe from the two leftmost graphs in Figure~\ref{fig:results1} that our optimal online
dynamic control algorithm outperforms the two alternative optimization approaches for all $\sigma > 0$.
The relative improvements in expected net-benefit grow in an
exponential manner with respect to increasing values of $\sigma$ over the range of $\mbox{\textsf{CoV}}$
values $(0,0.22]$ considered, with relative improvements up to $90\%$ in comparison with the offline LP algorithm
and more than $150\%$ in comparison with the CF algorithm.
Our results illustrate and quantify the fact that, even in discrete-time models with small time slot
lengths $\gamma$, nonnegligible volatility plays a critical role in the expected net-benefit of
any given resource allocation policy.
The significant relative improvements under the optimal online dynamic control algorithm then follow
from our stochastic optimal control approach that directly addresses the volatility of the demand process
in all primary and secondary resource allocation decisions.
This can be clearly seen in the results of Figure~\ref{fig:results1:sp} that illustrate the performance
of the three algorithms relative to demand over a representative interval of an individual sample path.
Figure~\ref{fig:results1:sp} represents a zoomed-in view of the results over a small segment of the time horizon (60 minutes of the 24 hour time horizon).
Although the offline LP algorithm based on \eqref{offline:obj},\eqref{offline:st1} would eventually
outperform our optimal online dynamic control algorithm as the time slot length $\gamma$ decreases
toward $0$, we note that the choice for $\gamma$ in our computational experiments is considerably
smaller than the $10$-minute intervals suggested in~\cite{LiWiAn+11}.
Moreover, as discussed in~\cite{CioFar12}, the optimal choice of $\gamma$ is a complex issue in and
of itself and it may need to vary over time depending upon the statistical properties of the demand
process $D(t)$.
A key advantage of our optimal online dynamic control algorithm is that such parameters are not needed.

We next observe from the two rightmost graphs in Figure~\ref{fig:results1} that the relative improvements
in expected net-benefit under our optimal online dynamic control algorithm similarly increases with respect
to increasing values of $\mathcal{C}_+$, though in a more linear fashion.
We also note that very similar trends were observed with respect to varying the value of $\mathcal{C}_{-}$,
though the magnitude of the relative improvement in expected net-benefit is smaller.
Within the limited range of parameter values considered,
our computational experiments suggest that the relative improvements in net-benefit under our optimal
dynamic control policy can be more sensitive to $\mathcal{C}_{+}$ than to $\mathcal{C}_{-}$.
Recall that $\mathcal{C}_+ = \mathcal{C}_p$ is the cost for the primary resource allocation capacity,
whereas $\mathcal{C}_- = \mathcal{N}_p - \mathcal{N}_s$ is the difference in net-benefit between the
primary and secondary resource allocation capacities.
%As noted above, many other model parameter values were varied to investigate the impact and sensitivity
%of these parameters on the relative improvements of our optimal dynamic control policy over the other
%two alternative optimization algorithms.
%We omit these results due to space restrictions.
%However, we
We also
note that similar trends were observed for changes in the values of
$f^1_{\mbox{\tiny min}}$ and $f^1_{\mbox{\tiny max}}$ when the relative improvement results are
considered as a function of $\mbox{\textsf{CoV}}$.

Now let us turn to our second set of results based on the average daily demand pattern $f^2(t)$
illustrated in the rightmost plot of Figure~\ref{fig:fx}, where the base parameter settings are given by
$\alpha=0.02$, $\sigma=7.0$, $\theta_{l}=-100$, $\theta_{u}=100$, $\mathcal{C}_{+}=20$, $\mathcal{C}_{-}=2$,
$\mathcal{D}_{p}=0.5$, $\mathcal{I}_{p}=0.5$, $f^2_{\mbox{\tiny min}}=15$, $f^2_{\mbox{\tiny max}}=90$,
$f^2_{\mbox{\tiny avg}}=61$, $x = X(0) = P(0) - D(0) = 0$.
In addition to these base settings, we vary the parameter values one at a time for $\sigma \in [0.01, 15]$,
$\mathcal{C}_{+} \in [10,40]$, $\mathcal{C}_{-} \in [1,10]$, $f^2_{\mbox{\tiny min}} \in [1,20]$ and
$f^2_{\mbox{\tiny max}} \in [9,120]$.
Once again, for each experiment comprising a specific workload, we generate $N=10,000$ sample paths
using a timescale of a couple of seconds and a $\gamma$ setting of five minutes, noting that a wide
variety of experiments with different timescale and $\gamma$ settings provided performance trends
that are identical to those presented herein.
We then apply our optimal dynamic control policy and the two alternative optimization approaches
to this set of $N$ sample paths as described above.
Our performance evaluation comparisons are based on the expectation of net-benefit realized under
each of the three algorithms, also as described above.
\begin{figure}[ht!]
\centering
\includegraphics[width=0.49\textwidth]{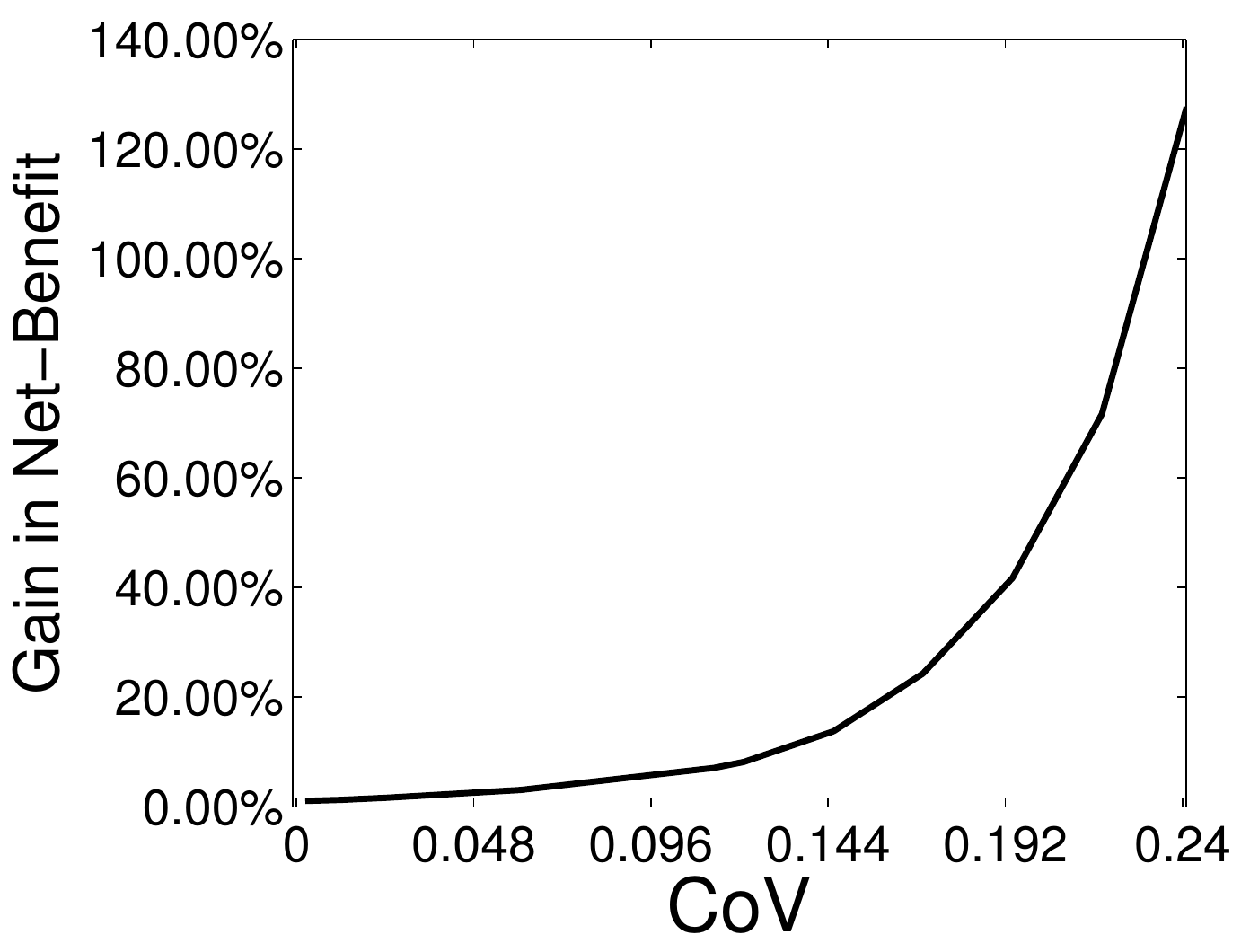}
\hspace*{0.0005in}
\includegraphics[width=0.49\textwidth]{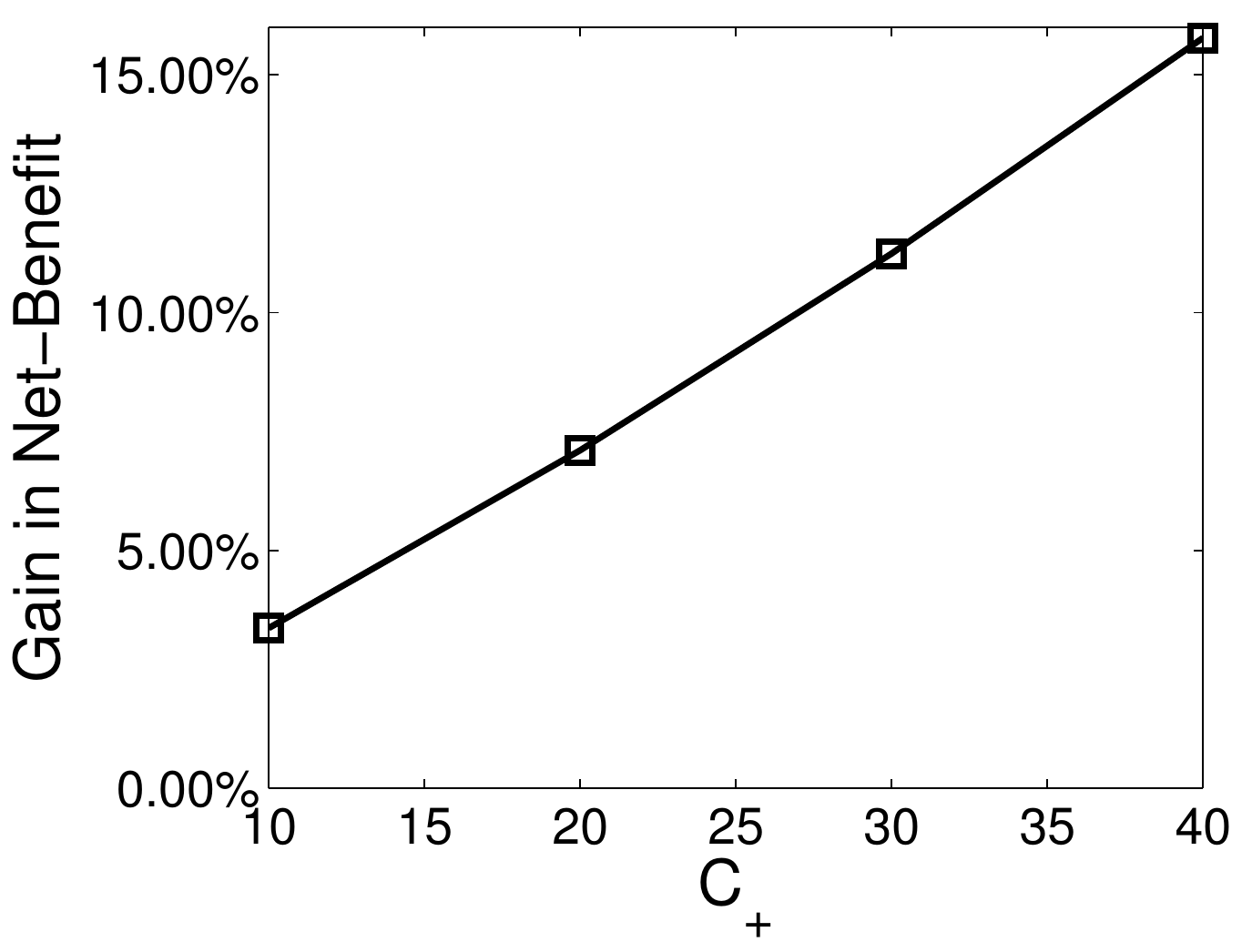}
\hspace*{0.0005in}
\includegraphics[width=0.49\textwidth]{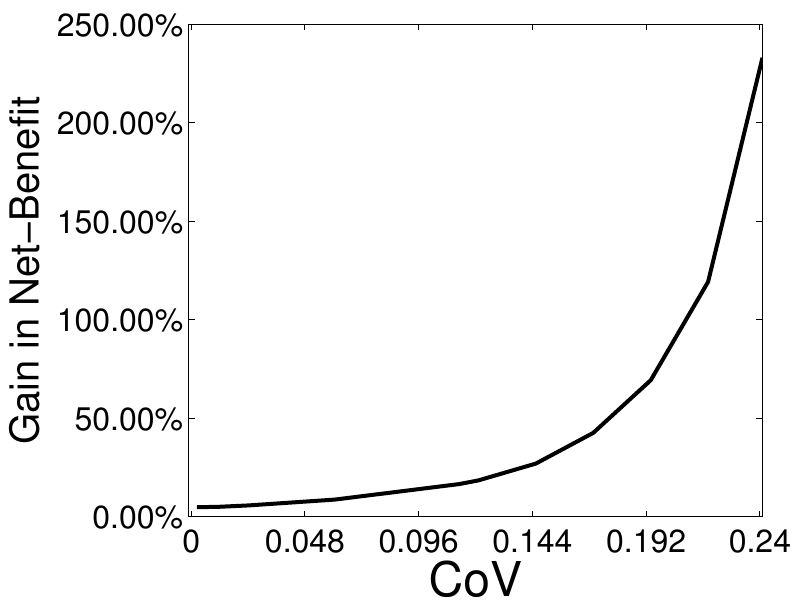}
\hspace*{0.0005in}
\includegraphics[width=0.49\textwidth]{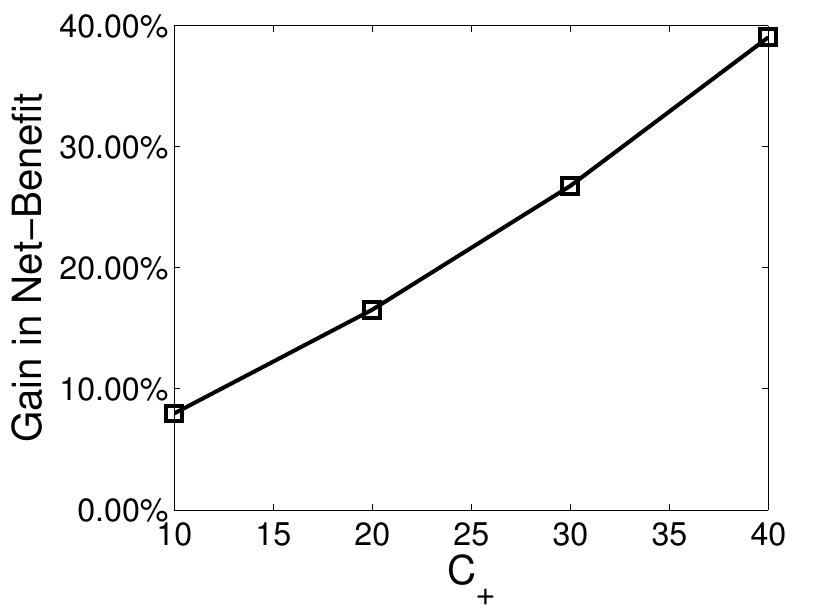}
\caption{Improvement in expected net-benefit under our optimal dynamic control policy relative to
the alternative offline LP (top two plots) and online CF algorithms (bottom two plots) for the second set of workloads based
on $f^2(t)$ and for varying values of $\sigma$ and $\mathcal{C}_{+}$.}
\label{fig:results2}
\end{figure}

Figure~\ref{fig:results2} presents a representative sample of our computational results for the demand process
based on $f^2(t)$, providing the analogous results that correspond
to those in Figure~\ref{fig:results1}.
We note that the larger range $[f^2_{\mbox{\tiny min}},f^2_{\mbox{\tiny max}}]$ exhibited in the
second average daily demand pattern as well as a higher value of $f^2_{\mbox{\tiny avg}}$ lead to
both a higher relative net-benefit for fixed $\sigma$ and a higher sensitivity to changes in $\sigma$,
thus improving the gains of our optimal online dynamic control algorithm over the alternative offline
LP and online CF algorithms.
This relative improvement in expected net-benefit as compared to the set of experiments for the
average daily demand pattern $f^1(t)$ can be understood to be caused in part by the sharp drop in average
demand from the maximum value of $90$ to a minimum of $15$ within a fairly short time span,
thus contributing to an increased effective volatility over and above that represented by $\sigma$.
Hence, the fact that the relative improvement exhibited by our optimal online dynamic control
algorithm is larger under the average daily demand pattern $f^2(t)$, up to $130\%$ in comparison
with the offline LP algorithm and over $230\%$ in comparison with the CF algorithm, can be viewed
in a sense to be very much an extension of the finding that the relative improvement provided by
our optimal online algorithm increases with an increase in $\mbox{\textsf{CoV}}$.
A similar gain in performance improvement can be seen in the rightmost two plots when we vary $\mathcal{C}_+$
with a fixed value of $\sigma$.

%In comparing the results in Figures~\ref{fig:results1} and \ref{fig:results2} as well as those from
%our many other computational experiments,  we observe that all of these net-benefit results are quite
%consistent in terms of qualitative and quantitative trends, where our optimal online dynamic control
%algorithm consistently provides superior results even when the demand uncertainty/volatility is very
%small.
%Our computational experiments further indicate that the primary resource allocation capacity decisions
%under the online CF algorithm are identical to those under the offline LP algorithm shifted by one
%time slot.
%Finally, we note that, upon applying our optimal online dynamic control algorithm and the offline
%LP and online CF algorithms directly to the actual trace data, the same performance trends among
%the three optimization approaches were observed as those shown above for the time-average over
%many sample paths from the Brownian motion demand process fitted to the same trace data.

\subsection{Multiple Primary Resources}
We now turn to investigate the relationship between the single primary resource formulation (SC-OPT:S) and the multiple primary resource formulation (SC-OPT:M).
One can envision, under appropriate circumstances and contractual agreements among the multiple primary resource options,
that the multiple primary resource allocation model can potentially render greater performance benefits to an organization
than the corresponding performance benefits from the single primary resource allocation model.
In this section we consider the performance trade-offs between both primary resource allocation models.

To this end, we present a comparison between the single primary resource formulation (SC-OPT:S)
and the multiple primary resource formulation (SC-OPT:M) under a two-dimensional contract vector $w$.
The particular representative set of results presented here is based on the average daily demand pattern $f^2(t)$ illustrated in the rightmost plot of Figure~\ref{fig:fx},
with the base parameter settings given by $\alpha=0.02$,
$\sigma=1$,
$\theta_{l,1}=-1000$, $\theta_{l,2}=-1000$,
$\theta_{u,1}=1000$, $\theta_{u,2}=1000$, 
$\mathcal{N}_{p,1} = 10$, $\mathcal{N}_{p,2} = 10$,
$\mathcal{C}_{p,1} = 2000$, $\mathcal{C}_{p,2} = 1200$,
$\mathcal{N}_s = 1$,
%$\mathcal{C}_{+}(w)=2000$, $\mathcal{C}_{-}(w)=9$,
$\mathcal{D}_{p,1}=0.001$, $\mathcal{D}_{p,2}=0.001$,
$\mathcal{I}_{p,1}=0.001$, $\mathcal{I}_{p,2}=0.001$,
$f^2_{\mbox{\tiny min}}=50$, $f^2_{\mbox{\tiny max}}=125$,
$f^2_{\mbox{\tiny avg}}=96.9$ and $x = X(0) = P(0) - D(0) = 0$, where
$f_{\mbox{\tiny min}} := \min_t\{f(t)\}$, $f_{\mbox{\tiny max}} := \max_t\{f(t)\}$
and
$f_{\mbox{\tiny avg}} := T^{-1} \int_0^T f(t)dt$.
In addition to these base settings, we vary the parameter values one at a time for
$\sigma \in [1.0, 5.0]$,
$\mathcal{N}_{p,1} \in [4,14]$, $\mathcal{N}_{p,2} \in [4,14]$,
$\mathcal{C}_{p,1} \in [1700,2300]$, $\mathcal{C}_{p,2} \in [1000,1400]$,
$f^2_{\mbox{\tiny min}} \in [40,60]$ and $f^2_{\mbox{\tiny max}} \in [115,135]$,
in order to investigate the impact and sensitivity of these parameters
on the performance trade-offs between the two formulations.
Within this experimental setting, as a representative example, we consider the corresponding single primary resource allocation model with $w_{s}=\begin{bmatrix}1&0\end{bmatrix}$
and consider a corresponding instance of the multiple primary resource allocation model with $w_{m}=\begin{bmatrix}0.7&0.3\end{bmatrix}$.
The values of $L(w)$ and $U(w)$ are separately obtained for each of the two primary resource formulations.
For every computational experiment under a given set of parameters, we generate $N=10,000$ daily sample paths using a timescale of a couple of seconds
and a $\gamma$ setting of five minutes.
We note that a wide variety of experiments with different parameter, timescale and $\gamma$ settings
were evaluated and shown to exhibit similar performance trends as those presented herein.
We then apply our optimal dynamic control policy for each formulation to this set of $N$ daily sample paths as described in Section~\ref{sec:exp:single}.

\begin{figure}[ht!]
\centering
\includegraphics[width=0.49\textwidth]{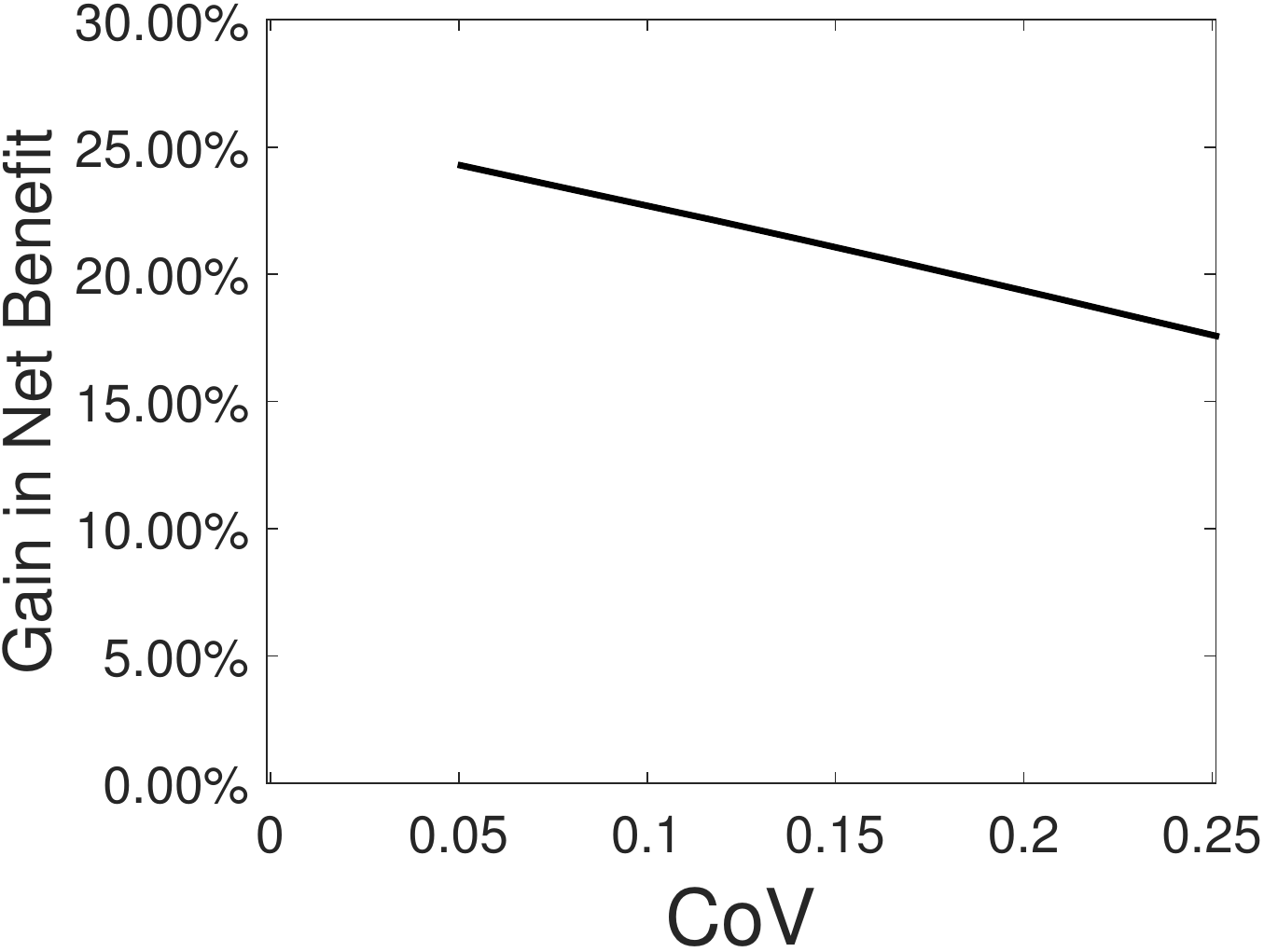}
\hspace*{0.0005in}
\includegraphics[width=0.49\textwidth]{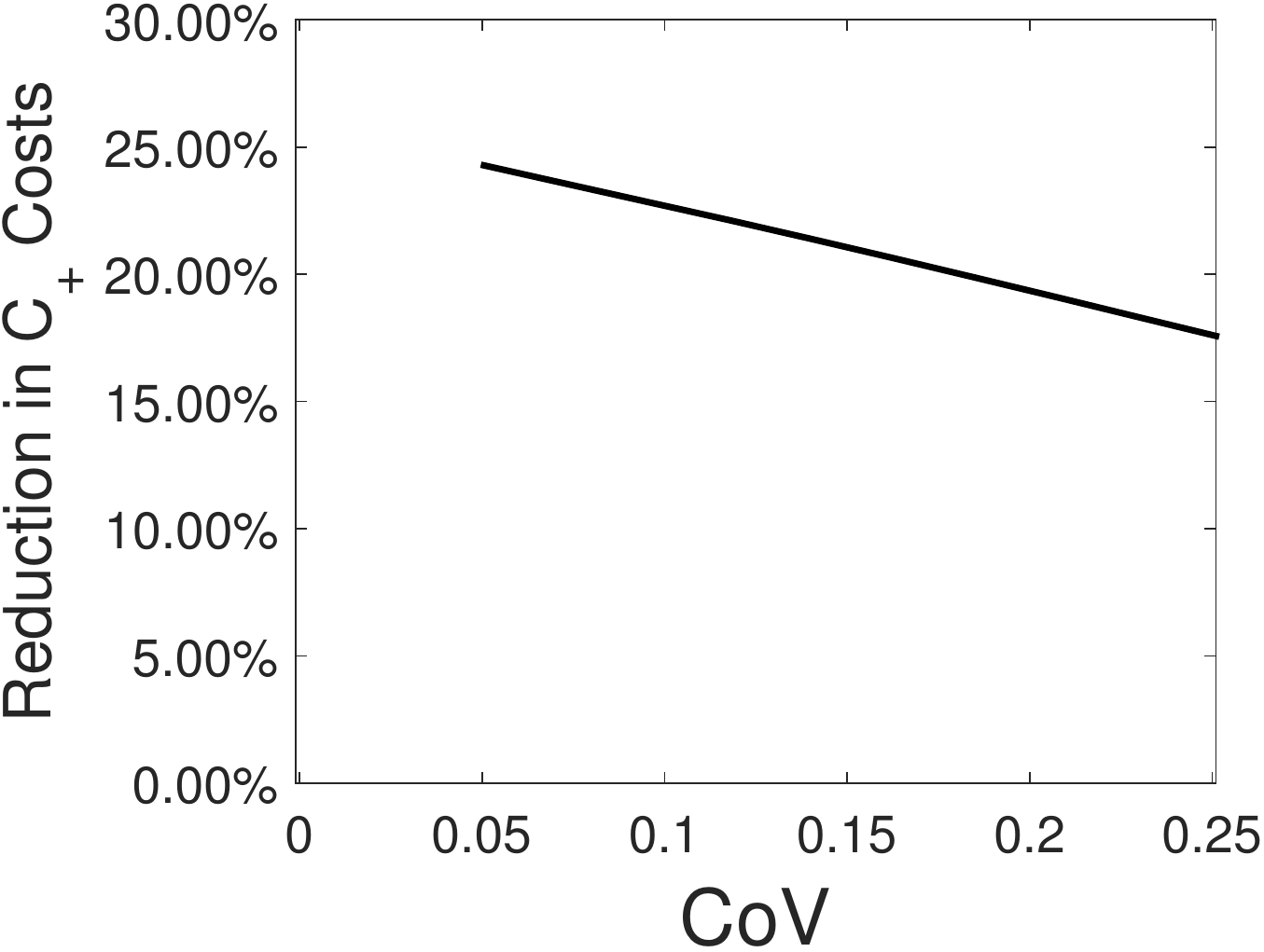}
\hspace*{0.0005in}
\caption{Improvements in expected net-benefit and expected cumulative discounted $\mathcal{C}_{+}(w)$ cost under our optimal dynamic control policy
for the multiple primary allocation model relative to that for the single primary allocation model based on $f^2(t)$ while increasing the demand uncertainty $\sigma$.}
\label{fig:results3}
\end{figure}

Figure~\ref{fig:results3} presents a representative sample of our computational results in which the leftmost graph provides comparisons based on
the relative improvements in expected net-benefit under our optimal control policy for the two primary resource formulations as a function of the
coefficient of variation $\mbox{\textsf{CoV}} = \sigma / f_{\mbox{\tiny avg}}$;
analogous to the leftmost graphs in Figures~\ref{fig:results1} and \ref{fig:results2},
the relative improvement is defined here as the difference between the expected net-benefit for the multiple primary resource model and the single primary resource model,
divided by the expected net-benefit for the single primary resource model.
The rightmost graph provides comparisons based on the relative improvements in expected cumulative discounted costs associated with $\mathcal{C}_{+}(w)$ 
under our optimal control policy for the two primary resource formulations as a function of $\mbox{\textsf{CoV}}$;
the relative improvement is defined here as the difference between the expected cumulative discounted contributions of $\mathcal{C}_{+}(w)$ over the infinite horizon
for the single primary resource model (i.e., in \eqref{optHF:obj}) and the multiple primary resource model (i.e., in \eqref{optHF:obj_w}),
divided by the expected net-benefit for the single primary resource model.
From these results we observe that the gain in expected net-benefit is significant, demonstrating the potential performance benefits to an organization under
the multiple primary resource allocation model for such problem instances.
These expected net-benefit improvements in the leftmost graph tend to decrease as the coefficient of variation increases,
while still remaining significant over the range of $\mbox{\textsf{CoV}}$ values.
This decrease in the expected net-benefit gain is primarily due to a similar decreasing trend in the relative difference in the expected cumulative discounted contributions
of $\mathcal{C}_{+}(w)$ over the infinite horizon as $\mbox{\textsf{CoV}}$ increases.
To help explain this, we note that the values of the second summand in \eqref{eq:simplify} and the expected cumulative discounted contributions of $\mathcal{C}_{+}(w)$
over the infinite horizon are of the same order of magnitude in these experiments, whereas the values of the remaining terms in \eqref{optHF:obj} and \eqref{optHF:obj_w}
are orders of magnitude smaller;
the similarity of the two graphs are due to the facts that the two higher order magnitude terms dominate the objective function value (expected net-benefit)
and the value of the second summand in \eqref{eq:simplify} is the same under both primary resource models.
Hence, the trends in the relative expected net-benefit improvements as a function of $\mbox{\textsf{CoV}}$ are directly related to the very similar trends in the
relative expected discounted cumulative $\mathcal{C}_{+}(w)$ cost improvements as a function of $\mbox{\textsf{CoV}}$.
These trends in turn are primarily due to the role of the risk-hedging interval that widens under each model as a function of the increasing $\mbox{\textsf{CoV}}$
in order to have the optimal dynamic control policy reduce the expected discounted cumulative costs associated with the primary resource(s) over the infinite horizon.
In other words, the optimal dynamic control policy under both models becomes somewhat more conservative due to the greater risks associated with having a primary
resource allocation position that is too large and that would otherwise result in larger expected discounted cumulative costs associated with $\mathcal{C}_{+}(w)$.
These factors have a somewhat stronger impact on the multiple primary resource model as $\mbox{\textsf{CoV}}$ increases.

Lastly, it can also be important to consider the role of the margin of rewards and costs when investigating the relationship between
the single primary resource formulation (SC-OPT:S) and the multiple primary resource formulation (SC-OPT:M).
In many practical applications, the overall rewards and costs are of similar magnitude with reasonably balanced margins,
and thus the relationship between the two formulations includes key trade-offs among the net-benefits from both the second summand in \eqref{eq:simplify}
and the term involving $\mathcal{C}_{-}$ and $\mathcal{C}_{-}(w)$ in \eqref{optHF:obj} and \eqref{optHF:obj_w} on the one hand,
and the relative costs from the terms involving $\mathcal{I}_p$, $\mathcal{D}_p$ and $\mathcal{C}_{+}$ in \eqref{optHF:obj} and
involving $\mathcal{I}_p(w)$, $\mathcal{D}_p(w)$ and $\mathcal{C}_{+}(w)$ in \eqref{optHF:obj_w} on the other hand.
These trade-offs are indeed reflected in the representative results illustrated in Figure~\ref{fig:results3}.
However, in situations where the net-benefits are considerably larger than the relative costs in \eqref{optHF:obj} and \eqref{optHF:obj_w},
then the relationship between the single primary resource formulation and the multiple primary resource formulation can depend in large part on the magnitude
and ordering between $\mathcal{N}_p$ and $\mathcal{N}_p(w)$;
in this case with $\mathcal{N}_p$ dominating $\mathcal{N}_p(w)$, the single primary resource model can outperform the multiple primary resource model .
Analogously, when the relative costs are considerably larger than the net-benefits in \eqref{eq:simplify}, \eqref{optHF:obj} and \eqref{optHF:obj_w},
then the relationship between the two formulations can depend in large part on the magnitude and ordering among
$\mathcal{C}_{+}$, $\mathcal{I}_p$, $\mathcal{D}_p$, $\mathcal{C}_{+}(w)$, $\mathcal{I}_p(w)$ and $\mathcal{D}_p(w)$;
in this case with, e.g., $\mathcal{C}_{+}$ dominating $\mathcal{C}_{+}(w)$, the multiple primary resource model can outperform the single primary resource model.

\section{Conclusions}\label{sec:conclusion}
In this paper we investigated a general class of dynamic resource allocation problems
arising across a broad spectrum of application domains that intrinsically involve
different types of resources and uncertain/variable demand.
With a goal of maximizing expected net-benefit based on rewards and costs from the
different resources, we derived the provably optimal dynamic control policy within a stochastic optimal control setting.
Our mathematical analysis includes obtaining simple expressions that govern the
dynamic adjustments to resource allocation capacities over time under the optimal
control policy.
A wide variety of extensive computational experiments demonstrates and quantifies the
significant benefits of our optimal dynamic control policy over recently proposed
alternative optimization approaches in addressing a general class of resource
allocation problems across a diverse range of application domains, including cloud
computing and data center environments, computer and communication networks, and
human capital supply chains.
Moreover, our results strongly suggest that the stochastic optimal control approach
taken in this paper can provide an effective means to develop easily-implementable
online algorithms for solving stochastic optimization problems.
Both single primary resource allocation and multiple primary resource allocation models
can be exploited, with the best option depending upon the system environment and model parameters.

Following along the lines of our computational experiments, our algorithm can exploit any
consistent seasonal patterns for $b(t)$ and $\sigma(t)$ observed from historical traces
in order to predetermine the threshold values $L$ and $U$ or $L(w)$ and $U(w)$.
In addition, various approaches such as statistical learning and/or model predictive control 
(e.g.,~\cite{CioFar12}) can be used to adjust these threshold values in real-time based on
identifying and learning any nonnegligible changes in the realized values for $b(t)$ and $\sigma(t)$.
Furthermore, this latter approach can be used directly for system/network environments whose
demand processes do not exhibit consistent seasonal patterns.

\section*{Acknowledgments.}
Xuefeng Gao acknowledges support from Hong Kong RGC ECS Grant 2191081 and CUHK Direct Grants for Research with project codes 4055035 and 4055054. A portion of Mark Squillante's research was sponsored by the US Army Research Laboratory and the UK Ministry of Defence and was accomplished under Agreement Number W911NF-06-3-0001. The views and conclusions contained in this document are those of the authors and should not be interpreted as representing the official policies, either expressed or implied, of the US Army Research Laboratory, the US Government, the UK Ministry of Defense, or the UK Government. The US and UK Governments are authorized to reproduce and distribute reprints for Government purposes notwithstanding any copyright notation hereon.
%The contribution from Joost Bosman has been carried out in the context of the IOP GenCom project Service Optimization and Quality (SeQual), supported by the Dutch Ministry of Economic Affairs, Agriculture and Innovation via its agency Agentschap NL.

%%%%%%%%%%%%%%%%%%%%%%%%%%%%%%%%%%%%%%%%%%%%%%%%%%%%%%%%%%%%%%%%%%%%%%%
%\newpage
\vspace*{1.0in}

\appendix
\section{Proofs of Results in Section~\ref{sec:main}} \label{sec:proofs}
In this appendix we collect the proofs of our main results in Section~\ref{sec:main}.
We start with a rigorous proof of Theorem~\ref{THM:CASE1} and Corollary~\ref{THM:SPECIAL}, which proceeds
in three main steps.
First, we express the optimality conditions for the stochastic dynamic
program, i.e., the Bellman equation corresponding to
\eqref{optHF:obj}~--~\eqref{optHF:st3}.
We then derive a solution of the Bellman equation and determine the
corresponding candidate value function and dynamic control policy,
establishing smoothness and convexity of the candidate value function and uniqueness
of the threshold values.
Finally, we verify that this dynamic control policy is indeed optimal
through a martingale argument.
Each of these main steps is presented in turn.
Then we present the proofs of our other main results.

\subsection{Proof of Theorem~\ref{THM:CASE1} and Corollary~\ref{THM:SPECIAL}: Step 1}\label{sec:proofs:case1:step1}
From the Bellman principle of optimality, we deduce that the
value function $V$ satisfies for each $t \ge 0$
\begin{align} \label{eq:prinofopt}
V(x) &= \min_{\theta_\ell \le \dot P(t)\le \theta_u} \E_x \bigg[ \int_0^t e^{-\alpha s} \Big\{\Big( \mathcal{C}_+ X(s)^+ + \mathcal{C}_- X(s)^{-}\Big) ds
+ ( \mathcal{I}_p \indi{\dot P(s)>0}-\mathcal{D}_p \indi{\dot P(s)<0} ) dP(s) \Big\} \nonumber \\
& \qquad\qquad\qquad\qquad + \; e^{-\alpha t} V(X(t))\bigg] ;
\end{align}
refer to~\cite[Chapter 4]{YonZho99}.
Suppose the value function $V$ is smooth, belonging to the set $C^2$ (i.e.,
the set of twice continuously differentiable functions) except for a finite
number of points, and $V'(x)$ is bounded for any $x$.
%Then we can apply Ito's formula (see, e.g.,~\cite{KarShr91}) to
%$e^{-\alpha t} V(X(t))$ and derive
%\begin{eqnarray}
%\lefteqn{e^{-\alpha t} V(X(t))- V(x)+ \int_{0}^{t} e^{-\alpha s}V'(X(s))
%    \sigma dw(s) = \int_{0}^{t} e^{-\alpha s}} \nonumber\\
%& [-\alpha V(X(s)) + \frac{1}{2} \sigma^2 V''(X(s)) +
%    ( \dot P(s)-b) V'(X(s)) ] ds . \qquad \label{eq:step1a}
%\end{eqnarray}
%
%Upon taking expectation on both sides of \eqref{eq:step1a} in combination
%with \eqref{eq:prinofopt}, we obtain
%\begin{align}
%0 = & \min_{\theta_\ell \le \dot P(t)\le \theta_u} \int_{0}^{t} e^{-\alpha s}
%    [-\alpha V(X(s)) + \frac{1}{2} \sigma^2 V''(X(s)) \nonumber \\
%& + (\theta-b) V'(X(t)) + \mathcal{C}_+ X(s)^+ + \mathcal{C}_- X(s)^-
%    \nonumber \\
%& + \mathcal{I}_p \dot P(s) \indi{\dot P(s)>0}-\mathcal{D}_p \dot P(s)
%    \indi{\dot P(s)<0} ] ds .
%\end{align}
%Dividing both sides by $t$ and letting $t$ tend to 0, we deduce the
%desired Bellman equation
Then, based on a standard application of Ito's formula as in~\cite[Chapter 1]{Kryl80},
we derive that the desired Bellman equation for the value function $V$
has the form
\begin{equation} \label{eq:bellmaneq}
-\alpha V(x) + \frac{1}{2} \sigma^2 V''(x) -b V'(x) + \mathcal{C}_+ x^+ + \mathcal{C}_- x^{-}
+ \inf_{ \theta_\ell \le \theta \le \theta_u} \cL(\theta, x) = 0,
\end{equation}
where
\begin{equation}
\cL(\theta, x)= \left\{\begin{matrix}
(V'(x)+\mathcal{I}_p) \theta \qquad \text{if} \quad \theta \ge 0 ,\\
(V'(x)-\mathcal{D}_p) \theta \qquad \text{if} \quad \theta<0 .
\end{matrix}\right. \label{eq:bellmaneq2}
\end{equation}

\subsection{Proof of Theorem~\ref{THM:CASE1} and Corollary~\ref{THM:SPECIAL}: Step 2}\label{sec:proofs:case1:step2}
Our next goal is to construct a convex function $Y$ that satisfies the
Bellman equation \eqref{eq:bellmaneq} and show that the threshold values
$L$ and $U$ are uniquely determined by the corresponding pair of nonlinear
equations in Theorem~\ref{THM:CASE1}.
Suppose a candidate value function $Y(x)$ satisfies \eqref{eq:bellmaneq}. Given \eqref{eq:bellmaneq2}, we expect a ``bang-bang'' type solution based on the signs of $Y'(x)+\mathcal{I}_p$ and $Y'(x)-\mathcal{D}_p$. In particular,
we seek to find $L$ and $U$ such that
\begin{equation} \label{eq:derofY}
Y'(x)= \left\{\begin{array}{ll}
 \ge \mathcal{D}_p  , & \qquad \text{if} \quad x \ge U ,\\
 \in (-\mathcal{I}_p, \mathcal{D}_p) & \qquad \text{if} \quad L<x < U ,\\
\le -\mathcal{I}_p & \qquad \text{if} \quad x \le L .
\end{array}\right.
\end{equation}
Moreover, we require that $Y$ meets
smoothly at the points $L$, $0$ and $U$ to order one, and $Y(x)=O(|x|)$ as $|x| \rightarrow \infty$ (i.e., $\lim_{|x| \rightarrow \infty} \frac{Y(x)}{|x|} \le C$ for some $C \ge 0$) so that $Y'$ is bounded.

For each of the three cases in Theorem~\ref{THM:CASE1}, we first solve the Bellman Equation \eqref{eq:bellmaneq} to derive the
corresponding pair of equations that $L$ and $U$ satisfy. Then we discuss conditions on model parameters under which the points
$L$ and $U$ are located in comparison with $0$. Finally, we show that the function $Y$ we construct has the property \eqref{eq:derofY}.

\subsubsection{Case I: $U \ge L \ge 0$}\label{sec:521}
%
%Focusing on the case $U > L \ge 0$ where $L$ and $U$ satisfy
%\eqref{eq:xhminusxf} and \eqref{eq:xH} subject to \eqref{eq:derofY},
We proceed to solve the Bellman equation \eqref{eq:bellmaneq}
depending on the value of $x$ in relation to $U$, $0$ and $L$.
There are four subcases to consider as follows.

\textit{(i).} ~If $x \ge U > 0$, we obtain from \eqref{eq:derofY} that
$$Y'(x) \; \ge \; \mathcal{D}_p \qquad \mbox{and} \qquad \inf_{ \theta_\ell \le \theta \le \theta_u} \cL(\theta, x) \; = \; \cL(\theta_\ell, x).$$
Then the Bellman equation \eqref{eq:bellmaneq} yields
$$
-\alpha Y(x) + \frac{1}{2} \sigma^2 Y''(x) -b Y'(x) + \mathcal{C}_+ x^+ + \mathcal{C}_- x^{-}+ \cL(\theta_\ell, x) \; = \; 0 ,
$$
or equivalently
$$
-\alpha Y(x) + \frac{1}{2} \sigma^2 Y''(x) -b Y'(x) +\mathcal{C}_+ x + (Y'(x)-\mathcal{D}_p)\theta_\ell \; = \; 0.
$$
Solving this second order linear nonhomogeneous differential equation, we obtain
\begin{equation*}
Y(x) \; = \; \frac{\mathcal{C}_+}{\alpha} x + \frac{1}{\alpha} \left(\frac{\mathcal{C}_+}{\alpha}(\theta_\ell-b) - \mathcal{D}_p \theta_\ell\right)+ l_1 e^{t_2x} + \bar l_1 e^{t_1 x},
\end{equation*}
where $t_1>0>t_2$ are given in (\ref{eq:t2}) and $l_1, \bar l_1$ are two constants to be determined. Since $Y(x)=O(|x|)$ when $|x|$ goes to $\infty$, one finds $\bar l_1 =0.$ Thus we derive for $x \ge U$
\begin{equation} \label{eq:caseI:YlargerU}
Y(x) \; = \; \frac{\mathcal{C}_+}{\alpha} x + \frac{1}{\alpha} \left(\frac{\mathcal{C}_+}{\alpha}(\theta_\ell-b) - \mathcal{D}_p \theta_\ell\right)+ l_1 e^{t_2x},
\end{equation}

\textit{(ii).} ~If $0 \le L < x < U$, we have
$$-\mathcal{I}_p<Y'(x) \; < \; \mathcal{D}_p, \qquad \mbox{and} \qquad \inf_{ \theta_\ell \le \theta \le \theta_u} \cL(\theta, x) \; = \; 0,  $$
and thus we obtain
$$
-\alpha Y(x) + \frac{1}{2} \sigma^2 Y''(x) -b Y'(x) +\mathcal{C}_+ x \; = \; 0.
$$
This implies for $0 \le L < x<U$
\begin{equation} \label{eq:caseI:YLU}
Y(x) \; = \; \frac{\mathcal{C}_+}{\alpha} x + \frac{-b \mathcal{C}_+}{\alpha^2}+ \lambda_1 e^{r_1 x} + \lambda_2 e^{r_2 x},
\end{equation}
where $r_1, r_2$ are given by (\ref{eq:r2}), and  $\lambda_1, \lambda_2$ are two generic constants to be determined.

\textit{(iii).} ~If $0 <x \le L$, we find from \eqref{eq:derofY} that
$$Y'(x) \; \le \; -\mathcal{I}_p \qquad \mbox{and} \qquad \inf_{ \theta_\ell \le \theta \le \theta_u} \cL(\theta, x) \; = \; \cL(\theta_u, x). $$
Then \eqref{eq:bellmaneq} becomes
$$
-\alpha Y(x) + \frac{1}{2} \sigma^2 Y''(x) -b Y'(x) + \mathcal{C}_+ x + (\mathcal{I}_p + Y'(x)) \theta_u \; = \; 0,
$$
and thus
\begin{equation}\label{eq:caseI:Y0L}
Y(x) \; = \; \frac{\mathcal{C}_+}{\alpha} x  + \frac{1}{\alpha} \left(\frac{\mathcal{C}_+}{\alpha}(\theta_u-b) + \mathcal{I}_p \theta_u\right) + \tilde \lambda_1 e^{s_1 x} + \tilde \lambda_2 e^{s_2 x},
\end{equation}
where $s_1, s_2$ are given by (\ref{eq:s2}), and $\tilde \lambda_1, \tilde \lambda_2$ are two generic constants to be
determined.

\textit{(iv).} ~If $x \le 0$, we have
$$Y'(x) \; \le \; -\mathcal{I}_p,$$
and we similarly derive that
\begin{eqnarray*}
-\alpha Y(x) + \frac{1}{2} \sigma^2 Y''(x) -b Y'(x) - \mathcal{C}_- x + (Y'(x)+\mathcal{I}_p) \theta_u \; = \; 0.
\end{eqnarray*}
Since $Y(x)=O(|x|)$ as $|x|\rightarrow \infty$, we deduce that the solution is then given by
\begin{equation} \label{eq:caseI:YlessL}
Y(x) \; = \; -\frac{\mathcal{C}_-}{\alpha} x + \frac{1}{\alpha} \left(-\frac{\mathcal{C}_-}{\alpha}(\theta_u-b) +\mathcal{I}_p \theta_u\right)+ l_2 e^{s_1 x},
\end{equation}
where $s_1$ is given by (\ref{eq:s2}), and $l_2$ is a generic constant to be determined.

Now we determine $L$ and $U$ as well as the other six unknown constants $l_1, l_2, \lambda_1, \lambda_2, \tilde \lambda_1, \tilde \lambda_2$. We do so by matching the value and the first-order derivative of $Y$ at the points
$U$, $0$ and $L$. This leads to eight nonlinear equations in total as illustrated below.
In addition, using such a construction, the function $Y$ will be twice continuously
differentiable with the exception of at most three points.
Let us first consider such matchings at the point $U$.
From \eqref{eq:derofY}, \eqref{eq:caseI:YlargerU} and \eqref{eq:caseI:YLU}, we obtain three equations:
\begin{eqnarray*}
Y'(U+) & = & \frac{\mathcal{C}_+}{\alpha}+l_1 t_2 e^{t_2 U}=\mathcal{D}_p, \\
Y'(U-) & = & \frac{\mathcal{C}_+}{\alpha}+\lambda_1 r_1 e^{r_1 U}+\lambda_2 r_2 e^{r_2 U}=\mathcal{D}_p,\\
Y(U+) & = & \frac{\mathcal{C}_+}{\alpha} U + \frac{1}{\alpha} \left(\frac{\mathcal{C}_+}{\alpha}(\theta_\ell-b) - \mathcal{D}_p \theta_\ell\right)+ l_1 e^{t_2U} \\
&=&Y(U-)  =  \frac{\mathcal{C}_+}{\alpha} U + \frac{-b \mathcal{C}_+}{\alpha^2}+ \lambda_1 e^{r_1 U} + \lambda_2 e^{r_2 U}.
\end{eqnarray*}
Upon simplifying, we can express $l_1, \lambda_1, \lambda_2$ in terms of $U$ and the model parameters as follows:
\begin{eqnarray}
l_1 & = & \frac{1}{t_2} (\mathcal{D}_p - \frac{\mathcal{C}_+}{\alpha}) e^{-t_2 U}, \label{eq:caseI:l1} \\
\lambda_1  &=& \frac{\sigma^2}{2 \alpha^2} e^{-r_1 U} \left(1-\frac{r_1}{r_2}\right)^{-1} (-B_1), \label{eq:caseI:lambda1}\\
\lambda_2  &=& \frac{\sigma^2}{2 \alpha^2} e^{-r_2 U} \left(1-\frac{r_2}{r_1}\right)^{-1} J_1. \label{eq:caseI:lambda2}
\end{eqnarray}
Similarly, matching at the point $L$, we obtain another three equations:
\begin{eqnarray*}
Y'(L+) & = & \frac{\mathcal{C}_+}{\alpha}+  \lambda_1 r_1 e^{r_1 L}+ \lambda_2 r_2 e^{r_2 L}=-\mathcal{I}_p, \\
Y'(L-) & = & \frac{\mathcal{C}_+}{\alpha}  +  \tilde \lambda_1 s_1 e^{s_1 L} + \tilde \lambda_2 s_2 e^{s_2 L}=-\mathcal{I}_p,\\
Y(L-) &= & \frac{\mathcal{C}_+}{\alpha} L + \frac{1}{\alpha} \left(\frac{\mathcal{C}_+}{\alpha}(\theta_u -b) + \mathcal{I}_p \theta_u\right)+ \tilde \lambda_1  e^{s_1 L} + \tilde \lambda_2  e^{s_2 L} \\
&=&Y(L+ ) = \frac{\mathcal{C}_+}{\alpha} L + \frac{- b \mathcal{C}_+}{\alpha^2}+ \lambda_1 e^{r_1 L} +  \lambda_2 e^{r_2 L}.
\end{eqnarray*}
Upon plugging the expressions for $\lambda_1$ and $\lambda_2$ into $Y'(L+)=-\mathcal{I}_p$,
we obtain \eqref{eq:xhminusxf}.
On the other hand, it follows from $Y(L+)=Y(L-)$ that
\begin{eqnarray} \label{eq:caseI:tildelambda1}
\lambda_1 e^{r_1 L} +  \lambda_2 e^{r_2 L} \; = \; \frac{\theta_u}{\alpha} \left(\mathcal{I}_p+ \frac{\mathcal{C}_+}{\alpha}\right) + \tilde \lambda_1  e^{s_1 L} + \tilde \lambda_2  e^{s_2 L} .
\end{eqnarray}
Combining this with $Y'(L-)=-\mathcal{I}_p$, we cancel $\tilde \lambda_1$ and derive
\begin{equation}\label{eq:L3}
\lambda_1 e^{r_1 L} +  \lambda_2 e^{r_2 L} \; = \; \left( \frac{\theta_u}{\alpha}-  \frac{1}{s_1}\right) \left(\mathcal{I}_p+ \frac{\mathcal{C}_+}{\alpha}\right) +
 \tilde \lambda_2 e^{s_2 L} \left(1- \frac{s_2}{s_1}\right) .
\end{equation}
Lastly, solving $\tilde \lambda_2$ by matching $Y$ at the point 0 to order one, we obtain two equations:
\begin{eqnarray*}
Y'(0+)\;=\;\frac{\mathcal{C}_+}{\alpha}+  \tilde \lambda_1 s_1 + \tilde \lambda_2 s_2&=&Y'(0-)\;=\; -\frac{\mathcal{C}_-}{\alpha}+ l_2 s_1, \\
Y(0-)\;=\;\frac{1}{\alpha} \left( -\frac{\mathcal{C}_-}{\alpha} (\theta_u -b) + \mathcal{I}_p \theta_u\right)  + l_2
&=& Y(0+)\;=\;\frac{1}{\alpha} \left( \frac{\mathcal{C}_+}{\alpha} (\theta_u -b) + \mathcal{I}_p \theta_u\right) + \tilde \lambda_1 + \tilde \lambda_2 .
\end{eqnarray*}
Solving these two equations renders
\begin{eqnarray}
\tilde \lambda_2 \; = \; \frac{\sigma^2}{2 \alpha^2} (\mathcal{C}_+ + \mathcal{C}_-) \frac{{s_1} ^2}{ s_1 -s_2}, &\qquad&
\tilde \lambda_1 - l_2 \; = \; \frac{\sigma^2}{2 \alpha^2} (\mathcal{C}_+ + \mathcal{C}_-)( \frac{ -{s_2} ^2}{ s_1 -s_2}). \label{eq:caseI:l2}
\end{eqnarray}
Upon substituting the expressions for $\lambda_1$, $\lambda_2$ and $\tilde \lambda_2$ into
\eqref{eq:L3} and simplifying the resulting expression, we conclude
that $L$ and $U$ satisfy \eqref{eq:xH}.
Once we determine $L$ and $U$ from \eqref{eq:xhminusxf} and \eqref{eq:xH},
the other unknown constants $l_1, l_2, \lambda_1, \lambda_2, \tilde \lambda_1$ can be derived from \eqref{eq:caseI:l1}-\eqref{eq:caseI:l2} accordingly.

Next, we provide conditions
under which \eqref{eq:xhminusxf} and \eqref{eq:xH} have a unique solution $L$
and $U$ that are both non-negative.
Define for $y>0$,
\begin{eqnarray}
h(y)=B_1 y + J_1 y^{\frac{r_2}{r_1}}+ A .  \label{eq:func-h}
\end{eqnarray}
We first consider the case $\mathcal{I}_p+\mathcal{D}_p>0.$ In this case,
one can readily verify that
\begin{eqnarray} \label{eq:h1}
h(1) = B_1 + J_1 + A = {\alpha} (r_1-r_2) {(-\mathcal{I}_p-\mathcal{D}_p) } <0.
\end{eqnarray}
In addition, $h$ is strictly convex and $\lim_{y \rightarrow 0+} h (y) = \infty$. Thus $h(y)=0$ has a unique solution in $(0,1)$. Since \eqref{eq:xhminusxf} is equivalent to $h(e^{r_1(L-U)}) = 0$, we deduce that $L-U$ is negative and it is uniquely determined by \eqref{eq:xhminusxf}.
On the other hand, observe that $L \ge 0$ holds if and only if $e^{s_2 L} \;\le\; 1$ holds since $s_2<0$. Using \eqref{eq:xH}, this condition is equivalent to
\begin{eqnarray*}
\frac{B_1 r_2}{r_1-r_2} e ^{r_1 (L-U)} + \frac{J_1 r_1}{r_1-r_2} e^{r_2 (L-U)} &\le&
(r_1+r_2-s_1)(\alpha \mathcal{I}_p+{\mathcal{C}_+}) + (\mathcal{C}_+ + \mathcal{C}_-) s_1,
\end{eqnarray*}
which, upon canceling $J_1 e^{r_2 (L-U)}$ using \eqref{eq:xhminusxf}, becomes
\begin{eqnarray*}
 e ^{r_1 (L-U)} \ge \frac{B_3 -B_2}{B_1}.
\end{eqnarray*}
Therefore, guaranteeing that \eqref{eq:xhminusxf} and \eqref{eq:xH} has a unique solution with $U>L \ge 0$ is equivalent to showing that
$h(y)=0$ has a unique solution in the interval
$$\left[\frac{B_3-B_2}{B_1}, 1\right).$$
This is true if either of the following two set of conditions hold:
\begin{eqnarray*}
(i) & \frac{B_3-B_2}{B_1} \in (0,1) \quad \text{and} \quad h\left(\frac{B_3-B_2}{B_1}\right) \ge 0;\\
(ii) &  B_3 \le B_2.
\end{eqnarray*}
Note that $\mathcal{I}_p+\mathcal{D}_p>0.$ We then find that these conditions are exactly condition (1a) and (1b) in Theorem~\ref{THM:CASE1} after applying simple algebraic manipulations. We next consider the case where $\mathcal{I}_p=\mathcal{D}_p=0$. It is clear from \eqref{eq:h1} that $h(1)=0$, and $h'(1)=B_1 +J_1 \frac{r_2}{r_1} = \frac{\mathcal{C}_+(r_1- r_2) t_2}{r_1} <0$. This implies the unique solution to \eqref{eq:xhminusxf} is $L=U := \delta$. Now we deduce from Equation \eqref{eq:xH} that
\[e^{s_2 \delta}= e^{s_2 L} \; = \; \frac{\mathcal{C}_+}{\mathcal{C}_+ +\mathcal{C}_-} \frac{s_1 -t_2}{s_1},\]
and $\delta \ge 0$ if and only if $\mathcal{C}_- s_1 + \mathcal{C}_+ t_2 \ge 0$. This condition is equivalent to $B_1 + B_2 -B_3 \ge 0$, which corresponds to condition (1c) in Theorem~\ref{THM:CASE1}. Note that $\delta=0$ if and only if $B_3-B_1-B_2=0$. Thus, as a byproduct, we have obtained the explicit form of $\delta$ in Corollary~\ref{THM:SPECIAL} when $\delta \ge 0$.

Finally, we verify that the candidate value function $Y$ satisfies the
required first-order properties in \eqref{eq:derofY}.
Since we have constructed the function $Y$ with $Y'(U)=\mathcal{D}_p$ and
$Y'(L)=-\mathcal{I}_p$, then to establish \eqref{eq:derofY} it suffices to
verify the convexity of the function $Y$.
To this end, we first consider $x \ge U$. One readily confirms from \eqref{eq:caseI:YlargerU} and \eqref{eq:caseI:l1} that
$$
Y''(x) \; = \; l_1 t_2^2 e^{t_2 x} \; = \; \left(\mathcal{D}_p- \frac{\mathcal{C}_+}{\alpha}\right) t_2 e^{t_2 (x-U)}.
$$
Given that $\mathcal{D}_p < \mathcal{C}_+ / \alpha$ and that
$t_2<0$ in \eqref{eq:t2}, we can conclude
$$Y''(x) \; > \; 0, \qquad \text{for $x \ge U$.}$$
Similarly, when $x \le 0$, we have
$$Y''(x) \; = \; l_2  s_1^2 e^{s_1 x}.$$
Hence, to show $Y''(x) \ge 0$, it suffices to show that $l_2 \ge 0$.
From \eqref{eq:caseI:l2}, it is equivalent to show
\begin{eqnarray}
\tilde \lambda_1 \ge \frac{\sigma^2}{2 \alpha^2} (\mathcal{C}_+ + \mathcal{C}_-)( \frac{ -{s_2} ^2}{ s_1 -s_2}). \label{ineq:tildelambda1}
\end{eqnarray}
Using the relationship between $\tilde \lambda_1$ and $\tilde \lambda_2$
in $Y'(L-)=-\mathcal{I}_p$ and substituting $\tilde \lambda_2$ given in \eqref{eq:caseI:l2},
we infer from \eqref{ineq:tildelambda1} that we simply need to establish
\begin{eqnarray*}
\frac{\sigma^2}{2 \alpha^2} (\mathcal{C}_+ + \mathcal{C}_-)( \frac{ -{s_2} ^2}{ s_1 -s_2}) s_1 e^{s_1L} +
\frac{\sigma^2}{2 \alpha^2} (\mathcal{C}_+ + \mathcal{C}_-)( \frac{ {s_1} ^2}{ s_1 -s_2}) s_2 e^{s_2L} + (\frac{\mathcal{C}_+}{\alpha} + \mathcal{I}_p) \le 0.
\end{eqnarray*}
Given $\mathcal{I}_p < \frac{\mathcal{C}_-}{\alpha}$, it suffices to show
\begin{eqnarray} \label{eq:caseI:case2xpos3}
s_2 e^{s_1 L} - s_1 e^{s_2 L} + (s_1 -s_2) \le 0.
\end{eqnarray}
This is readily proved since the left-hand side of \eqref{eq:caseI:case2xpos3} is 0 when $L=0$, and it is non-increasing on $[0, \infty)$ as a function of $L.$ Thus we have established the convexity of $Y$ for $x \le 0$.
Turning to the case $L<x<U$, we can verify
$$
Y''(x) \; = \; \lambda_1 r_1^2 e^{r_1x} +\lambda_2 r_2^2 e^{r_2x}.
$$
Upon substituting \eqref{eq:caseI:lambda1} and \eqref{eq:caseI:lambda2}
for $\lambda_1$ and $\lambda_2$, we obtain when $L<x<U$
\begin{eqnarray*}
Y''(x) &=& \left(e^{r_1(x-U)} \frac{r_1^2r_2}{r_1-r_2} B_1 + e^{r_2(x-U)} \frac{r_2^2r_1}{r_1-r_2} J_1 \right) \frac{\sigma^2}{2 \alpha^2} \\
&\ge& \left(\frac{r_1^2r_2}{r_1-r_2} B_1 + \frac{r_2^2r_1}{r_1-r_2} J_1 \right) \frac{\sigma^2}{2 \alpha^2} \\
&=& ( \frac{\mathcal{C}_+} {\alpha} - \mathcal{D}_p) (- t_2) \;\; > \;\; 0.
\end{eqnarray*}
Finally, for the case $0 <x \le L$, we derive
\begin{eqnarray*}
Y''(x)&=& \tilde \lambda_1 s_1^2 e^{s_1x} + \tilde \lambda_2 s_2^2 e^{s_2x} \\
& = & \left(-\mathcal{I}_p - \frac{\mathcal{C}_+}{\alpha} - \tilde \lambda_2 s_2 e^{s_2 L}\right) s_1 e^{s_1 (x-L)}   + \tilde \lambda_2 s_2^2 e^{s_2x}
\end{eqnarray*}
where in the second equality we cancel $\tilde \lambda_1$ by using the relationship between $\tilde \lambda_1$ and $\tilde \lambda_2$
in $Y'(L-)=-\mathcal{I}_p$. After substituting $\tilde \lambda_2$ given in \eqref{eq:caseI:l2}, and simplifying the resulting expression, we obtain
\begin{align} \label{eq:case3Ycov}
Y''(x) = \left[\left(-\mathcal{I}_p - \frac{\mathcal{C}_+}{\alpha}\right) s_1 + \frac{\mathcal{C}_+ + \mathcal{C}_-}{\alpha} s_1 e^{s_2 L} \right] e^{s_1 (x-L)} + \frac{\mathcal{C}_+ + \mathcal{C}_-}{\alpha}  \frac{-s_1 s_2}{s_1 -s_2} e^{s_2 L} (e^{s_2 (x-L)} - e^{s_1 (x-L)}) . \qquad
\end{align}
Suppose we have
\begin{equation} \label{eq:case3cov}
Y''(L-) = \left(-\mathcal{I}_p - \frac{\mathcal{C}_+}{\alpha}\right) s_1 + \frac{\mathcal{C}_+ + \mathcal{C}_-}{\alpha} s_1 e^{s_2 L} \; > \; 0 ,
\end{equation}
then we can show $Y''(x) \; \ge \; 0$ {for all $x \in (0,L]$} by discussing two cases as follows. If
\begin{equation} \label{eq:ineqcase3}
\left(-\mathcal{I}_p - \frac{\mathcal{C}_+}{\alpha}\right) s_1 + \frac{\mathcal{C}_+ + \mathcal{C}_-}{\alpha} s_1 e^{s_2 L} \; \ge \; \frac{\mathcal{C}_+ + \mathcal{C}_-}{\alpha}  \frac{-s_1 s_2}{s_1 -s_2} e^{s_2 L},
\end{equation}
it follows from \eqref{eq:case3Ycov} that $Y''(x) \ge 0$.
On the other hand, if \eqref{eq:ineqcase3} does not hold, one can readily
verify from \eqref{eq:case3Ycov} and the fact $s_1>0 > s_2$ that $Y''(x)$ is non-increasing on $(0, L]$. This also implies that
$Y''(x) > 0$ for $x\in (0,L]$ due to \eqref{eq:case3cov}.
Hence, to show the convexity of the candidate value function $Y$ on
$(0, L]$, it remains to establish \eqref{eq:case3cov}. We rely on the fact that $L$ and $U$ satisfy the Equations \eqref{eq:xhminusxf} and \eqref{eq:xH}.
From \eqref{eq:xH}, we deduce that in order for \eqref{eq:case3cov} to hold,
it suffices to establish that
$${{ \frac{B_1 r_2}{r_1-r_2} e ^{r_1 (L-U)} +  \frac{J_1 r_1}{r_1-r_2} e ^{r_2 (L-U)}}} \; > \; (r_1 +r_2)( \alpha \mathcal{I}_p+ {\mathcal{C}_+}),$$
which becomes
$$B_1 r_2 e ^{r_1 (L-U)} + J_1 r_1 e ^{r_2 (L-U)}  + (r_1 +r_2) A \; > \; 0.$$
Since $L-U$ satisfies \eqref{eq:xhminusxf}, we can equivalently show
$$B_1 r_1 e ^{r_1 (L-U) } + J_1 r_2  e ^{r_2 (L-U)} \; < \; 0 .$$
Write
$$f (y) \; = \; h(e^{r_1 y}),$$
where the function $h$ is given by (\ref{eq:func-h}). We can deduce from the property of $h$ that $f$ is strictly convex, $f(0) <0$, and $f'(0)<0$. This immediately implies that $f$ is strictly decreasing on $(-\infty, 0]$ and
$$f'(L-U) \; = \; B_1 r_1 e^{r_1 (L-U) } + J_1 r_2  e^{r_2 (L-U)} \; < \; 0.$$
Thus we have also established that $Y''$ is nonnegative on $(0, L]$. The proof of convexity of function $Y$ is therefore completed.

%To establish necessary and sufficient conditions under which $L \ge 0$,
%we note that
%\begin{eqnarray*}
%L \ge 0 &\hspace*{-0.1in}\Leftrightarrow&\hspace*{-0.1in} e^{s_2 L} \le 1,\\
%&\hspace*{-0.1in}\Leftrightarrow&\hspace*{-0.1in} \text{RHS \eqref{eq:xH}} \ge (r_1 +r_2 -s_1)(\alpha \mathcal{I}_p +\mathcal{C}_+)+ (\mathcal{C}_+ + \mathcal{C}_- )s_1 ,\\
%&\hspace*{-0.1in}\Leftrightarrow&\hspace*{-0.1in} \text{LHS \eqref{eq:xH}} \ge (r_1 +r_2 -t_2)(\alpha \mathcal{D}_p+\mathcal{C}_-)+ (\mathcal{C}_+ + \mathcal{C}_- )t_2 ,\\
%&\hspace*{-0.1in}\Leftrightarrow&\hspace*{-0.1in} B_2 \frac{r_2}{r_1-r_2} e ^{r_1 (U-L)} +  \frac{r_1}{r_1-r_2}{[-K-B_2 e^{r_1 (U-L)}]} \\
%&\hspace*{-0.1in}&\hspace*{-0.1in}\ge  (r_1 +r_2 -t_2)(\alpha \mathcal{D}_p+\mathcal{C}_-)+ (\mathcal{C}_+ + \mathcal{C}_- )t_2,\\
%&\hspace*{-0.1in}\Leftrightarrow&\hspace*{-0.1in} e ^{r_1 (U-L)} \le \frac{B_3-B_1}{B_2}.
%\end{eqnarray*}
%
%\modify{ADD: more details for two cases: $L \ge 0$ and $L < 0 < U$.}

\subsubsection{Case II: $0 \ge U \ge L$}\label{sec:522}
Similar to Case I, we proceed to solve the Bellman equation \eqref{eq:bellmaneq}
depending on the value of $x$ in relation to $U$, $0$ and $L$.
One readily obtains %For the sake of the clarity of the presentation, we only provide the sketch of proof.
\begin{eqnarray}
Y(x)=
\begin{cases}
\frac{\mathcal{C}_+}{\alpha} x + \frac{1}{\alpha}
\left(\frac{\mathcal{C}_+}{\alpha}(\theta_\ell-b) - \mathcal{D}_p \theta_\ell\right)+ l_3 e^{t_2x},  & x \ge 0,\\
Y(x)=-\frac{\mathcal{C}_-}{\alpha} x + \frac{1}{\alpha}
\left(-\frac{\mathcal{C}_-}{\alpha}(\theta_\ell-b) - \mathcal{D}_p \theta_\ell\right)+ \lambda_3 e^{t_1 x}
+ \lambda_4 e^{t_2 x}, & U < x < 0,\\
Y(x)=- \frac{\mathcal{C}_-}{\alpha} x + \frac{b \mathcal{C}_-}{\alpha^2}+
 \tilde \lambda_3 e^{r_1 x} + \tilde \lambda_4 e^{r_2 x}, & L < x < U,\\
Y(x)=-\frac{\mathcal{C}_-}{\alpha} x + \frac{1}{\alpha}
\left(-\frac{\mathcal{C}_-}{\alpha}(\theta_u-b) +\mathcal{I}_p \theta_u\right)+ l_4
 e^{s_1 x}, & x<L,
\end{cases}
\end{eqnarray}
where $t_i, s_i, r_i$ are given in Section~\ref{sec:main:prelim:Single} and $l_3, l_4, \lambda_3, \lambda_4, \tilde \lambda_3, \tilde \lambda_4$ are constants to be determined.

To determine $L$ and $U$ together with other six unknown constants $l_3, l_4, \lambda_3, \lambda_4, \tilde \lambda_3, \tilde \lambda_4$, we match the value and the first-order derivative of $Y$ at the points
$U$, $0$ and $L$, and thus obtain eight nonlinear equations of those eight unknown numbers as given below.
The first three equations are derived from matching at the point $U$:
\begin{align}
Y'(U+)&=-\frac{\mathcal{C}_-}{\alpha}+\lambda_3 t_1 e^{t_1 U}+\lambda_4 t_2 e^{t_2 U}=\mathcal{D}_p, \nonumber \\
Y'(U-)&=-\frac{\mathcal{C}_-}{\alpha}+\tilde \lambda_3 r_1 e^{r_1 U}+\tilde \lambda_4 r_2 e^{r_2 U}=\mathcal{D}_p, \nonumber\\
Y(U+)&=-\frac{\mathcal{C}_-}{\alpha} U + \frac{1}{\alpha} \left(-\frac{\mathcal{C}_-}{\alpha}(\theta_\ell-b) - \mathcal{D}_p \theta_\ell\right)+ \lambda_3 e^{t_1 U}+\lambda_4 e^{t_2 U} \nonumber \\
&=Y(U-)=-\frac{\mathcal{C}_-}{\alpha} U + \frac{ b \mathcal{C}_-}{\alpha^2}+ \tilde \lambda_3 e^{r_1 U} + \tilde \lambda_4 e^{r_2 U} . \nonumber
\end{align}
We also obtain three equations from matching at $L$:
\begin{align}
Y'(L+) & =  - \frac{\mathcal{C}_-}{\alpha}+ \tilde \lambda_3 r_1 e^{r_1 L}+\tilde \lambda_4 r_2 e^{r_2 L}=-\mathcal{I}_p, \nonumber \\
Y'(L-) & = -\frac{\mathcal{C}_-}{\alpha}  + l_4 s_1 e^{s_1 L}=-\mathcal{I}_p, \nonumber \\
Y(L-) & = -\frac{\mathcal{C}_-}{\alpha} L + \frac{1}{\alpha} \left(-\frac{\mathcal{C}_-}{\alpha}(\theta_u -b) + \mathcal{I}_p \theta_u\right)+ l_4 e^{s_1 L} \nonumber\\
&=Y(L+ ) \; = \; - \frac{\mathcal{C}_-}{\alpha} L + \frac{ b \mathcal{C}_-}{\alpha^2}+ \tilde \lambda_3 e^{r_1 L} + \tilde \lambda_4 e^{r_2 L}. \nonumber
\end{align}
Two additional equations arise from matching at $0$:
\begin{align}
Y'(0+) & \; = \; \frac{\mathcal{C}_+}{\alpha}+ l_3 t_2  =  Y'(0-)= -\frac{\mathcal{C}_-}{\alpha}+ \lambda_3 t_1 +  \lambda_4 t_2 , \nonumber \\
Y(0+) &\; = \; \frac{1}{\alpha} \left( (\theta_\ell -b) \frac{\mathcal{C}_+}{\alpha} -\mathcal{D}_p \theta_\ell\right) + l_3 \nonumber \\
&=Y(0-) \; = \; \frac{1}{\alpha} \left( (\theta_\ell -b) -\frac{\mathcal{C}_-}{\alpha} -\mathcal{D}_p \theta_\ell\right)  +  \lambda_3  +  \lambda_4. \nonumber
\end{align}
Solving these equations, we can conclude that $L$ and $U$
satisfy \eqref{eq:xfminusxh} and \eqref{eq:xF}.

Next, we provide necessary and sufficient conditions
under which \eqref{eq:xfminusxh} and \eqref{eq:xF} have unique solution $L$
and $U$ that are both non-positive. Define for $y>0$,
\begin{eqnarray}
\bar h(y)=B_2 y + J_2 y^{\frac{r_2}{r_1}}+ K. \label{eq:func-hbar}
\end{eqnarray}
We first consider the case $\mathcal{I}_p +\mathcal{D}_p>0$.
Using a similar argument as for Case I, one readily checks that guaranteeing \eqref{eq:xfminusxh} and \eqref{eq:xF} have a unique solution with $ 0 \ge U \ge L$ is equivalent to showing that $\bar h(y)=0$ has a unique solution in the interval $[1, \frac{B_3 -B_1}{B_2}]$.
Since $\bar h$ is strictly increasing on $[1, \infty)$, and $\bar h(1) = B_2+J_2 +K <0$, one deduces that it is equivalent to having the conditions
\begin{eqnarray} \label{ineq:xfnego}
\frac{B_3-B_1}{B_2} > 1 \qquad \text{and} \qquad
\bar h\left(\frac{B_3-B_1}{B_2}\right) \ge 0.
\end{eqnarray}
This can be confirmed by simple algebraic manipulations establishing
that the condition \eqref{ineq:xfnego} are the same as condition (2a) in Theorem~\ref{THM:CASE1}. We next consider the case $\mathcal{I}_p =\mathcal{D}_p=0$. One infers from \eqref{eq:func-hbar} that $\bar h(1) =0$ and thus we obtain $L=U := \delta$. In addition, we derive from \eqref{eq:xF} that
\begin{eqnarray} \label{eq:deltaneg}
 e^{t_1 \delta} = e^{t_1 U} = \frac{\mathcal{C}_-}{\mathcal{C}_++\mathcal{C}_-} \frac{s_1-t_2}{-t_2},
\end{eqnarray}
and $\delta \le 0$ if and only if $\mathcal{C}_- s_1 + \mathcal{C}_+ t_2 \le 0$. This condition is equivalent to $B_3-B_1 - B_2 \ge 0$, which is condition (2b) in Theorem~\ref{THM:CASE1}. One checks that when $\mathcal{I}_p =\mathcal{D}_p=0$, we have $\delta =0$ hold when $B_3 - B_1 - B-2 =0$. In addition, we obtain the explicit form for $\delta \ge 0$ in Corollary~\ref{THM:SPECIAL} from (\ref{eq:deltaneg}).

Finally, we verify the convexity of the function $Y$ which implies that
 $Y$ satisfies the required
first-order properties in \eqref{eq:derofY}. This can be proved by establishing the non-negativity of second-order derivative $Y''$ on each interval
$(-\infty, L), [L, U], (U, 0],$ and $(0, \infty)$ separately. The proof is similar as in Case I, and thus is omitted here.

\subsubsection{Case III: $U \ge 0 \ge L$}\label{sec:523}
%
%Focusing on the case $U \ge 0 \ge L$ where $L$ and $U$ satisfy
%\eqref{eq:xhxf2} and \eqref{eq:xhxf2c}
%subject to \eqref{eq:derofY},
We proceed to solve the Bellman equation \eqref{eq:bellmaneq}
depending on the value of $x$ in relation to $U$, $0$ and $L$.
We obtain
\begin{eqnarray*}
Y(x)=
\begin{cases}
\frac{\mathcal{C}_+}{\alpha} x + \frac{1}{\alpha}
\left(\frac{\mathcal{C}_+}{\alpha}(\theta_\ell-b) - \mathcal{D}_p \theta_\ell\right)+ l_5 e^{t_2x}, & x > U \ge 0,\\
\frac{\mathcal{C}_+}{\alpha} x + \frac{-b \mathcal{C}_+}{\alpha^2}+
 \lambda_5 e^{r_1 x} + \lambda_6 e^{r_2 x}, & 0 \le x<U,\\
- \frac{\mathcal{C}_-}{\alpha} x + \frac{b \mathcal{C}_-}{\alpha^2}+
 \tilde \lambda_5 e^{r_1 x} + \tilde \lambda_6 e^{r_2 x}, & L<x < 0,\\
Y(x)=-\frac{\mathcal{C}_-}{\alpha} x + \frac{1}{\alpha} \left(-\frac{\mathcal{C}_-}{\alpha}(\theta_u-b) +\mathcal{I}_p \theta_u\right)+ l_6 e^{s_1 x}, &x \le L,
\end{cases}
\end{eqnarray*}
where $l_5, l_6, \lambda_5, \lambda_6, \tilde \lambda_5, \tilde \lambda_6$ are unknown constants to be determined. To find these constants,
%\textit{(i).} ~If $$, we have
%$$
%Y'(x) \; \ge \; \mathcal{D}_p \qquad \mbox{and} \qquad \inf_{ \theta_\ell \le \theta \le \theta_u} \cL(\theta, x) \; = \; \cL(\theta_\ell, x).
%$$
%The Bellman equation \eqref{eq:bellmaneq} then yields for $x \ge U$
%\begin{equation} \label{eq:caseIII:YlargerU}
%
%\end{equation}
%where $l_5$ is a generic constant to be determined.
%
%\textit{(ii).} ~If $ 0 \le x \le U$, we have
%$$-\mathcal{I}_p \; < \; Y'(x) \; < \; \mathcal{D}_p, \qquad \mbox{and} \qquad \inf_{ \theta_\ell \le \theta \le \theta_u} \cL(\theta, x) \; = \; 0.
%$$
%Thus we obtain from \eqref{eq:bellmaneq} that for $0 \le x<U$
%\begin{equation} \label{eq:caseIII:Y0U}
%Y(x)=\frac{\mathcal{C}_+}{\alpha} x + \frac{-b \mathcal{C}_+}{\alpha^2}+
% \lambda_5 e^{r_1 x} + \lambda_6 e^{r_2 x},
%\end{equation}
%where $\lambda_5$ and $\lambda_6$ are generic constants to be determined.
%
%\textit{(iii).} ~If $ L<x < 0$, we similarly derive
%\begin{equation} \label{eq:caseIII:YL0}
%Y(x)=- \frac{\mathcal{C}_-}{\alpha} x + \frac{b \mathcal{C}_-}{\alpha^2}+
% \tilde \lambda_5 e^{r_1 x} + \tilde \lambda_6 e^{r_2 x},
%\end{equation}
%where $\tilde \lambda_5$ and $\tilde \lambda_6$ are generic constants to be
%determined.
%
%\textit{(iv).} ~If $ x \le L$, we have
%$$Y'(x) \; \le \; -\mathcal{I}_p$$
%and
%the Bellman equation \eqref{eq:bellmaneq} implies
%\begin{equation} \label{eq:caseIII:YlessL}
%Y(x)=-\frac{\mathcal{C}_-}{\alpha} x + \frac{1}{\alpha} \left(-\frac{\mathcal{C}_-}{\alpha}(\theta_u-b) +\mathcal{I}_p \theta_u\right)+ l_6 e^{s_1 x},
%\end{equation}
%where $l_6$ is a generic constant to be determined.
we match the value and the first-order derivative of $Y$ at the points
$U$, $0$ and $L$. This enables us to establish eight equations for the eight unknowns
$L, U, l_5, l_6, \lambda_5, \lambda_6, \tilde \lambda_5, \tilde \lambda_6$.
Let us first consider such matchings at the point $U$.
%From \eqref{eq:caseIII:YlargerU} and \eqref{eq:caseIII:Y0U}, we obtain
%\begin{eqnarray*}
%Y'(U+)\;=\;\frac{\mathcal{C}_+}{\alpha}+l_5 t_2 e^{t_2 U}=\mathcal{D}_p, &\qquad&
%Y'(U-)\;=\;\frac{\mathcal{C}_+}{\alpha}+\lambda_5 r_1 e^{r_1 U}+\lambda_6 r_2 e^{r_2 U}=\mathcal{D}_p,\\
%Y(U+)\;=\;\frac{\mathcal{C}_+}{\alpha} U + \frac{1}{\alpha} \left(\frac{\mathcal{C}_+}{\alpha}(\theta_\ell-b) - \mathcal{D}_p \theta_\ell\right)+ l_5 e^{t_2U} &=&
%Y(U-)\;=\;\frac{\mathcal{C}_+}{\alpha} U + \frac{-b \mathcal{C}_+}{\alpha^2}+ \lambda_5 e^{r_1 U} + \lambda_6 e^{r_2 U}.
%\end{eqnarray*}
One readily verifies that we can express $\lambda_5$ and $\lambda_6$ in terms of $U$ and the model parameters as follows:
\begin{eqnarray}
\lambda_5 &=& \frac{\sigma^2}{2 \alpha^2} e^{-r_1 U} \left(1-\frac{r_1}{r_2}\right)^{-1} (-B_1), \label{eq:caseIII:lambda1}\\
\lambda_6 &=& \frac{\sigma^2}{2 \alpha^2} e^{-r_2 U} \left(1-\frac{r_2}{r_1}\right)^{-1} J_1. \label{eq:caseIII:lambda2}
\end{eqnarray}
Note that $\lambda_5$ and $\lambda_6$ have the same forms as $\lambda_1$ and $\lambda_2$ given in \eqref{eq:caseI:lambda1} and \eqref{eq:caseI:lambda2} in Case I, but their values could be different since the unknown number $U$ maybe different for Case I and Case III.
Matching at the point $L$, we deduce that
\begin{eqnarray}
\tilde \lambda_5 &=& \frac{\sigma^2}{2 \alpha^2} e^{-r_1 L} \left(1-\frac{r_1}{r_2}\right)^{-1} B_2, \label{eq:caseIII:tlamb1}\\
\tilde \lambda_6 &=& \frac{\sigma^2}{2 \alpha^2} e^{-r_1 L} \left(1-\frac{r_2}{r_1}\right)^{-1} (-J_2). \label{eq:caseIII:tlamb2}
\end{eqnarray}
Matching $Y$
at the point 0 to order one renders
\begin{eqnarray*}
Y'(0+) \; = \; \frac{\mathcal{C}_+}{\alpha}+ \lambda_5 r_1 + \lambda_2 r_2 & = & Y'(0-) \; = \; -\frac{\mathcal{C}_-}{\alpha}+ \tilde \lambda_5 r_1 + \tilde \lambda_6 r_2 ,  \\
Y(0+) \; = \; \frac{-b \mathcal{C}_+}{\alpha^2} + \lambda_5  + \lambda_6 & = & Y(0-) \; = \; \frac{b \mathcal{C}_-}{\alpha^2} + \tilde \lambda_5  + \tilde \lambda_6.
\end{eqnarray*}
Upon substituting the equations
\eqref{eq:caseIII:lambda1}~--~\eqref{eq:caseIII:tlamb2} for
$\lambda_5$, $\lambda_6$, $\tilde \lambda_5$, $\tilde \lambda_6$
and simplifying the resulting expressions, we conclude that $L$ and $U$
satisfy the two equations given in \eqref{eq:xhxf2} and \eqref{eq:xhxf2c}.

Next, we provide conditions
under which \eqref{eq:xhxf2} and \eqref{eq:xhxf2c} have a unique solution $L$
and $U$ with $L \le 0 \le U$. To this end, we define for $y>0$,
\begin{eqnarray}
\tilde h(y)=J_1 \left(\frac{B_3 -B_2 y}{B_1} \right)^{\frac{r_2}{r_1}} + J_2 y^{\frac{r_2}{r_1}} -J_3. \label{eq:func-htilde}
\end{eqnarray}
One checks that in order for the existence of a unique pair $L$
and $U$ with $L \le 0 \le U$ solving \eqref{eq:xhxf2} and \eqref{eq:xhxf2c}, it is equivalent to
$\tilde h (y) = 0$ having a unique solution on the interval $\Big[\max\{1, \frac{B_3 -B_1}{B_2} \}, \frac{B_3}{B_2} \Big)$.
This is true if and only if
\begin{eqnarray} \label{eq:tildeh-ineq}
\tilde h \Big(\max\{1, \frac{B_3 -B_1}{B_2} \} \Big) \le 0,
\end{eqnarray}
after noting that $\lim_{y \rightarrow \frac{B_3}{B_2}-} \tilde h (y) = \infty$ and $\tilde h$ is strictly convex. Simplifying \eqref{eq:tildeh-ineq} we arrive at the following two conditions:
\begin{eqnarray*}
(i) &  B_3 > B_2, B_3 -B_1 -B_2 <0, \quad \mbox{and} \quad \left(\frac{B_3 -B_2}{B_1} \right)^{\frac{r_2}{r_1}} \le \frac{J_3 -J_2}{J_1}, \\
(ii) & B_3 > B_2, B_3 -B_1 -B_2 \ge 0, \quad \mbox{and} \quad \left(\frac{B_3 -B_1}{B_2} \right)^{\frac{r_2}{r_1}} \le \frac{J_3 -J_1}{J_2}.
\end{eqnarray*}
These conditions contain the complement of the union of Conditions (1a)-(1c) and Conditions (2a)-(2b) in Theorem~\ref{THM:CASE1}. %In the special case where $\mathcal{I}_p= \mathcal{D}_p = 0$, we have $B_1+B_2 =B_3$ and $J_1+J_2 =J_3$. Thus we infer from \eqref{eq:xhxf2} and \eqref{eq:xhxf2c} that $L=U := \delta = 0$.

Finally, we verify that the candidate value function $Y$ satisfies the
required first-order properties in \eqref{eq:derofY} by establishing
the convexity of the function $Y$.
Note for $x \ge U$, one first readily confirms that
$$
Y''(x) \;\; = \;\; \left(\mathcal{D}_p- \frac{\mathcal{C}_+}{\alpha}\right) t_2 e^{t_2 (x-U)}> 0,
$$
since $\mathcal{D}_p < \mathcal{C}_+ / \alpha$ and $t_2<0$ in \eqref{eq:t2}.
Similarly, when $x \le L, $ we have
$$
Y''(x) \;\; = \;\; \left(\frac{\mathcal{C}_-}{\alpha}-\mathcal{I}_p\right) s_1 e^{s_1 (x-L)}>0,
$$
due to the facts that $\mathcal{I}_p< \mathcal{C}_- / \alpha$ and
$s_1>0$ in \eqref{eq:s1}.
Turning to the case $0<x<U$, one can readily verify that
$$
Y''(x) \;\; = \;\; \lambda_5 r_1^2 e^{r_1x} +\lambda_6 r_2^2 e^{r_2x}.
$$
Upon substituting \eqref{eq:caseIII:lambda1} and \eqref{eq:caseIII:lambda2}
for $\lambda_5$ and $\lambda_6$, we obtain when $0<x<U$
\begin{eqnarray*}
Y''(x) &=& \left(e^{r_1(x-U)} \frac{r_1^2r_2}{r_1-r_2} B_1 + e^{r_2(x-U)} \frac{r_2^2r_1}{r_1-r_2} J_1 \right) \frac{\sigma^2}{2 \alpha^2} \\
&\ge& \left(\frac{r_1^2r_2}{r_1-r_2} B_1 + \frac{r_2^2r_1}{r_1-r_2} J_1 \right) \frac{\sigma^2}{2 \alpha^2}\\
&=& ( \mathcal{D}_p - \frac{\mathcal{C}_+} {\alpha})  t_2 \;\; > \;\; 0.
\end{eqnarray*}
Similarly, for $L<x \le 0$, we have
\begin{eqnarray*}
Y''(x) \;\; = \;\; \tilde \lambda_5 r_1^2 e^{r_1x} + \tilde \lambda_6 r_2^2 e^{r_2x}
\;\;=\;\;(\frac{\mathcal{C}_-}{\alpha} -  \mathcal{I}_p)  s_1 \;\; > \;\; 0.
\end{eqnarray*}
Hence, the convexity of the candidate value function $Y$ has been established.

\subsection{Proof of Theorem~\ref{THM:CASE1} and Corollary~\ref{THM:SPECIAL}: Step 3}\label{sec:proofs:case1:step3}
The final step of our proof of Theorem~\ref{THM:CASE1} consists of
verifying that the proposed two-threshold dynamic control policy is
optimal and that $Y(x)=V(x)$ for all $x$.
We take a martingale argument approach where the key idea is to construct
a submartingale to prove that the candidate value function $Y$ is a lower
bound for the stochastic optimization problem \eqref{optHF:obj}~--~\eqref{optHF:st3}.
To this end, first consider an admissible process $\{P(t): t \ge 0\}$ adapted to the filtration $\mathcal{F}_t$ generated by
$\{D(s): 0 \le s \le t\}$ and $dP(t)=\theta dt$,
where $\theta \in [\theta_\ell, \theta_u]$.
Recalling $X(t)=P(t)-D(t)$, we define for $t \ge 0$
\begin{align*}
M(t) &:= e^{-\alpha t} Y(X(t)) + \int_{0}^{t}e^{-\alpha s} (\mathcal{C}_+ X(s)^+
+ \mathcal{C}_- X(s)^- + \mathcal{I}_p \theta \indi{\theta \ge 0} - \mathcal{D}_p \theta \indi{\theta<0}) ds,
\end{align*}
with our goal being to show that $\{M(t): t \ge 0 \}$ is a submartingale.

Since $Y$ is twice continuously differentiable with the exception of at most three points, we can apply Ito's formula to $e^{-\alpha t} Y(X(t))$ and obtain
for any $0 \le t_1 \le t_2$
\begin{align} \label{eq:submat}
M(t_2)-M(t_1) = & \int_{t_1}^{t_2} e^{-\alpha s} \Big(-\alpha Y(X(s)) + \frac{1}{2} \sigma^2 Y''(X(s)) + (\theta-b) Y'(X(s))
%\nonumber \\ & \qquad + \;
+ \mathcal{C}_+ X(s)^+ + \mathcal{C}_- X(s)^- \nonumber \\
& \qquad + \; \mathcal{I}_p \theta \indi{\theta \ge 0} - \mathcal{D}_p \theta \indi{\theta<0} \Big) ds -\int_{t_1}^{t_2} e^{-\alpha s}Y'(X(s)) \sigma dW(s).
\end{align}
We have established in Section~\ref{sec:proofs:case1:step2} that $Y$
satisfies the Bellman equation
$$
-\alpha Y(x) + \frac{1}{2} \sigma^2 Y''(x) -b Y'(x) + \mathcal{C}_+ x^+ + \mathcal{C}_- x^{-} + \inf_{ \theta_\ell \le \theta \le \theta_u} \cL(\theta, x) \; = \; 0,
$$
%\eqref{eq:bellmaneq}
where
\begin{equation*}
\cL(\theta, x) \; = \; \left\{\begin{matrix}
(Y'(x)+\mathcal{I}_p) \theta \qquad \text{if} \quad \theta \ge 0 ,\\
(Y'(x)-\mathcal{D}_p) \theta \qquad \text{if} \quad \theta<0 .
\end{matrix}\right.
\end{equation*}
%together with \eqref{eq:bellmaneq2}.
This implies that, for any given $x$ and any $\theta\in [\theta_\ell, \theta_u]$,
\begin{eqnarray*}
{-\alpha Y(x) + \frac{1}{2} \sigma^2 Y''(x) -b Y'(x) + \mathcal{C}_+ x^+ +
\mathcal{C}_- x^{-}+\cL(\theta, x) } \ge 0.
\end{eqnarray*}
Since $Y'(\cdot)$ is bounded, upon taking the conditional expectation in
\eqref{eq:submat} with respect to the filtration $\mathcal{F}_{t_1}$,
we deduce for any $t_1 \le t_2$ that
\[\E_x [M(t_2)| F_{t_1}]-M(t_1) \ge 0.\]
Namely, $M(t)$ is a submartingale and therefore we have
\begin{eqnarray*}
\E_x [M(t)] \ge M(0)=Y(x), \qquad \text{for any $t \ge 0$.}
\end{eqnarray*}
Letting $t$ go to infinity, we can conclude that $Y$ is a lower bound
for the optimal value of the stochastic optimization problem
\eqref{optHF:obj}~--~\eqref{optHF:st3},
and thus $Y(x)=V(x)$ for all $x$.
Hence, the dynamic control policy characterized by the two threshold
values $L$ and $U$ is indeed optimal.
Our proof of Theorem~\ref{THM:CASE1} and Corollary~\ref{THM:SPECIAL} is complete.

\subsection{Proof of Theorem~\ref{THM:CASE2}}\label{sec:proofs:main:case2}
The results follow along similar lines to the rigorous proof of Theorem~\ref{THM:CASE1}.
We describe some of the key aspects herein.

The Bellman equation \eqref{eq:bellmaneq} still applies and the parameter settings (i.e., ${\mathcal{D}_p \ge \mathcal{C}_+ / \alpha}$ and ${\mathcal{I}_p <  \mathcal{C}_- / \alpha}$)
indicate that we need to construct a function $Y(x)$ and find $L$ such that
\begin{equation} \label{eq:derofY_case_2}
Y'(x)= \left\{\begin{array}{ll}
> -\mathcal{I}_p , & \qquad \text{if} \quad x> L ,\\
 \le -\mathcal{I}_p & \qquad \text{if} \quad x \le L .
\end{array}\right.
\end{equation}
We require $Y(x)=O(|x|)$ as $|x|\rightarrow \infty$, and $Y$ meets
smoothly at the points $L$, $0$ and $U$ to order one.
We consider two cases $L<0$ and $L\ge 0$ separately.

When $L<0$, solving the Bellman equation for different values of $x$
as in Section~\ref{sec:proofs:case1:step2}, we derive
\begin{equation*}
Y(x) \;\; = \;\; \left\{\begin{array}{ll}
   \frac{\mathcal{C}_+}{\alpha} x + \frac{-b \mathcal{C}_+}{\alpha^2}+
 l_7 e^{r_2 x} & \qquad \text{if} \quad x \ge 0 ,\\
 - \frac{\mathcal{C}_-}{\alpha} x + \frac{b \mathcal{C}_-}{\alpha^2}+
  \lambda_7 e^{r_1 x} +  \lambda_8 e^{r_2 x} & \qquad \text{if} \quad L<x<0 ,\\
-\frac{\mathcal{C}_-}{\alpha} x + \frac{1}{\alpha}
\left(-\frac{\mathcal{C}_-}{\alpha}(\theta_u-b) +\mathcal{I}_p \theta_u\right)+ l_8 e^{s_1 x}, &
\qquad \text{if} \quad  x \le L ,
\end{array}\right.
\end{equation*}
where $l_7, l_8, \lambda_7, \lambda_8$ are
four unknown constants. By matching the value and the first-order derivative of the function $Y$
at the points 0 and $L$, and simplifying the resulting expression, we obtain
\[B_2 e^{-r_1 L} \; = \; B_3,\]
thus rendering \eqref{eq:thm2Lneg} for the case $L<0$. This also implies that $L<0$ if and only if $B_3>B_2$. Moreover, direct calculations establish the convexity of function $Y$, which guarantees \eqref{eq:derofY_case_2}.

When $L \ge 0$, we can proceed in a similar fashion to solve the Bellman equation and find
\begin{eqnarray*}\label{eq:dxK}
(r_2 -s_1)( \alpha \mathcal{I}_p+ {\mathcal{C}_+})+ (\mathcal{C}_+ + \mathcal{C}_- ) s_1 e^{s_2 L}=0.
\end{eqnarray*}
Thus we obtain \eqref{eq:thm2Lpos} for the case $L \ge 0$. In addition, $L$ is nonnegative if and only if
$$(\alpha \mathcal{I}_p+ \mathcal{C}_+) (s_1-r_2 ) \; \le \; (\mathcal{C}_+ +\mathcal{C}_-) s_1,$$
which is equivalent to $B_2 \ge B_3$. In addition, the convexity of the function $Y$ can be readily verified.

Finally, a repetition of the arguments in Section~\ref{sec:proofs:case1:step3}
guarantees that the function $Y$ we construct is indeed the optimal value function. This completes the proof.

%\begin{remark}
%Equation \eqref{eq:dxK} has a close relationship with \eqref{eq:xH}.
%Roughly speaking, the case we consider here can be thought of as $U=\infty$
%and $L \ge 0$, which solves \eqref{eq:xhminusxf} and \eqref{eq:xH} in
%Case I of Theorem~\ref{THM:CASE1}.
%Specifically, if we substitute $L$ (as in \eqref{eq:dxK}) and $U=\infty$
%into \eqref{eq:xhminusxf} and \eqref{eq:xH}, we obtain $B_1 e^{r_1 (L-U)}=0$
%and then $(L, U)$ (with $L$ as in \eqref{eq:dxK} and $U=\infty$) solves the
%equations by pretending $J_1 e^{r_2 (L-U)} + A =0$.
%\end{remark}
%If $L <0$, we have
%\[ e^{-r_1 L} = \frac{B_3}{B_2} = \frac{(\mathcal{C}_+ +\mathcal{C}_-) (-r_2 )}{(\mathcal{C}_- -\alpha \mathcal{I}_p)(s_1 -r_2)} .\]
%If $L \ge 0$, we have
%\[e^{s_2 L}= \frac{(\alpha \mathcal{I}_p +\mathcal{C}_+)(s_1 -r_2)}{ (\mathcal{C}_+ +\mathcal{C}_-) s_1}.\]

\subsection{Proof of Theorem~\ref{THM:CASE3}}\label{sec:proofs:main:case3}
The results follow along similar lines to the rigorous proof of Theorem~\ref{THM:CASE1}.
We describe some of the key aspects herein.

Theorem~\ref{THM:CASE3} is the antithesis of Theorem~\ref{THM:CASE2}.
Our task is to construct a function $Y(x)$ that satisfies the Bellman equation
\eqref{eq:bellmaneq}
together with a number $U$ such that
\begin{equation*} \label{eq:derofY_case3}
Y'(x)= \left\{\begin{array}{ll}
 \ge \mathcal{D}_p  , & \qquad \text{if} \quad x \ge U ,\\
< \mathcal{D}_p &  \qquad \text{if} \quad x <U .
\end{array}\right.
\end{equation*}
We require $Y(x)=O(|x|)$ as $|x|\rightarrow \infty$, and $Y$ meets
smoothly at the points $L$, $0$ and $U$ to order one.
We again consider two cases $U \ge 0$ and $U<0$.

When $U \ge 0$, solving the Bellman equation for different values of $x$
as in Section~\ref{sec:proofs:case1:step2}, we derive
\begin{equation*}
Y(x) \;\;=\;\; \left\{\begin{array}{ll}\frac{\mathcal{C}_+}{\alpha} x + \frac{1}{\alpha}
\left(\frac{\mathcal{C}_+}{\alpha}(\theta_\ell-b) - \mathcal{D}_p \theta_\ell\right)+ l_9 e^{t_2x},  & \qquad \text{if} \quad x \ge U ,\\\frac{\mathcal{C}_+}{\alpha} x + \frac{-b \mathcal{C}_+}{\alpha^2}+
 \lambda_9 e^{r_1 x} + \lambda_{10} e^{r_2 x}, & \qquad \text{if} \quad U> x \ge 0 , \\- \frac{\mathcal{C}_-}{\alpha} x + \frac{b \mathcal{C}_-}{\alpha^2}+
 l_{10} e^{r_1 x}, & \qquad \text{if} \quad x<0 ,
\end{array}\right.
\end{equation*}
where $l_9, l_{10}, \lambda_9,\lambda_{10}$ are
four constants to be determined.
By matching the value and derivative of the function $Y$ at the points
0 and $U$, we get four equations and further obtain
\[J_1 e^{-r_2 U} \; = \; J_3.\]
This implies that $U \ge 0$ if and only if $J_3 \ge J_1$. Hence we have derived \eqref{eq:thm3U} for $U \ge 0$.

When $U<0$, proceeding in a similar fashion, we deduce
$$
(r_1 -t_2)(\alpha \mathcal{D}_p+\mathcal{C}_-)+ (\mathcal{C}_+ + \mathcal{C}_- )t_2 e^{t_1 U} \; = \; 0,
$$
%\begin{remark}
%This essentially solves \eqref{eq:xfminusxh} and
%\eqref{eq:xF} in nondegenerate Case II by the above $U$ and setting
%$L=-\infty$.
%\end{remark}
from which we obtain
\[ e^{t_1 U} \; = \; \frac{(\alpha \mathcal{D}_p + \mathcal{C}_-)(r_1-t_2)}{(\mathcal{C}_+ + \mathcal{C}_-) (-t_2)}, \]
thus rendering \eqref{eq:thm3U} for $U < 0$.

The convexity of the candidate value function $Y(x)$ in both cases can be directly verified and the function $Y$ can be shown to be the optimal value function as in the proof of Section~\ref{sec:proofs:case1:step3}. Therefore we have completed the
proof of Theorem~\ref{THM:CASE3}.

\subsection{Proof of Theorem~\ref{THM:CASE4}}\label{sec:proofs:main:case4}
The results follow along similar lines to the rigorous proof of Theorem~\ref{THM:CASE1}.
Some of the key aspects are described herein.

We construct a function $Y(x)$ with $-\mathcal{I}_p < Y'(x) \; < \; \mathcal{D}_p$ for any $x$.
The Bellman equation \eqref{eq:bellmaneq} together with the linear growth of $Y$ then implies
\begin{equation*}
Y(x) \;\;=\;\; \left\{\begin{array}{ll}
\frac{\mathcal{C}_+}{\alpha} x + \frac{- b \mathcal{C}_+}{\alpha^2}+
 l_{11} e^{r_2 x}
 & \qquad \text{if} \quad x \ge 0 , \\ \frac{-\mathcal{C}_-}{\alpha} x + \frac{b \mathcal{C}_-}{\alpha^2}+
 l_{12} e^{r_1 x}, & \qquad \text{if} \quad x<0 .
\end{array}\right.
\end{equation*}
Matching the value and
the first-order derivative of the function $Y$ at the point 0, we obtain
\begin{eqnarray*}
l_{11} \;=\; \frac{\sigma^2 (\mathcal{C}_+ +\mathcal{C}_-)}{2\alpha^2} \frac{r_1^2}{(r_1 -r_2)} >0
& \qquad \mbox{and} \qquad &
l_{12} \;=\; \frac{\sigma^2 (\mathcal{C}_+ +\mathcal{C}_-)}{2\alpha^2} \frac{r_2^2}{(r_1 -r_2)}>0.
\end{eqnarray*}
It is clear $Y$ is convex and $-\mathcal{I}_p \; \le \; -\frac{\mathcal{C}_-}{\alpha}< Y'(x) \; < \; \frac{\mathcal{C}_+}{\alpha} \; \le \; \mathcal{D}_p$ \text{for any $x$}.
Therefore, the proof is complete after applying the argument in~Section~\ref{sec:proofs:case1:step3}.

\section{Proofs of Results in Section~\ref{sec:main2}} \label{sec:proofs2}
In this appendix we provide proofs of our main results in Section~\ref{sec:main2}.
We start with Theorem~\ref{THM:CASE1:Multiple}, leveraging the proofs of Appendix~\ref{sec:proofs}.
Then we present the proof of Theorem~\ref{thm:contract}.

\subsection{Proof of Theorem~\ref{THM:CASE1:Multiple}}
The results follow from the arguments establishing Theorem~\ref{THM:CASE1} applied to the aggregate primary resource capacity $P$
according to the given contract vector $w$.

\subsection{Proof of Theorem~\ref{thm:contract}}
%
%This subsection is devoted to proving Theorem~\ref{thm:contract}.
%The key idea of the proof is to show that the optimal threshold values $L(w)$ and $U(w)$ are continuous functions of $w$. %, which leads to the continuity of the value function $\mathcal{V}_w(x)$ with respect to $w$.
%Fix $x \in \R$. Since $\Omega$ defined in (\ref{eq:omega}) is a compact set in $\R^{\mathcal{P}}$, we deduce that in order to prove the existence of the optimal contract $w^{*}(x)$, it suffices to show that $\mathcal{V}_w(x)$ is continuous with respect to $w \in \Omega$. Write ${{V}}_w(x)$ to represent the value function of the stochastic dynamic program \eqref{optHF:obj}-\eqref{optHF:st3} with parameters given in \eqref{eq:Rpw}-\eqref{eq:cplusminusw}. We then infer from \eqref{eq:simplify} that
%\begin{eqnarray*}
%\mathcal{V}_w(x) = (\mathcal{R}_p (w) - \mathcal{C}_p (w)) \cdot \E \left[ \int_0^\infty e^{-\alpha t} D(t) dt \right] - V_{w}(x),
%\end{eqnarray*}
% where $\mathcal{R}_p (w), \mathcal{C}_p (w)$ given in \eqref{eq:Rpw} are linear in $w$. Thus it is enough to show $V_w(x)$ is continuous with respect to $w \in \Omega$.
%
We first prove that the optimal threshold values $L(w)$ and $U(w)$, as functions of $w$, are continuous. We start with the first case in Theorem~\ref{THM:CASE1:Multiple} where $U(w) \ge L(w) \ge 0$.
Since $U(w)$ and $L(w)$ are uniquely determined by \eqref{eq:xhminusxf_w} and \eqref{eq:xH_w}, it suffices to show that the solutions to these two equations depend continuously on $w$. %If $\mathcal{D}_p (w) = \mathcal{I}_p (w)=0$, we have $L(w)=U(w)$ and \eqref{eq:xH} immediately implies the continuity of these two functions with respect to $w$. Now assume $\mathcal{D}_p (w) + \mathcal{I}_p (w)>0$.
We first derive the continuity of $L(w)-U(w)$ from \eqref{eq:xhminusxf_w} by considering the following optimization problem parameterized by $w\in \Omega$:
\begin{eqnarray}\label{eq:func-hw}
\min_{y} |h(y, w)|=|B_1(w) y + J_1(w) y^{\frac{r_2}{r_1}}+ A(w)| \quad \text{subject to $y \in [0,1]$},
\end{eqnarray}
where we recall
\begin{eqnarray*}
B_1(w) &:=& (\mathcal{C}_+(w)-\alpha \mathcal{D}_p(w)) (\tilde t_2 -r_2), \qquad
J_1(w):= (\mathcal{C}_+(w)-\alpha \mathcal{D}_p(w)) (r_1 - \tilde t_2), \\
A(w) &:= &(\mathcal{C}_+(w) + \alpha \mathcal{I}_p(w)) (r_2-r_1).
\end{eqnarray*}
For fixed $w$, one readily verifies that $h(1, w) \le 0$ similarly as in \eqref{eq:h1}. In conjunction with the facts that $h(0, w)= \infty$ and $h$ is strictly convex, we deduce that the optimization problem \eqref{eq:func-hw} has a unique minimizer $y^*(w)$ with $|h(y^*(w), w)| = 0$.

We next argue that $y^*(w)$ is continuous at $w$. We apply Proposition 4.4 of Bonnans and Shapiro~\cite{bonnans2000perturbation} and check the four conditions there. Note that the feasible set in the optimization problem \eqref{eq:func-hw} is constant $[0,1]$ (independent of $w$) and this feasible set is closed. It readily follows that the conditions (ii) and (iv) in Proposition 4.4 of \cite{bonnans2000perturbation} hold. Moreover, it is clear that the function $|h(y, w)|$ is jointly continuous in its two arguments, so condition (i) holds. Finally, the level set $\{y \in[0,1]: |h(y, w)| \le c\}$ is nonempty for each $c \ge 0$ and it is contained in the compact set $[0,1]$, so condition (iii) also holds.
 Thus we deduce from Proposition 4.4 of \cite{bonnans2000perturbation} that the optimal solution $y^*(w)$ is continuous with respect to $w$. This implies the solution $L(w)-U(w)$ to \eqref{eq:xhminusxf_w} is continuous in $w$. Using this result, one infers from \eqref{eq:xH_w} that $L(w)$ also depends continuously on $w$, thus completing the proof of continuity of $L(w)$ and $U(w)$ in the first case of Theorem~\ref{THM:CASE1:Multiple}. For the second and third case of Theorem~\ref{THM:CASE1:Multiple}, we can apply similar arguments and conclude that $L(w)$ and $U(w)$ are both continuous at each $w \in \Omega$.

We are now ready to prove the value function $V_w(x)$ is continuous with respect to $w \in \Omega$. To see this, we again illustrate it using the first case of Theorem~\ref{THM:CASE1:Multiple}. The other two cases are similar. Since $L(w)$ and $U(w)$ depend continuously on $w$, one can verify from \eqref{eq:caseI:l1}-\eqref{eq:caseI:l2} that the same continuity property also applies to
$l_1, l_2, \lambda_1, \lambda_2, \tilde \lambda_1, \tilde \lambda_2$ (these constants now depend on $w$) appeared in Section~\ref{sec:521}. Thus we conclude from the explicit form of value function \eqref{eq:caseI:YlargerU}-\eqref{eq:caseI:YlessL} that $V_w(x)$ is a continuous function of $w$.

Finally, the continuity of $J_w(x)$ with respect to $w$ readily follows from the continuity of $V_w(x)$ and Equation \eqref{eq:J-V}. The proof is therefore complete.

%\bibliographystyle{ormsv080}
%\bibliography{gao}

%\newpage
\vspace*{1.0in}
\bibliographystyle{plain}
%\bibliography{gao,main}

\end{document}